\documentclass[a4paper,12pt,reqno]{amsart}
\usepackage[T1]{fontenc}
\usepackage[utf8]{inputenc}
\usepackage[english]{babel}
\usepackage{amsmath, amsthm, amssymb}
\usepackage{siunitx}
\usepackage{fullpage}
\usepackage{xcolor}
\usepackage{enumerate}
\usepackage{caption}
\usepackage{subcaption}
\usepackage{chemformula}
\usepackage{pgfplotstable}
\usetikzlibrary{
	matrix,
}
\pgfplotsset{
compat = newest,
tick label style = {font = \tiny},
legend style = {font = \tiny},
xlabel style={yshift=+0.5ex},
ylabel style={yshift=-1.0ex}
}
\usepackage[width=0.95\textwidth]{caption}
\usepackage[width=0.95\textwidth,labelfont=rm]{subcaption}
\usepackage{doi, hyperref}
\makeatletter
\def\@seccntformat#1{%
  \protect\textup{\protect\@secnumfont
    \ifnum\pdfstrcmp{subsection}{#1}=0 \bfseries\fi
    \csname the#1\endcsname
    \protect\@secnumpunct
  }%
}
\makeatother
\newcommand*\patchAmsMathEnvironmentForLineno[1]{%
  \expandafter\let\csname old#1\expandafter\endcsname\csname #1\endcsname
  \expandafter\let\csname oldend#1\expandafter\endcsname\csname end#1\endcsname
  \renewenvironment{#1}%
     {\linenomath\csname old#1\endcsname}%
     {\csname oldend#1\endcsname\endlinenomath}}%
\newcommand*\patchBothAmsMathEnvironmentsForLineno[1]{%
  \patchAmsMathEnvironmentForLineno{#1}%
  \patchAmsMathEnvironmentForLineno{#1*}}%
\AtBeginDocument{%
\patchBothAmsMathEnvironmentsForLineno{equation}%
\patchBothAmsMathEnvironmentsForLineno{align}%
\patchBothAmsMathEnvironmentsForLineno{flalign}%
\patchBothAmsMathEnvironmentsForLineno{alignat}%
\patchBothAmsMathEnvironmentsForLineno{gather}%
\patchBothAmsMathEnvironmentsForLineno{multline}%
}
\usepackage[mathlines]{lineno}

\theoremstyle{definition}
\newtheorem{theorem}{Theorem}[section]
\newtheorem{proposition}[theorem]{Proposition}
\newtheorem{lemma}[theorem]{Lemma}
\newtheorem{algorithm}[theorem]{Algorithm}
\newtheorem{definition}[theorem]{Definition}
\newtheorem{remark}[theorem]{Remark}
\newcommand{\comment}[1]{{\color{black}#1}}
\renewcommand\AA{\boldsymbol{A}}
\renewcommand\gg{\boldsymbol{g}}
\newcommand\ee{\boldsymbol{e}}
\newcommand\hh{\boldsymbol{h}}
\newcommand\mm{\boldsymbol{m}}
\newcommand\nn{\boldsymbol{n}}
\newcommand\uu{\boldsymbol{u}}
\newcommand\UU{\boldsymbol{U}}
\newcommand\vv{\boldsymbol{v}}
\newcommand\ww{\boldsymbol{w}}
\newcommand\CC{\boldsymbol{C}}
\newcommand\HH{\boldsymbol{H}}
\newcommand\LL{\boldsymbol{L}}
\newcommand\ff{\boldsymbol{f}}
\newcommand\E{\mathcal{E}}
\newcommand\FF{\boldsymbol{F}}
\newcommand\MM{\boldsymbol{\mathcal{M}}}
\newcommand\KK{\boldsymbol{\mathcal{K}}}
\newcommand\II{\mathcal{I}}
\newcommand\NN{\mathcal{N}}

\renewcommand\S{\mathcal{S}}
\newcommand\T{\mathcal{T}}
\newcommand{\eps}{\varepsilon}
\newcommand{\boldvar}{\boldsymbol{\varepsilon}}
\newcommand{\boldvarm}{\boldvar_{\mathrm{m}}}
\newcommand{\boldvarel}{\boldvar_{\mathrm{el}}}
\newcommand{\boldsig}{\boldsymbol{\sigma}}
\newcommand\pphi{\boldsymbol{\phi}}
\newcommand\ppsi{\boldsymbol{\psi}}
\newcommand\vvphi{\boldsymbol{\varphi}}
\newcommand\xxi{\boldsymbol{\xi}}
\newcommand\0{\boldsymbol{0}}
\newcommand\sphere{\mathbb{S}^2}
\newcommand\Grad{\boldsymbol{\nabla}}
\newcommand\Lapl{\boldsymbol{\Delta}}

\newcommand\A{\mathbb{A}}
\newcommand\C{\mathbb{C}}
\newcommand\R{\mathbb{R}}
\newcommand\dt{\mathrm{d}_t}

\newcommand\Z{\mathbb{Z}}
\newcommand{\abs}[1]{\left\lvert #1 \right\rvert}

\newcommand{\inner}[3][]{\langle #2,#3 \rangle_{#1}}
\newcommand{\norm}[2][]{\left\lVert #2 \right\rVert_{#1}}
\DeclareMathOperator{\diam}{diam}
\DeclareMathOperator{\trace}{tr}
\newcommand\ellex{\ell_{\textnormal{ex}}}
\newcommand\heff{\hh_{\mathrm{eff}}}
\newcommand\hhm{\hh_{\mathrm{m}}}
\newcommand\mmt{\partial_t \mm}
\newcommand\uut{\partial_{t}\uu}
\newcommand\uutt{\partial_{tt}\uu}
\newcommand{\weakstarto}{\overset{\ast}{\rightharpoonup}}
\newcommand{\weakto}{\rightharpoonup}

\newenvironment{customlegend}[1][]{%
    \begingroup
    \csname pgfplots@init@cleared@structures\endcsname
    \pgfplotsset{#1}%
}{%
    \csname pgfplots@createlegend\endcsname
    \endgroup
}%

\def\addlegendimage{\csname pgfplots@addlegendimage\endcsname}

\begin{document}
\title{A decoupled, convergent and fully linear algorithm
for the Landau--Lifshitz--Gilbert equation\\ with magnetoelastic effects}
\author{Hywel~Normington}
\author{Michele~Ruggeri}
\address{Department of Mathematics and Statistics,
University of Strathclyde,
26 Richmond Street, Glasgow G1~1XH, United Kingdom}
\email{hywel.normington@strath.ac.uk}
\email{michele.ruggeri@strath.ac.uk}
\date{\today}

\begin{abstract}
We consider the coupled system of the Landau--Lifshitz--Gilbert equation
and the conservation of linear momentum law
to describe magnetic processes in ferromagnetic materials including magnetoelastic effects
in the small-strain regime.
For this nonlinear system of time-dependent partial differential equations, we present a decoupled integrator
based on first-order finite elements in space and an implicit one-step method in time.
We prove unconditional convergence of the sequence of discrete approximations towards a weak solution of the system
as the mesh size and the time-step size go to zero.
Compared to previous numerical works on this problem, for our method,
we prove a discrete energy law that mimics that of the continuous problem and,
passing to the limit, yields an energy inequality satisfied by weak solutions.
Moreover, our method does not employ a nodal projection to impose the unit length constraint on the discrete magnetisation,
so that the stability of the method does not require weakly acute meshes.
Furthermore, our integrator and its analysis hold for a more general setting, including body forces and traction,
as well as a more general representation of the magnetostrain.
Numerical experiments underpin the theory and showcase the applicability of the scheme
for the simulation of the dynamical processes involving magnetoelastic materials at submicrometer length scales.
\end{abstract}

\maketitle


\section{Introduction}

Magnetoelastic (or magnetostrictive) materials
are smart materials characterised by a strong interplay
between their mechanical and magnetic properties~\cite{brown1966}.
On the one hand, they change shape when subject to applied magnetic fields (direct magnetostrictive effect), and on the other,
they undergo a change in their magnetic state
when subject to externally applied mechanical stresses
(inverse magnetostrictive effect).
Because of these properties,
magnetoelastic materials currently find use in many technological applications requiring a magnetomechanical transducer, e.g.\ actuators or sensors~\cite{pasquale2003}.

In this work,
we design and analyse a fully discrete numerical scheme
for a coupled nonlinear system of partial differential equations
(PDEs) modelling the dynamics of magnetisation
and displacement in magnetoelastic materials in the small-strain regime.
\comment{The small-strain assumption is
well justified for many ferromagnetic materials,
a variety of which experience strains on the order of $10^{-5}$;
see Table~\ref{tab:valuetable} below for an example of such a material,
and Figure~\ref{fig:illustmagnetostriction} for an illustration of a ferromagnetic cube with those material parameters.

\begin{figure}
	\centering
	\includegraphics[height=2.5in]{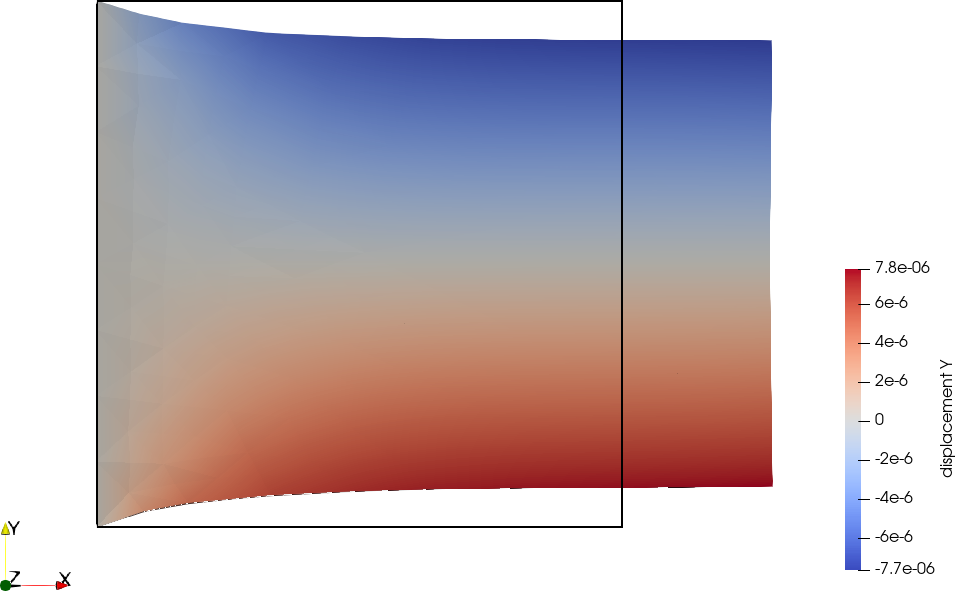}
	\caption{\comment{Illustration of isotropic, isochoric magnetostriction of a cube shaped ferromagnetic material.\label{fig:illustmagnetostriction} The ferromagnetic unit cube is initially uniformly magnetised along the $x$-direction in the undeformed state the edges of which are shown by the black square. The final state shows the body elongated in the $x$-direction, and contracted in the $y$ and $z$ directions. The colour indicates the strength of the displacement in the $y$ coordinate.
	The mesh displacement scaled by a factor of $10^4$.}}
\end{figure}

The validity of the small-strain regime
for general magnetostrictive materials is the subject of the seminal work \cite{brown1966}.}
The proposed system consists of the Landau--Lifshitz--Gilbert (LLG) equation
for the magnetisation
and the conservation of linear momentum law for the displacement
(see~\eqref{eq:newton}--\eqref{eq:llg} below).
The two equations are nonlinearly coupled to each other:
One of the contributions to the effective field appearing in the LLG equation depends on the mechanical stress in the body
(and thus on the displacement)
and there is a magnetisation-dependent contribution to the strain (usually referred to as the magnetostrain) 
in the conservation of momentum law.
One additional difficulty is represented
by a nonconvex pointwise constraint on the magnetisation,
which is a vector field of constant unit length.

Several versions of this PDE system 
have been used for physical investigations of magnetoelastic materials; see e.g.\
\cite{shu2004micromagnetic,Mballa2014,bhbvfs2014,pwhcn2015,rj2021,renuka2021solution,dw2023}.
As far as the mathematical literature is concerned,
we refer to~\cite{visintin1985landau,cef2011},
in which existence of weak solutions has been established,
and to a series of works by L.\ Ba\v{n}as and coauthors~\cite{bs2005,banas2005a,bs2006,banas2008,bppr2013}
in numerical analysis.
In~\cite{bs2005,banas2005a,bs2006,banas2008}, the focus is on finite element methods for the approximation of strong solutions.
More recently, \cite{bppr2013} extended the tangent plane scheme proposed in~\cite{alouges2008a} for the LLG equation
to this PDE system.
The integrator, based on first-order finite elements in space and on an implicit first-order time-stepping method in time, decouples the system
and only requires the solution of two linear systems per time-step.
Under the assumption that all meshes used for the spatial discretisation are weakly acute
(needed to guarantee the stability of the nodal projection used to impose the unit length constraint on the magnetisation~\cite{bartels2005,alouges2008a}),
the authors proved unconditional convergence of the finite element approximations towards a weak solution of the problem.

In this work, we generalise the PDE system considered in~\cite{visintin1985landau,cef2011,bppr2013}
by including volume and surfaces forces,
as well as a more general expression for the magnetostrain~\cite{federico2019tensor},
which allows the description of a larger class of magnetoelastic materials.
For this generalised system, we propose an integrator which resembles the one in~\cite{bppr2013}
(same finite element approximation spaces,
same time discretisation method,
same decoupled approach).
However, following~\cite{bartels2016,abert_spin-polarized_2014}
(and differently from~\cite{bppr2013}),
we remove the nodal projection from the update of
the magnetisation (but we keep it in the discretisation of the elastic contributions).
By doing this, we can avoid the requirement
of weakly acute meshes at the expense of not maintaining
the unit length constraint on the magnetisation
at the vertices of the meshes.
However, like in~\cite{bartels2016,abert_spin-polarized_2014},
we can uniformly control the violation of the constraint
by the time-step size.
Despite the strong nonlinearity of the problem,
the resulting integrator is fully linear
(in the sense that it involves only linear operations
like solving linear systems and updating
the approximations using a linear time-stepping).
For this generalised and modified integrator,
we show unconditional well-posedness,
a discrete energy law satisfied by the approximations,
unconditional stability,
and unconditional convergence of the approximations
towards a weak solution of the problem
(here, the adjective `unconditional' refers to the fact that the analysis does not require any restrictive coupling condition between the time and spatial discretisation parameters).
Moreover, assuming a (very restrictive, but artificial) Courant--Friedrichs--Lewy (CFL) condition
on the time-step size and the spatial mesh size,
we can pass the discrete energy law to the limit
and show that the weak solution
towards which our finite element approximation is converging
satisfies an energy inequality.

Summarising, the contribution of the present paper over the existing literature
(and, in particular, over~\cite{visintin1985landau,cef2011,bppr2013}) is threefold:
\begin{itemize}
\item We consider a more general setting than in~\cite{visintin1985landau,cef2011,bppr2013}
\comment{including volume/surface forces and a
more general magnetoelastic contribution based
upon the magnetostriction tensor $\Z$ (see equation~\eqref{eq:magnetostrain} below),
which allows for more general crystal classes to be considered.
In particular, previous papers have described
this tensor as fully symmetric and positive definite, e.g.\ \cite{visintin1985landau,bppr2013,banas2005a,banas2008,cef2011},
which is not true in general (such as in nickel, which has negative magnetostriction).}
Since our convergence proof is constructive, a byproduct of our analysis is a proof of existence of weak solutions for a more general model of magnetoelastic materials in the small-strain regime.
\item Our integrator is energetically `mindful', in the sense that our approximations satisfy a discrete energy law which resembles the one satisfied by solutions of the continuous problem (cf.\ Proposition~\ref{prop:energy} below).
Under a \comment{restrictive CFL condition on the discretisation parameters, specifically $k=o(h^9)$},
we can pass the result to the limit
and obtain an energy inequality for weak solutions.
This aspect was not considered in~\cite{bppr2013}, where only boundedness of energy was proven.
\item The spatial meshes used by our integrator are assumed to be
only shape-regular (and do not need to be weakly acute as in~\cite{bppr2013}).
This allows for the use of general mesh generators.
\comment{This is especially useful in three dimensions,
as weakly acute meshes are difficult to generate for arbitrary shapes~\cite{brandts2020simplicial}.
Following~\cite{bartels2016,abert_spin-polarized_2014}, the assumption on the meshes is removed by omitting the nodal projection from the magnetisation update.
However, the nodal projection is kept for the magnetisation appearing in the magnetoelastic terms for the sake of unconditional stability.
This modifications of the original algorithm of~\cite{bppr2013} give rise to additional errors that need to be controlled,
which makes the analysis more involved (e.g.\ a more accurate estimate of the projection error is needed; see Lemma~\ref{lem:projection_error_Lp} below).
Finally,} the discrete variational problems appearing in our integrator are standard and therefore easy to implement
in standard finite element packages.
For example, in the numerical experiments included in this work,
we use Netgen/NGSolve~\cite{netgen}.
\end{itemize}

The remainder of this work is organised as follows:
In Section~\ref{sec:model},
we present the PDE system we are interested in;
In Section~\ref{sec:ingredients},
we introduce the `ingredients' that are necessary for the definition of our numerical scheme and for its analysis;
In Section~\ref{sec:main},
we present our numerical integrator
(Algorithm~\ref{algorithm})
and state the main results of the work;
Section~\ref{sec:numerics} is devoted to numerical experiments.
In Section~\ref{sec:proofs}, we collect the proofs of all results.
For the convenience of the reader,
we conclude the paper with two appendices,
Appendix~\ref{sec:linear_algebra},
in which we collect several linear algebra definitions
and results used throughout the work,
and Appendix~\ref{sec:physics},
in which we show how to pass from the fully dimensional model
considered in the physics literature
to the dimensionless setting we study.

\section{Model problem}
\label{sec:model}

Let $\Omega\subset\R^3$ be a bounded Lipschitz domain
representing the volume occupied by a ferromagnetic body.
We assume the boundary $\partial\Omega$ is split into
two disjoint relatively open parts $\Gamma_D$
(of positive measure)
and $\Gamma_N$,
i.e.\
$\partial\Omega = \overline{\Gamma}_D \cup \overline{\Gamma}_N$
and $\Gamma_D \cap \Gamma_N = \emptyset$.
Let $T>0$ denote some final time.

The magnetomechanical state of the material
is described by two vector fields:
the displacement $\uu: \Omega \times (0,T) \to \R^3$
and the magnetisation $\mm : \Omega \times (0,T) \to \sphere$.
The total strain $\boldvar$ is made up of the elastic strain $\boldvarel$
and the magnetisation-dependent generally incompatible
(in the sense that it does not satisfy the Saint-Venant compatibility conditions~\cite{amrouche2006saint,Mballa2014})
magnetostrain $\boldvarm$, i.e.\ $\boldvar = \boldvarel + \boldvarm$.
The total strain is given by
\begin{equation*}
    \boldvar(\uu)
    = \frac{1}{2}\left(\Grad\uu + \Grad\uu^{\top}\right)
\end{equation*}
(strain-displacement relation).
Following~\cite{federico2019tensor},
we consider the expression
\begin{equation} \label{eq:magnetostrain}
    \boldvarm(\mm)
    = \Z : (\mm \otimes \mm),
\end{equation}
where $\Z \in \R^{3^4}$ is a fourth-order tensor, 
which we assume to be minorly symmetric
(i.e.\ $\Z_{ij\ell m} = \Z_{ji\ell m} = \Z_{ijm \ell}$ for all $i,j,\ell,m = 1,2,3$, cf.\ Appendix~\ref{sec:linear_algebra}).
It follows that
\begin{equation*}
    \boldvarel(\uu,\mm) = \boldvar(\uu) - \boldvarm(\mm).
\end{equation*}
The elastic part of the strain compensates for the magnetic part to make the total strain compatible~\cite{Mballa2014}.
The elastic strain is related to the stress tensor $\boldsig$ by Hooke's law
\begin{equation*}
    \boldsig(\uu,\mm) = \C : \boldvarel(\uu,\mm),
\end{equation*}
where $\C \in \R^{3^4}$ is the fourth-order, fully symmetric
(i.e.\ $\C_{ij\ell m} = \C_{\ell mij} = \C_{ji\ell m} = \C_{ij m \ell}$ for all $i,j,\ell,m = 1,2,3$, cf.\ Appendix~\ref{sec:linear_algebra}), positive definite stiffness tensor.
The elastic energy reads as
\begin{equation*}
    \E_{\mathrm{el}}[\uu,\mm]
    =
    \frac{1}{2} \int_{\Omega}[\boldvar(\uu) - \boldvarm(\mm)]
    : \{ \C : [\boldvar(\uu) - \boldvarm(\mm)] \}
    - \int_{\Omega}\ff\cdot\uu
    - \int_{\Gamma_{N}}\gg\cdot\uu,
\end{equation*}
where the last two terms
model the work done by
a volume force
$\ff : \Omega \to \R^3$
and a surface force
$\gg : \Gamma_N \to \R^3$
(traction),
both assumed to be constant in time.
The magnetic energy,
for simplicity assumed to comprise only the 
Heisenberg exchange contribution,
is given by
\begin{equation} \label{eq:energy_mag}
\E_{\mathrm{m}}[\mm]
= \frac{1}{2}\int_{\Omega}|\Grad \mm|^2.
\end{equation}
The total free energy of the system is defined as
the sum of the magnetic and elastic energies,
i.e.
\begin{equation}\label{eqn:definition_of_energy}
\begin{split}
    &\E[\uu,\mm]
    = \E_{\mathrm{m}}[\mm] + \E_{\mathrm{el}}[\uu,\mm] \\
    & \quad = \frac{1}{2}\int_{\Omega}|\Grad \mm|^2
    + \frac{1}{2} \int_{\Omega}[\boldvar(\uu) - \boldvarm(\mm)] : \{ \C : [\boldvar(\uu) - \boldvarm(\mm)] \}
    - \int_{\Omega}\ff\cdot\uu
    - \int_{\Gamma_{N}}\gg\cdot\uu.
\end{split}
\end{equation}
The dynamics of $\uu$ and $\mm$
is governed by the
coupled system of 
the conservation of (linear) momentum law
and the LLG equation
\begin{alignat}{2}
\label{eq:newton}
    \uutt
    & = \nabla \cdot \boldsig (\uu,\mm) + \ff
    && \text{in } \Omega \times (0,T),\\
\label{eq:llg}
    \mmt
    & = - \mm \times \heff[\uu,\mm]
    + \alpha \, \mm \times \mmt
    &\quad& \text{in } \Omega \times (0,T),
\end{alignat}
supplemented with the initial and boundary conditions
\begin{subequations} \label{eq:ibc}
\begin{alignat}{2}
    \label{eq:ic_u}
    \uu(0) &= \uu^0 &\quad&\text{in } \Omega,\\
    \label{eq:ic_ut}
    \uut(0) &= \dot\uu^0 &&\text{in } \Omega,\\    
     \label{eq:ic_m}
    \mm(0) &= \mm^0 &&\text{in } \Omega,\\
    \label{eq:bc_u_d}
    \uu &= \0 &&\text{on } \Gamma_D \times (0, T),\\
    \label{eq:bc_u_n}
    \boldsig\nn &= \gg && \text{on } \Gamma_{N} \times (0, T),\\
    \label{eq:bc_m}
    \partial_{\nn} \mm &= \0 && \text{on } \partial\Omega \times (0, T),    
\end{alignat}
\end{subequations}
where $\uu^0, \dot\uu^0 : \Omega \to \R^3$
and $\mm^0 : \Omega \to \sphere$
are suitable initial data,
while $\nn:\partial\Omega \to \sphere$
denotes the outward-pointing unit normal vector to $\partial\Omega$.
In~\eqref{eq:llg},
$\alpha>0$ denotes the Gilbert damping parameter, whereas
the effective field $\heff[\uu,\mm]$
is the variational derivative of the free energy
with respect to the magnetisation, i.e.
\begin{equation*}
    \heff[\uu,\mm]
    = -\frac{\delta \E[\uu,\mm]}{\delta \mm}
    = \Lapl\mm + \hhm[\uu,\mm],
\end{equation*}
where the elastic field reads as
\begin{equation} \label{eqn:MagnetostrictiveEffectiveField}
    \hhm[\uu,\mm]
    = 2 \, [\Z^\top : \boldsig(\uu,\mm)] \mm
    = 2 \, (\Z^\top : \{\C : [\boldvar(\uu) - \boldvarm(\mm)] \}) \mm,
\end{equation}
with $\Z^\top$ being the transpose of $\Z$ (cf.\ Appendix~\ref{sec:linear_algebra}).
Note that~\eqref{eq:newton} can be rewritten as
\begin{equation*}
    \uutt = -\frac{\delta \E[\uu,\mm]}{\delta \uu}.
\end{equation*}
A simple formal calculation reveals that sufficiently smooth solutions
to~\eqref{eq:newton}--\eqref{eq:ibc} satisfy the energy law
\begin{equation} \label{eq:energy_law}
\frac{\mathrm{d}}{\mathrm{d}t}
\left( \E[\uu(t),\mm(t)] + \frac{1}{2}\norm{\uut(t)}^2 \right)
=
- \alpha \norm{\mmt(t)}^2
\le 0,
\end{equation}
i.e.\ the sum of the total energy~\eqref{eqn:definition_of_energy}
(which can be understood as a potential energy)
and the kinetic energy $\norm{\uut}^2/2$
decays over time,
with the decay being modulated by $\alpha$.

For the data of the problem,
we assume that
$\C \in \LL^{\infty}(\Omega)$ is uniformly positive definite,
i.e.\ there exists $C_0>0$ such that
\begin{equation} \label{eq:tensor_coercivity}
\AA:(\C : \AA) \ge C_0 \norm{\AA}^2
\quad
\text{for all }
\AA \in \R^{3 \times 3},
\end{equation}
$\Z \in \LL^{\infty}(\Omega)$,
$\ff\in\LL^2(\Omega)$,
$\gg\in\LL^2(\Gamma_{N})$,
$\uu^0\in\HH^1(\Omega)$,
$\dot\uu^0\in\LL^2(\Omega)$,
and
$\mm^0\in\HH^1(\Omega;\sphere)$.
In the following definition,
we state the notion
of a weak solution of the initial boundary value problem~\eqref{eq:newton}--\eqref{eq:ibc};
see~\cite{cef2011}.
Hereafter, we shall denote $L^2$-integrals in space over some domain $D$
with $\inner[D]{\cdot}{\cdot}$,
omitting the subscript if $D = \Omega$.
Moreover, we denote by $\Omega_T$ the space-time cylinder $\Omega \times (0,T)$.

\begin{definition} \label{def:weak}
We say that a pair $(\uu,\mm) :\Omega_T\to\R^3 \times \R^3$
is a weak solution to the initial boundary value problem~\eqref{eq:newton}--\eqref{eq:ibc}
if the following conditions hold:
\begin{enumerate}[(i)]
    \item $\uu \in L^\infty(0,T;\HH^1_D(\Omega))$ with $\uut \in L^\infty(0,T;\LL^2(\Omega))$
    and
    $\mm \in L^\infty(0,T;\HH^1(\Omega;\sphere))$ with $\mmt \in L^2(0,T;\LL^2(\Omega))$;
    \item for all $\xxi \in \CC_{c}^{\infty}([0,T);\CC^{\infty}(\overline{\Omega}))$
    and $\vvphi \in \CC^{\infty}(\overline{\Omega_T})$,
    we have
    \begin{align}
    \label{eq:weak_u}
    \begin{split}
    & -\int_{0}^{T}\inner{\uut(t)}{\partial_{t}\xxi(t)}\mathrm{d}t
    + \int_{0}^{T} \inner{\mathbb{C} : [\boldvar(\uu(t)) - \boldvarm(\mm(t))]}{\boldvar(\xxi(t))} \mathrm{d}t \\
    & \qquad
    =
    \int_{0}^{T}\inner{\ff}{\xxi(t)} \mathrm{d}t
    + \int_{0}^{T}\inner[\Gamma_{N}]{\gg}{\xxi(t)}\mathrm{d}t
    + \inner{\dot{\uu}^{0}}{\xxi(0)},
    \end{split}\\
    \label{eq:weak_m}
    \begin{split}
    & \int_{0}^{T}\inner{\mmt(t)}{\vvphi(t)}\mathrm{d}t
    - \alpha \int_{0}^{T}\inner{\mm(t)\times\mmt(t)}{\vvphi(t)}\mathrm{d}t\\
    & \qquad
    =
    \int_{0}^{T}\inner{\mm(t)\times\Grad\mm(t)}{\Grad\vvphi(t)}\mathrm{d}t
    -\int_{0}^{T}\inner{\mm(t)\times\hhm[\uu(t),\mm(t)]}{\vvphi(t)}\mathrm{d}t;
    \end{split}
    \end{align}
    \item the initial conditions $\uu(0)=\uu^0$ and $\mm(0)=\mm^0$ hold in the sense of traces;
    \item for almost all $t' \in (0,T)$, it holds that
    \begin{equation}\label{eqn:EnergyLawDefinition}
        \E[\uu(t'),\mm(t')]
        + \frac{1}{2}\norm{\partial_t \uu(t')}^2
        + \alpha \int_{0}^{t'}\norm{\partial_t \mm(t)}^2 \mathrm{d}t
        \le
        \E[\uu^0,\mm^0]
        + \frac{1}{2}\norm{\dot{\uu}^{0}}^2.
    \end{equation}
\end{enumerate}
\end{definition}

Equations~\eqref{eq:weak_u} and~\eqref{eq:weak_m} are space-time
variational formulations of~\eqref{eq:newton} and~\eqref{eq:llg}, respectively.
The initial condition~\eqref{eq:ic_ut}
and the boundary conditions~\eqref{eq:bc_u_n} and~\eqref{eq:bc_m}
are imposed as natural boundary conditions in the variational formulations;
The initial conditions~\eqref{eq:ic_u} and~\eqref{eq:ic_m} are imposed in the sense
of traces in~{(iii)};
The Dirichlet boundary condition~\eqref{eq:bc_u_d} is imposed as essential boundary condition.
Equation~\eqref{eqn:EnergyLawDefinition} is the weak counterpart
of the energy law~\eqref{eq:energy_law} satisfied by strong solutions.

\begin{remark}
Formula~\eqref{eq:magnetostrain} is the general expression
of the magnetostrain for anisotropic ferromagnets~\cite{federico2019tensor}
and covers the typical forms of the magnetostrain found in literature.
These usually assume that the magnetostrain is \emph{isochoric}~\cite[Section~3.2.6]{hubert1998magnetic}
(i.e.\ it has zero trace).
In an isochoric material,
the magnetic body elongates (contracts) in the magnetisation direction,
and contracts (elongates) in the other
two for positive (negative) magnetostriction.
\comment{An example of positive magnetostriction
is shown in Figure \ref{fig:illustmagnetostriction}.}
Importantly, formula~\eqref{eq:magnetostrain} covers the common \emph{cubic} case,
considered e.g.\ in~\cite{james1998magnetostriction,shu2004micromagnetic,Mballa2014,rj2021,renuka2021solution}
and given by
\begin{equation*}
\boldvarm(\mm)
= \frac{3}{2} \Bigg\{
\lambda_{100} \left( \mm\otimes\mm - \frac{I}{3}\right)
+ 
(\lambda_{111} - \lambda_{100})
\sum_{\substack{i,j=1 \\ i \neq j}}^3 (\mm\cdot\ee_i^\mathrm{c})(\mm\cdot\ee_j^\mathrm{c})
(\ee_i^\mathrm{c} \otimes \ee_j^\mathrm{c})
\Bigg\},
\end{equation*}
where $I \in \R^{3 \times 3}$ denotes the $3$-by-$3$ identity matrix,
$\lambda_{100},\lambda_{111} \in \R$ are material constants,
and $\{\ee_1^\mathrm{c}, \ee_2^\mathrm{c}, \ee_3^\mathrm{c}\}$ denotes
an orthonormal set yielding the crystal basis.
When $\lambda_{100}=\lambda_{111}$, the latter reduces to the so-called \emph{isotropic} case
\begin{equation}\label{eqn:magnetostrain_isotropic}
\boldvarm(\mm)
= \frac{3}{2} \lambda_{\textnormal{100}} \left( \mm\otimes\mm - \frac{I}{3}\right),
\end{equation}
considered e.g.\ in~\cite{bhbvfs2014,pwhcn2015,dw2023}.
For further details regarding specific crystal classes and their magnetostrain representation,
we refer to~\cite{federico2019tensor}.
\end{remark}

\begin{remark}
For the sake of simplicity
(and since the focus of this work is on the design of a numerical method
for the coupled system~\eqref{eq:newton}--\eqref{eq:llg}),
we neglect from the magnetic energy~\eqref{eq:energy_mag}
all lower-order contributions
(magnetocrystalline anisotropy, Zeeman energy, magnetostatic energy, Dzyaloshinskii–-Moriya interaction).
However, we note that their numerical integration is well understood;
see e.g.\ \cite{bffgpprs2014,dpprs2019,hpprss2019}.
\end{remark}

\section{Preliminaries}
\label{sec:ingredients}

In this section, we collect some notation and preliminary results that will be necessary
to introduce and analyse the fully discrete algorithm we propose to approximate solutions
to the initial boundary value problem~\eqref{eq:newton}--\eqref{eq:ibc}.
Hereafter, as customary in numerical analysis, given $A,B \in \R$,
we shall write $A \lesssim B$ if there exists a constant $c>0$, clear from the context
and always independent of the discretisation parameters, such that $A \le c \, B$.

\subsection{Time discretisation}
Let $0=t_{0}<t_{1}<\dots<t_{N}=T$ be a uniform partition of the time interval
into $N$ uniform intervals with constant time-step size $k = T/N$,
i.e.\ $t_i = ik$ for all $i=0,\dots,N$.
Given values $\{\phi^i\}_{0 \le i \le N}$ and $\dot\phi^0$,
we define the discrete time derivatives by
\begin{align}
\dt \phi^i &:=
\begin{cases}
\dot\phi^0, & \text{if } i=0,\\
\displaystyle
\frac{\phi^i - \phi^{i-1}}{k}, & \text{if } 1 \le i \le N,
\end{cases} \\
\label{eq:second_derivative}
\dt^2 \phi^{i+1}
&:= \frac{\dt \phi^{i+1} - \dt\phi^i}{k}
=
\begin{cases}
\displaystyle
\frac{\phi^1 - \phi^0 - k \dot\phi^0}{k^2}, & \text{if } i=0,\\
\displaystyle
\frac{\phi^{i+1} - 2\phi^i+\phi^{i-1}}{k^2}, & \text{if } 1 \le \comment{i} \le N-1.
\end{cases}
\end{align}
Moreover,
we define the time reconstructions
$\phi_{k}$, $\phi^-_{k}$, $\phi^+_{k}$, $\dot\phi_{k}$, $\dot\phi^-_{k}$, $\dot\phi^+_{k}$,
defined, for all $0 \le i \le N-1$ and $t \in [t_i,t_{i+1})$, by
\begin{subequations} \label{eq:reconstructions}
\begin{gather}
\phi_{k}(t) := \frac{t-t_i}{k}\phi^{i+1} + \frac{t_{i+1} - t}{k}\phi^i,
\quad
\phi_{k}^-(t) := \phi^i,
\quad
\phi_{k}^+(t) := \phi^{i+1}, \\
\dot\phi_{k}(t) := \dt\phi^i + (t - t_i)\dt^2\phi^{i+1},
\quad
\dot\phi_{k}^-(t) := \dt\phi^i,
\quad
\dot\phi_{k}^+(t) := \dt\phi^{i+1}.
\end{gather}
\end{subequations}
Note that $\partial_t \phi_k (t) = \dot\phi^+_{k}(t) = d_t \phi^{i+1}$ for all $t \in [t_i,t_{i+1})$.

\subsection{Space discretisation}
Let $\Omega$ be a polyhedral domain.
Let $\{\T_h\}_{h>0}$ be a shape-regular family of meshes of $\Omega$ into tetrahedra,
where $h=\max_{K\in\T_{h}}h_K$ denotes the mesh size of $\T_{h}$
and $h_K = \diam K$ for all $K \in \T_h$.
We denote by $\NN_{h}$ the set of nodes in the triangulation $\T_{h}$.
For all $K \in \T_h$, we denote by $\mathcal{P}_{1}(K)$ the space of polynomials of degree at most 1 over $K$.
We denote by $\S^1(\T_{h})$ the space of piecewise affine and globally continuous functions from $\Omega$ to $\R$,
i.e.
\begin{equation*}
\S^1(\T_{h}) = \{\phi_{h}\in C(\overline{\Omega}):\phi_{h}\vert_{K}\in \mathcal{P}_{1}(K)
\text{ for all } K\in\T_{h}\} \subset H^1(\Omega).
\end{equation*}
We denote by $\II_{h}:C(\overline{\Omega})\to \S^1(\T_{h})$ the nodal interpolant
satisfying $\II_{h}[\phi](z) = \phi(z)$ for each $z\in\NN_{h}$,
where $\phi$ is a continuous function.
Moreover, we consider the space
$\S^{1}_{D}(\T_{h}) = \S^{1}(\T_{h})\cap H_{D}^{1}(\Omega)$,
where homogeneous Dirichlet boundary conditions on $\Gamma_{D}$ are imposed explicitly.

Since the unknowns of the problem in which we are interested are vector fields,
we consider the vector-valued finite element space $\S^1(\T_{h})^3$
and use the same notation adopted in the scalar case
to denote the vector-valued nodal interpolant $\II_{h}:\CC(\overline{\Omega})\to \S^1(\T_{h})^3$.
For all $0 \le i \le N$, the approximate displacement at time $t_i$, $\uu_h^i \approx \uu(t_i)$,
will be sought in the finite element space $\S^1_D(\T_{h})^3$,
whereas the approximate magnetisation, $\mm_h^i \approx \mm(t_i)$, will be sought in the set
\begin{equation} \label{eq:discrete_magnetisation}
\MM_{h,\delta}
= \big\{\pphi_{h}\in \S^1(\T_{h})^3: \abs{\pphi_{h}(z)}\geq 1\text{ for all } z\in\NN_{h} \text{ and } \norm[L^1(\Omega)]{\II_{h}\left[|\pphi_h|^2\right] - 1}\leq \delta \big\}
\end{equation}
for some $\delta>0$.
Note that discrete magnetisations in $\MM_{h,\delta}$ generally do not satisfy the unit length constraint,
not even at the vertices of the mesh,
but the error is controlled in the $L^1$-sense by $\delta$.
For the case $\delta=0$, we obtain the set
\begin{equation*}
\MM_{h,0} = \{\pphi_{h}\in \S^1(\T_{h})^3:\abs{\pphi_{h}(z)}=1\text{ for all } z\in\NN_{h}\},
\end{equation*}
in which the constraint holds at the vertices of the mesh.
We define the nodal projection operator
$\Pi_{h}:\MM_{h,\delta}\to\MM_{h,0}$ by $\Pi_{h}\pphi_{h}(z) = \pphi_{h}(z)/|\pphi_{h}(z)|$
for all $z\in\NN_{h}$ and $\pphi_{h}\in\MM_{h,\delta}$.

Another important property of solutions to the LLG equation is
the orthogonality $\mmt\cdot\mm = 0$.
To realise it at the discrete level,
given an approximation $\mm_h^i \approx \mm(t_i)$ in $\MM_{h,\delta}$,
we consider the discrete tangent space
\begin{equation*}
\KK_{h}[\mm_{h}^i]
= \{\ppsi_{h}\in \S^1(\T_{h})^3: \mm_{h}^i(z)\cdot\ppsi_{h}(z) = 0 \text{ for all } z\in\NN_{h}\},
\end{equation*}
where approximations $\vv_h^i \approx \mmt(t_i)$ will be sought.
Note that the desired orthogonality property is imposed only at the vertices of the mesh.

To conclude, we recall the definition of mass-lumped $L^2$-product
$\inner[h]{\cdot}{\cdot}$, i.e.
\begin{equation} \label{eq:mass-lumping}
\inner[h]{\ppsi}{\pphi}
= \int_\Omega \II_h[\ppsi \cdot \pphi]
\quad \text{for all } \ppsi, \pphi \in \CC^0(\overline{\Omega}),
\end{equation}
which is a scalar product on $\S^1(\T_h)^3$.

\section{Algorithm and main results}
\label{sec:main}

In the following algorithm,
we state the fully discrete numerical scheme we propose
to approximate solutions
to the initial boundary value problem~\eqref{eq:newton}--\eqref{eq:ibc}.

\begin{algorithm}[decoupled algorithm for the LLG equation with magnetostriction]
\label{algorithm}
\underline{Discretisation parameters:}
Mesh size $h>0$, time-step size $k>0$, $\theta\in(1/2,1]$.\\
\underline{Input:}
Approximate initial conditions
$\mm_{h}^{0} \in \MM_{h,0}$,
$\uu_{h}^{0} \in \S^1_D(\T_h)^3$,
$\dot\uu_{h}^{0} \in \S^1(\T_h)^3$. \\
\underline{Loop:}
For all integers $0 \le i \le N-1$,
iterate {\rm(i)--(iii)}:
\begin{enumerate}[(i)]
\item Compute $\vv_{h}^{i}\in\KK_{h}[\mm_{h}^{i}]$ such that, for all $\pphi_{h}\in\KK_{h}[\mm_{h}^{i}]$, it holds that
\begin{multline}\label{alg1:magnetisation_update}
    \alpha \inner[h]{\vv_{h}^{i}}{\pphi_{h}}
    + \inner[h]{\mm_{h}^{i}\times \vv_{h}^{i}}{\pphi_{h}}
    + \theta k \inner{\Grad\vv_{h}^{i}}{\Grad \pphi_{h}} \\
    = -\inner{\Grad\mm_{h}^{i}}{\Grad \pphi_{h}}
    + \inner{\hhm[\uu_{h}^{i},\Pi_h\mm_{h}^{i}]}{\pphi_{h}}.
\end{multline}
\item Define
\begin{equation} \label{eq:alg_update_m}
    \mm_{h}^{i+1} := \mm_{h}^{i} + k \vv_{h}^{i} \in \S^1(\T_h)^3.
\end{equation}
\item Compute $\uu_{h}^{i+1}\in\S^{1}_{D}(\T_{h})^3$ such that,
for all $\ppsi_{h}\in\S^{1}_{D}(\T_{h})^3$, it holds that
\begin{multline}\label{alg1:displacement_update}
    \inner{\dt^2 \uu_{h}^{i+1}}{\ppsi_h} +  \inner{\C:\boldvar(\uu_{h}^{i+1})}{\boldvar(\ppsi_{h})} \\
    = \inner{\C:\boldvarm(\Pi_h\mm_{h}^{i+1})}{\boldvar(\ppsi_{h})}
    + \inner{\ff}{\ppsi_{h}}
    + \inner[\Gamma_{N}]{\boldsymbol{g}}{\ppsi_{h}}.
\end{multline}
\end{enumerate}
\underline{Output:}
Approximations $\{ (\uu_h^i,\mm_h^i) \}_{0 \le i \le N}$.
\end{algorithm}

Algorithm~\ref{algorithm} resembles the decoupled algorithm proposed in~\cite{bppr2013}.
The discrete initial data
$\mm_{h}^{0} \in \MM_{h,0}$, $\uu_{h}^{0} \in \S^1_D(\T_h)^3$ and $\dot\uu_{h}^{0} \in \S^1(\T_h)^3$
denote suitable approximations of the initial conditions $\mm^0$, $\uu^0$ and $\dot\uu^0$, respectively.
For every time-step,
given current approximations of the magnetisation and the displacement,
we compute the new magnetisation first,
and then the updated displacement using this.

Specifically, to compute the new magnetisation,
we use the tangent plane scheme~\cite{aj2006,bkp2008,alouges2008a}:
In step~(i), given $\uu_{h}^i$ and $\mm_{h}^i$,
we compute an approximation $\vv_{h}^{i} \approx \mmt(t_i)$
residing in the discrete tangent space $\KK_{h}[\mm_{h}^{i}]$.
The variational problem~\eqref{alg1:magnetisation_update} solved by $\vv_{h}^{i}$
is a discretisation of the equivalent formulation of the LLG equation
\begin{equation} \label{eq:equivalent_llg}
\alpha \, \mmt + \mm \times \mmt = \heff[\uu,\mm] - (\heff[\uu,\mm]\cdot\mm)\mm,
\end{equation}
which can be obtained from~\eqref{eq:llg} via simple algebraic manipulations; cf.\ \cite{aj2006}.
Looking at~\eqref{alg1:magnetisation_update},
we note that the discrete variational formulation of the left-hand side of~\eqref{eq:equivalent_llg}
makes use of the mass-lumped $L^2$-product~\eqref{eq:mass-lumping}.
The two terms constituting the effective field $\heff[\uu,\mm]$ are treated differently:
The exchange contribution is treated implicitly and therefore contributes to the left-hand side
of~\eqref{alg1:magnetisation_update}.
The `degree of implicitness' is modulated by the parameter $\theta \in (1/2,1]$.
The elastic field is treated explicitly.
In step~(ii),
with $\vv_h^i$ at hand,
we compute the new magnetisation $\mm_{h}^{i+1}$ using a first-order time-stepping;
cf.\ \eqref{eq:alg_update_m}.
Differently from the seminal papers on the tangent plane schemes~\cite{aj2006,bkp2008,alouges2008a}
and from~\cite{bppr2013},
we follow the approach of~\cite{bartels2016,abert_spin-polarized_2014}
and in our update we do not use the nodal projection.
In particular, it holds that $\dt\mm_h^{i+1} = \vv_h^i$.
Finally, in step~(iii),
we compute the new displacement $\uu_{h}^{i+1}$
using a standard finite element discretisation of~\eqref{eq:newton}.
We use the backward Euler method in time 
(the second time derivative in~\eqref{eq:newton} is approximated
using the different quotient~\eqref{eq:second_derivative}).

In Algorithm~\ref{algorithm},
we apply the nodal projection to all approximate magnetisations
arising from the elastic energy,
i.e.\
in the elastic field on the right-hand side of~\eqref{alg1:magnetisation_update}
and in the magnetostrain term on the right-hand side of~\eqref{alg1:displacement_update},
whereas the nodal projection is omitted from the magnetisation
in the exchange field on the right-hand side of~\eqref{alg1:magnetisation_update}, the cross product on the left hand side of~\eqref{alg1:magnetisation_update},
and from the update~\eqref{eq:alg_update_m}.

Notably,
despite the nonlinearity of the LLG equation
and its nonlinear coupling with the conservation of momentum law,
Algorithm~\ref{algorithm} is \emph{fully linear}
and only requires the solution of two linear systems per time-step.

\comment{\begin{remark}
In Algorithm~\ref{algorithm},
the magnesation update is based on the projection-free tangent plane scheme~\cite{bartels2016,abert_spin-polarized_2014}.
As our analysis below will show, the unit length constraint is imposed inexactly,
but the constraint violation error can be controlled by the time-step size.
Moreover, passing the estimate to the limit, we can show that the constraint is satisfied by the weak solution towards which the finite element approximations are converging.
Other approaches in the literature aim to impose the constraint exactly
(at least at the vertices of the underlying finite element mesh).
This can be achieved by projecting the magnetisation onto the sphere after each update~\cite{aj2006,bkp2008,alouges2008a},
using constraint-preserving variational formulation~\cite{bp2006}
or designing magnetisation updates based on exponential map~\cite{ln2003}.
In Algorithm~\ref{algorithm}, we refrain from these approaches
as they lead to geometric restrictions on the finite element meshes~\cite{aj2006,bkp2008,alouges2008a}
or lead to the solution of nonlinear systems of equations at each time-step~\cite{ln2003,bp2006}.
\end{remark}}

In the following proposition, we show the well-posedness of Algorithm~\ref{algorithm}.
The proof, based on standard arguments, is postponed to Section~\ref{sec:proofs_wellposedness}.

\begin{proposition}
\label{prop:wellposedness}
Algorithm~\ref{algorithm} is well defined for every $\theta\in(1/2,1]$,
i.e.\ for every integer $0 \le i \le N-1$,
there exists a unique
$(\vv_{h}^{i},\mm_{h}^{i+1},\uu_{h}^{i+1}) \in \KK_{h}[\mm_{h}^{i}] \times \S^{1}(\T_{h})^3 \times \S^{1}_{D}(\T_{h})^3$
satisfying~\eqref{alg1:magnetisation_update}--\eqref{alg1:displacement_update}.
\end{proposition}

In the following proposition,
we establish a discrete counterpart of the energy law~\eqref{eq:energy_law}
satisfied by smooth solutions of the continuous problem
(see also~\eqref{eqn:EnergyLawDefinition} for the corresponding property for weak solutions).
Its proof is postponed to Section~\ref{sec:proofs_energy_law}.

\begin{proposition} \label{prop:energy}
For every integer $0 \le i \le N-1$,
the iterates of Algorithm~\ref{algorithm} satisfy the discrete energy law
\begin{equation} \label{eq:discrete_energy_law}
\E[\uu_h^{i+1},\mm_h^{i+1}] + \frac{1}{2} \norm{\dt\uu_h^{i+1}}^2
-
\E[\uu_h^i,\mm_h^i] - \frac{1}{2} \norm{\dt\uu_h^i}^2
=
- \alpha k \norm[h]{\vv_h^i}^2
- D_{h,k}^i
- E_{h,k}^i,
\end{equation}
where $D_{h,k}^i$
and $E_{h,k}^i$ are given by
\begin{multline} \label{eq:dissipation}
D_{h,k}^i
= k^2 (\theta-1/2) \norm{\Grad\vv_h^i}^2
+ \frac{1}{2} \norm{\dt\uu_h^{i+1} - \dt\uu_h^i}^2 \\
+ \frac{1}{2} \norm[\C]{[\boldvar(\uu_{h}^{i+1})-\boldvarm(\mm_{h}^{i+1})] - [\boldvar(\uu_{h}^i)-\boldvarm(\mm_{h}^i)]}^2 \ge 0
\end{multline}
and
\begin{equation} \label{eq:error_energy}
\begin{split}
E_{h,k}^i
& = k^2 \inner{\C:[\boldvar(\uu_{h}^{i+1}) - \boldvarm(\mm_{h}^{i+1})]}{\boldvarm(\vv_{h}^{i})} \\
& \quad 
+ 2k \inner{\C:\{[\boldvar(\uu_{h}^{i+1}) - \boldvarm(\mm_{h}^{i+1})] - [\boldvar(\uu_{h}^{i}) - \boldvarm(\mm_{h}^{i})] \}}{\Z(\mm_{h}^{i} \otimes \vv_{h}^{i})} \\
& \quad
+ 2k\inner{\C:[\boldvar(\uu_{h}^{i}) - \boldvarm(\mm_{h}^{i})]}{\Z[(\mm_h^i - \Pi_h \mm_h^i) \otimes \vv_{h}^{i}]}\\
& \quad
+ \inner{\C:[\boldvarm(\mm_{h}^{i+1}) - \boldvarm(\Pi_h\mm_{h}^{i+1})]}{\boldvar(\uu_h^{i+1}) - \boldvar(\uu_h^{i})}.
\end{split}
\end{equation}
respectively.
\end{proposition}

In~\eqref{eq:dissipation},
we use the norm $\norm[\C]{\cdot}^2 = \inner{\C:(\cdot)}{\cdot}$
for matrix-valued functions in $L^2(\Omega)^{3 \times 3}$.
Thanks to our assumptions on $\C$ (cf.\ \eqref{eq:tensor_coercivity}),
this norm is equivalent to the standard $L^2$-norm.

Looking at the right-hand side of~\eqref{eq:discrete_energy_law},
we see that the inherent $\alpha$-modulated energy dissipation
of the model (cf.\ \eqref{eq:energy_law})
is spoiled by two terms:
\begin{itemize}
\item
the artificial damping $D_{h,k}^i$,
arising from
the implicit treatment of the exchange contribution of the effective field in~\eqref{alg1:magnetisation_update} (the first term)
and the use of the backward Euler method
in~\eqref{alg1:displacement_update} (the last two terms),
\item
the error $E_{h,k}^i$ due to
linearisation (the first term)
decoupling (the second term),
and use of the nodal projection to impose the unit length constraint
on the magnetisations appearing in the elasticity terms
(the third and fourth terms).
\end{itemize}

\begin{remark}
Our argument to show Proposition~\ref{prop:energy}
for Algorithm~\ref{algorithm} can be transferred to
the algorithm of~\cite{bppr2013},
hence a by-product of our analysis is a discrete energy law for that algorithm.
Due to the use of the nodal projection in~\cite{bppr2013},
the counterpart of~\eqref{eq:discrete_energy_law}
is only an inequality (not an identity),
its proof requires to assume that the mesh is weakly acute,
and the error term $E_{h,k}^i$ 
does not include the last two terms in~\eqref{eq:error_energy}.
\end{remark}

Now, we discuss the stability and the convergence of Algorithm~\ref{algorithm}.
To this end, we consider the following convergence assumption on the approximate initial conditions:
\begin{equation} \label{eq:convergence_initial_data}
    \uu_h^0 \to \uu^0 \text{ in } \HH^1(\Omega), \ \
    \dot\uu_h^0 \to \dot\uu^0 \text{ in } \LL^2(\Omega), \ \
    \text{and} \ \
    \mm_h^0 \to \mm^0 \text{ in } \HH^1(\Omega), \ \
    \text{ as } h \to 0.
\end{equation}
Firstly,
we can show that Algorithm~\ref{algorithm} is unconditionally stable
and that the error in the unit length constraint
can be controlled by the time-step size.

\begin{proposition} \label{prop:stability}
Suppose that assumption~\eqref{eq:convergence_initial_data} is satisfied.
There exists a threshold $k_0>0$ such that, if $k < k_0$,
for every integer $1 \le j \le N$, the iterates of Algorithm~\ref{algorithm}
satisfy
\begin{multline} \label{eq:stability}
\norm{\dt \uu_{h}^j}^2
+ \norm{\boldvar(\uu_{h}^{j})}^2
+ \sum_{i=0}^{j-1}\norm{\dt \uu_{h}^{i+1} - \dt\uu_{h}^i}^2
+ \sum_{i=0}^{j-1}\norm{\boldvar(\uu_{h}^{i+1})-\boldvar(\uu_{h}^{i})}^2 \\
+ \norm[\HH^1(\Omega)]{\mm_{h}^{j}}^2
+ k \sum_{i=0}^{j-1} \norm{\vv_{h}^{i}}^2
+ \left(\theta - \frac{1}{2}\right) k^2\sum_{i=0}^{j-1} \norm{\Grad\vv_{h}^{i}}^2
\leq C
\end{multline}
and
\begin{equation}
\label{eq:general_constraint}
\big\lVert \II_h\big[\abs{\mm_h^j}^2\big]-1 \big\rVert_{L^1(\Omega)}
\le C k.
\end{equation}
The threshold $k_0>0$ and the constant $C>0$ depend only on
the shape-regularity parameter of $\T_h$,
the problem data $\alpha$, $T$, $\Omega$, $\C$, $\Z$, $\ff$ and $\gg$,
and the uniform bounds of the energy of the approximate initial data guaranteed by~\eqref{eq:convergence_initial_data}.
\end{proposition}

For the proof of the result, we refer to Section~\ref{sec:proofs_stability}.
Note that~\eqref{eq:general_constraint} implies that,
if the time-step size is sufficiently small,
the approximate magnetisations generated by the algorithm
belong to the set $\MM_{h,\delta}$ from~\eqref{eq:discrete_magnetisation}
with $\delta = Ck$.

With the approximations generated by Algorithm~\ref{algorithm},
we can construct the piecewise affine time reconstructions
$\uu_{hk} : (0,T) \to \S^1(\T_h)^3$
and
$\mm_{hk} : (0,T) \to \S^1(\T_h)^3$;
see~\eqref{eq:reconstructions}.
In the following theorem,
we show that
the sequences $\{\uu_{hk}\}$ and $\{\mm_{hk}\}$
converge in a suitable sense towards a weak solution of
the initial boundary value problem~\eqref{eq:newton}--\eqref{eq:ibc}
as $h$, $k$ go to $0$.
Its proof is postponed to Sections~\ref{sec:proofs_convergence}--\ref{sec:proofs_energy_inequality}.

\begin{theorem} \label{thm:convergence}
Suppose that assumption~\eqref{eq:convergence_initial_data} is satisfied.
\begin{enumerate}[(i)]
\item
There exist a weak solution $(\uu,\mm)$
of~\eqref{eq:newton}--\eqref{eq:ibc}
in the sense of Definition~\ref{def:weak}(i)--(iii)
and a (nonrelabeled) subsequence of $\{ (\uu_{hk},\mm_{hk}) \}$
which converges towards $(\uu,\mm)$ as $h,k \to 0$.
In particular, as $h,k \to 0$, it holds that
$\uu_{hk} \weakstarto \uu$ in $L^{\infty}(0,T; \HH^1_D(\Omega))$,
$\partial_t\uu_{hk} \weakstarto \comment{\partial_t}\uu$ in $L^{\infty}(0,T; \LL^2(\Omega))$,
$\mm_{hk} \weakstarto \mm$ in $L^{\infty}(0,T; \HH^1(\Omega ;\sphere))$,
and
$\partial_t\mm_{hk} \weakto \comment{\partial_t}\mm$ in $\LL^2(\Omega_T)$.
\item
If the discretisation parameters additionally satisfy the CFL condition $k=o(h^9)$,
the weak solution from part~(i) satisfies
the energy inequality~\eqref{eqn:EnergyLawDefinition} from Definition~\ref{def:weak}(iv).
\end{enumerate}
\end{theorem}

The proof of Theorem~\ref{thm:convergence} is constructive
and provides also a proof of existence of weak solutions.
We recall that, due to the non-convex nature of the problem,
uniqueness of weak solutions cannot be expected
(cf.\ the explicit proof of non-uniqueness of weak solutions to the pure LLG equation in~\cite{as1992}).
Moreover, if $\theta \in [0,1/2]$, then Theorem~\ref{thm:convergence} still holds,
but with an additional CFL condition for part~(i),
i.e.\ $k=o(h^2)$ if $\theta \in [0,1/2)$ and $k=o(h)$ if $\theta = 1/2$;
see~\cite{alouges2008a}.

\begin{remark}
The application of the nodal projection to all approximate magnetisations
arising from the elastic energy
is responsible for two of the error terms in~\eqref{eq:error_energy}
and for the severe CFL condition in Theorem~\ref{thm:convergence}(ii)
(cf.\ the analysis in Section~\ref{sec:proofs_energy_inequality} below),
so one would be tempted to completely remove it.
However, we believe that a fully projection-free approach would not lead
to an unconditionally stable method.
In particular, the use of the nodal projection on the outermost magnetisation
in the elastic field (cf. \eqref{eqn:MagnetostrictiveEffectiveField})
is non-negotiable as the total strain $\boldvar(\uu)$ is only in $\LL^2(\Omega)$.
For a stable method, it would be sufficient to take only one projection, not two,
within the magnetostrain as this would yield the estimate
$\norm{\Z: (\Pi_{h}\mm_{h} \otimes \mm_{h})} \lesssim \norm{\mm_{h}}$,
which would allow for the stability estimate of Proposition~\ref{prop:stability}.
However, we prefer not to use this approach as it would introduce some `unnatural' non-symmetry.
\end{remark}

\begin{remark}
The proof of the energy inequality typically requires extra assumptions to be proven.
In~\cite[Appendix A]{bffgpprs2014}, in the case of the LLG equation (with full effective field),
its proof requires higher regularity and stronger convergence assumptions on the applied field and general contribution terms.
In~\cite[Theorem 3.2]{dpprs2019}, in the case of the coupled system of the LLG equation and the eddy current equation,
a similar situation arises with a CFL condition $k=o(h^{3/2})$.
The very severe CFL condition in Theorem~\ref{thm:convergence}(ii) is an artifact of the analysis
and is due to the nonlinearity of the coupling and the fact that our proof requires explicit estimates of the error associated
with the use of the nodal projection in the elastic terms.
In particular, we need to estimate this error in a norm that is stronger than the $L^1$-norm,
which leads to a reduced convergence rate with respect to the time-step size $k$ (see Lemma~\ref{lem:projection_error_Lp} below).
This, combined with the fact that we need inverse estimates to obtain quantities we are able to control, leads to the CFL condition.
For more details, we refer to the proof of the result in Section~\ref{sec:proofs_energy_inequality} below.
However, we stress that this restriction does not show up within the numerics (see, in particular, the experiment Section~\ref{sec:numerics_cfl}).
\end{remark}

\section{Numerical experiments}
\label{sec:numerics}

In this section, to show the applicability of our algorithm,
we present a collection of numerical experiments.
The implementation of Algorithm~\ref{algorithm} was written
using the Netgen/NGSolve package \cite{netgen}
using version 6.2.2302.
The solution of the constrained linear system~\eqref{alg1:magnetisation_update}
in Algorithm~\ref{algorithm} is based on the null-space method given
in~\cite{ramage2013preconditioned,kraus2019iterative}.
The resulting system is solved using GMRES with an incomplete
LU decomposition preconditioner, with the previous linear update
$\vv_{h}^{i-1}$ as a starting guess. The elastic equation~\eqref{alg1:displacement_update}
is solved using a Jacobi preconditioned conjugate gradient method.
All computations were made on an i5-9500 CPU with 16GB of installed memory. 

\subsection{Material parameters}

In the upcoming numerical experiments,
we use material parameters estimated for \ch{(Fe90Co10)78Si12B10} (which we shall call FeCoSiB) from~\cite{dw2022}.
For the mass density and the Gilbert damping parameter (needed in our model, but not in~\cite{dw2022}),
we take the values used in~\cite{kohmoto1985mass} and~\cite{hu2023temperature}, respectively.
The resulting exchange length is $\ellex = \sqrt{2A/(\mu_0 M_s^2)} \approx 3\cdot 10^{-9}$ \si{m}.
The stiffness tensor $\C$ is assumed to be isotropic and
acts on symmetric matrices $\boldvar$ (the only type required) as
\begin{equation*}
\C:\boldvar = 2\mu \, \boldvar + \lambda \trace(\boldvar) \,I,
\end{equation*}
where $\mu$ and $\lambda$ are referred to as Lam\'e constants
(for FeCoSiB after non-dimensionalisation we have $\mu\approx 6.89$ and $\lambda \approx 21.96$).
For the magnetostrain, we consider the expression in~\eqref{eqn:magnetostrain_isotropic}.
In some experiments, the magnetic energy~\eqref{eq:energy_mag} will be supplemented
with the term $-\inner{\hh_{\mathrm{ext}}}{\mm}$ (Zeeman energy),
modelling the interaction of the magnetisation with an applied external field $\hh_{\mathrm{ext}}$.
For the sake of reproducibility, the values used are reported in Table~\ref{tab:valuetable}
(we refer to Appendix~\ref{sec:physics} for the relationship between the fully dimensional model
and the dimensionless setting of this paper).

\begin{table}[ht]
\begin{tabular}{|c|l|c|}
\hline
Symbol          & Name                       & Value                                        \\ \hline
$A$             & Exchange constant          & $1.5\cdot 10^{-11}$ \si{J.m^{-1}}            \\ \hline
$\alpha$        & Gilbert damping parameter  & $0.005$                                       \\ \hline
$\gamma$        & Gyromagnetic ratio         & $1.761 \cdot 10^{11}$ \si{rad.s^{-1}.T^{-1}} \\ \hline
$\mu_0$         & Permeability of free space & $1.25663706 \cdot 10^{-6}$                   \\ \hline
$M_s$           & Saturation magnetisation   & $1.5 \cdot 10^{6}$ \si{A.m^{-1}}             \\ \hline
$\lambda_{100}$ & Saturation magnetostrain   & $30 \cdot 10^{-6}$                           \\ \hline
$\rho$          & Density                    & $7900$ \si{kg.m^{-3}}                \\ \hline
$g$             & Gravitational acceleration & $9.81$ \si{m.s^{-2}}                         \\ \hline
$\mu$           & First Lam\'e constant      & $172$ \si{GPa}                               \\ \hline
$\lambda$       & Second Lam\'e constant     & $54$ \si{GPa}                                \\ \hline
\end{tabular}
\vspace{0.1cm}
\caption{Estimated material parameters for FeCoSiB taken from~\cite{dw2022,kohmoto1985mass,hu2023temperature}.}
\label{tab:valuetable}
\end{table}

\subsection{Magnetoelastic coupling} \label{sec: experiments_physics}

In this section, we present two numerical experiments
aimed at showcasing the capability of Algorithm~\ref{algorithm} to simulate
physical processes involving magnetoelastic materials.

\begin{figure}[t!]
    \centering
\includegraphics[height=2.5cm]{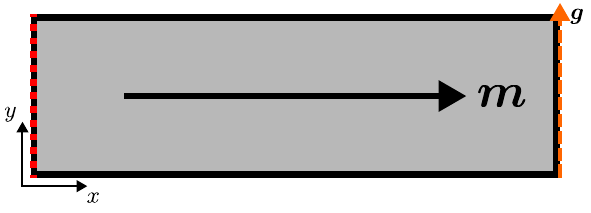}
    \caption{Experiments of Section~\ref{sec: experiments_physics}:
    View from above of the FeCoSiB bar of dimensions $(20\ellex,6\ellex,6\ellex)$.}
    \label{fig:xy_box_schematic}
\end{figure}

The simulation object is a bar of FeCoSiB,
clamped at one end ($y=0$ plane),
shown in Figure \ref{fig:xy_box_schematic}.
The bar has a physical length
of $20\ellex$ and width/height of $6\ellex$.
The maximum mesh size is
$h_{\max} \approx 0.9 \ellex$
(thereby being below the exchange length).
The initial magnetisation is uniformly
in the $x$-direction $\mm_{h}^{0} = (1,0,0)$,
whereas we set zero initial displacement $\uu_{h}^{0} = \boldsymbol{0}$ with
zero initial velocity $\dot{\uu}_{h}^{0} = \boldsymbol{0}$.
Gravity is enabled and implemented as a volume force $\ff = (0,0,-g)$, with
a value of $-g = -2.97 \cdot 10^{-14}$ after non-dimensionalisation.
If enabled, tractions (represented by a surface force $\gg$ applied on $\Gamma_N$)
and applied external fields $\hh_{\mathrm{ext}}$
are applied along the $+y$ direction.
Simulations are run
for $1$ \si{ns}, using time-steps of size $2\cdot 10^{-12}$ \si{s}.
This corresponds to a non-dimensional time length of
$T \approx 330$ and time-step $k \approx 0.66$.

\subsubsection{Direct magnetostrictive effect}\label{sec:applied_magnetic_field}

In this experiment,
we show that changes in the magnetisation yield changes in the mechanical state of the body.
To this end, we neglect traction and apply a uniform applied external field $\hh_{\mathrm{ext}}$ along the $+y$ direction
with low values of $0,1\cdot 10^{-4},3\cdot 10^{-4},5\cdot 10^{-4},7\cdot 10^{-4}$,
which corresponds to fields of strength $0$, $0.2$, $0.6$, $0.9$, $1.3$ \si{\milli\tesla}.
The fields are weak so that
the dynamics is not too fast.

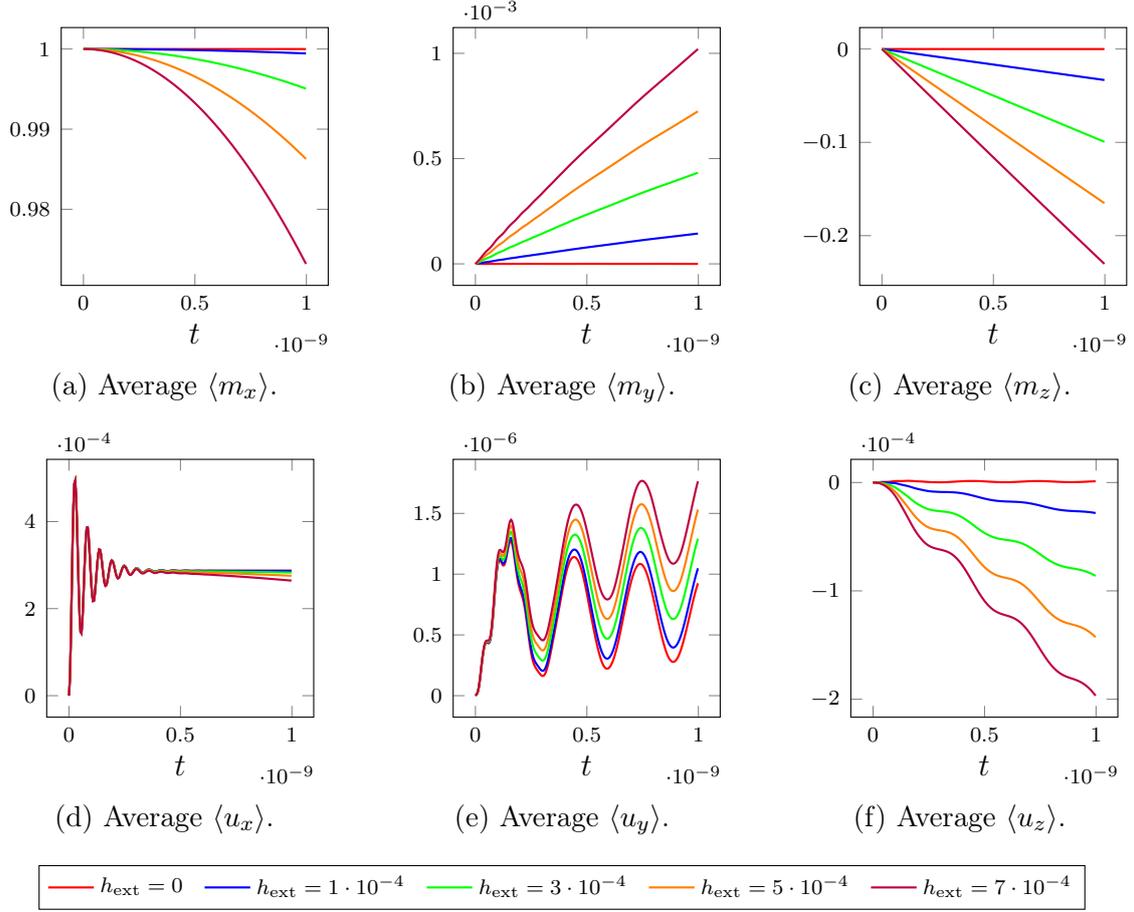
\begin{figure}[ht]
\begin{subfigure}{0.32\textwidth}
\centering
\begin{tikzpicture}
\pgfplotstableread{data/Applied_Field/0em4_zeeman.dat}{\zerozeeman}
\pgfplotstableread{data/Applied_Field/1em4_zeeman.dat}{\onezeeman}
\pgfplotstableread{data/Applied_Field/3em4_zeeman.dat}{\threezeeman}
\pgfplotstableread{data/Applied_Field/5em4_zeeman.dat}{\fivezeeman}
\pgfplotstableread{data/Applied_Field/7em4_zeeman.dat}{\sevenzeeman}
\begin{axis}
[
width = \textwidth, height=5cm,
xlabel = {$t$ (in \si{s})},
legend style={at={(0.5,1)},anchor=south},
legend columns=-1
]
\addplot[red, thick, solid]	table[x=t, y=x_mag_avg]{\zerozeeman};
\addplot[blue, thick, dashed]	table[x=t, y=x_mag_avg]{\onezeeman};
\addplot[green, thick, dotted]	table[x=t, y=x_mag_avg]{\threezeeman};
\addplot[orange, thick, dashdotted]	table[x=t, y=x_mag_avg]{\fivezeeman};
\addplot[purple, thick, densely dotted]	table[x=t, y=x_mag_avg]{\sevenzeeman};
\end{axis}
\end{tikzpicture}
\caption{Average $\langle m_x \rangle$.\label{zeeman_avg_x_mag}}
\end{subfigure}
\begin{subfigure}{0.32\textwidth}
\centering
\begin{tikzpicture}
\pgfplotstableread{data/Applied_Field/0em4_zeeman.dat}{\zerozeeman}
\pgfplotstableread{data/Applied_Field/1em4_zeeman.dat}{\onezeeman}
\pgfplotstableread{data/Applied_Field/3em4_zeeman.dat}{\threezeeman}
\pgfplotstableread{data/Applied_Field/5em4_zeeman.dat}{\fivezeeman}
\pgfplotstableread{data/Applied_Field/7em4_zeeman.dat}{\sevenzeeman}
\begin{axis}
[
width = \textwidth, height=5cm,
xlabel = {$t$ (in \si{s})},
legend style={at={(0.5, 1.1)},anchor=south},
]
\addplot[red, thick, solid]	table[x=t, y=y_mag_avg]{\zerozeeman};
\addplot[blue, thick, dashed]	table[x=t, y=y_mag_avg]{\onezeeman};
\addplot[green, thick, dotted]	table[x=t, y=y_mag_avg]{\threezeeman};
\addplot[orange, thick, dashdotted]	table[x=t, y=y_mag_avg]{\fivezeeman};
\addplot[purple, thick, densely dotted]	table[x=t, y=y_mag_avg]{\sevenzeeman};
\end{axis}
\end{tikzpicture}
\caption{Average $\langle m_y \rangle$.\label{zeeman_avg_y_mag}}
\end{subfigure}
\begin{subfigure}{0.32\textwidth}
\centering
\begin{tikzpicture}
\pgfplotstableread{data/Applied_Field/0em4_zeeman.dat}{\zerozeeman}
\pgfplotstableread{data/Applied_Field/1em4_zeeman.dat}{\onezeeman}
\pgfplotstableread{data/Applied_Field/3em4_zeeman.dat}{\threezeeman}
\pgfplotstableread{data/Applied_Field/5em4_zeeman.dat}{\fivezeeman}
\pgfplotstableread{data/Applied_Field/7em4_zeeman.dat}{\sevenzeeman}
\begin{axis}
[
width = \textwidth, height=5cm,
xlabel = {$t$ (in \si{s})},
]
\addplot[red, thick, solid]	table[x=t, y=z_mag_avg]{\zerozeeman};
\addplot[blue, thick, dashed]	table[x=t, y=z_mag_avg]{\onezeeman};
\addplot[green, thick, dotted]	table[x=t, y=z_mag_avg]{\threezeeman};
\addplot[orange, thick, dashdotted]	table[x=t, y=z_mag_avg]{\fivezeeman};
\addplot[purple, thick, densely dotted]	table[x=t, y=z_mag_avg]{\sevenzeeman};
\end{axis}
\end{tikzpicture}
\caption{Average $\langle m_z \rangle$.\label{zeeman_avg_z_mag}}
\end{subfigure}\\
\bigskip
\begin{subfigure}{0.32\textwidth}
\centering
\begin{tikzpicture}
\pgfplotstableread{data/Applied_Field/0em4_zeeman.dat}{\zerozeeman}
\pgfplotstableread{data/Applied_Field/1em4_zeeman.dat}{\onezeeman}
\pgfplotstableread{data/Applied_Field/3em4_zeeman.dat}{\threezeeman}
\pgfplotstableread{data/Applied_Field/5em4_zeeman.dat}{\fivezeeman}
\pgfplotstableread{data/Applied_Field/7em4_zeeman.dat}{\sevenzeeman}
\begin{axis}
[
width = \textwidth, height=5cm,
xlabel = {$t$ (in \si{s})},
]
\addplot[red, thick, solid]	table[x=t, y=x_disp_avg]{\zerozeeman};
\addplot[blue, thick, dashed]	table[x=t, y=x_disp_avg]{\onezeeman};
\addplot[green, thick, dotted]	table[x=t, y=x_disp_avg]{\threezeeman};
\addplot[orange, thick, dashdotted]	table[x=t, y=x_disp_avg]{\fivezeeman};
\addplot[purple, thick, densely dotted]	table[x=t, y=x_disp_avg]{\sevenzeeman};
\end{axis}
\end{tikzpicture}
\caption{Average $\langle u_x \rangle$.\label{zeeman_avg_x_disp}}
\end{subfigure}
\begin{subfigure}{0.32\textwidth}
\centering
\begin{tikzpicture}
\pgfplotstableread{data/Applied_Field/0em4_zeeman.dat}{\zerozeeman}
\pgfplotstableread{data/Applied_Field/1em4_zeeman.dat}{\onezeeman}
\pgfplotstableread{data/Applied_Field/3em4_zeeman.dat}{\threezeeman}
\pgfplotstableread{data/Applied_Field/5em4_zeeman.dat}{\fivezeeman}
\pgfplotstableread{data/Applied_Field/7em4_zeeman.dat}{\sevenzeeman}
\begin{axis}
[
width = \textwidth, height=5cm,
xlabel = {$t$ (in \si{s})},
]
\addplot[red, thick, solid]	table[x=t, y=y_disp_avg]{\zerozeeman};
\addplot[blue, thick, dashed]	table[x=t, y=y_disp_avg]{\onezeeman};
\addplot[green, thick, dotted]	table[x=t, y=y_disp_avg]{\threezeeman};
\addplot[orange, thick, dashdotted]	table[x=t, y=y_disp_avg]{\fivezeeman};
\addplot[purple, thick, densely dotted]	table[x=t, y=y_disp_avg]{\sevenzeeman};
\end{axis}
\end{tikzpicture}
\caption{Average $\langle u_y \rangle$.\label{zeeman_avg_y_disp}}
\end{subfigure}
\begin{subfigure}{0.32\textwidth}
\centering
\begin{tikzpicture}
\pgfplotstableread{data/Applied_Field/0em4_zeeman.dat}{\zerozeeman}
\pgfplotstableread{data/Applied_Field/1em4_zeeman.dat}{\onezeeman}
\pgfplotstableread{data/Applied_Field/3em4_zeeman.dat}{\threezeeman}
\pgfplotstableread{data/Applied_Field/5em4_zeeman.dat}{\fivezeeman}
\pgfplotstableread{data/Applied_Field/7em4_zeeman.dat}{\sevenzeeman}
\begin{axis}
[
width = \textwidth, height=5cm,
xlabel = {$t$ (in \si{s})},
legend style = {fill=none},
]
\addplot[red, thick, solid]	table[x=t, y=z_disp_avg]{\zerozeeman};
\addplot[blue, thick, dashed]	table[x=t, y=z_disp_avg]{\onezeeman};
\addplot[green, thick, dotted]	table[x=t, y=z_disp_avg]{\threezeeman};
\addplot[orange, thick, dashdotted]	table[x=t, y=z_disp_avg]{\fivezeeman};
\addplot[purple, thick, densely dotted]	table[x=t, y=z_disp_avg]{\sevenzeeman};
\end{axis}
\end{tikzpicture}
\caption{Average $\langle u_z \rangle$.\label{zeeman_avg_z_disp}}
\end{subfigure}\\
\bigskip
\begin{tikzpicture}
    \begin{customlegend}[legend entries={
    $h_{\text{ext}}=0$\hspace*{2mm},
    $h_{\text{ext}}=1\cdot10^{-4}$\hspace*{2mm},
    $h_{\text{ext}}=3\cdot10^{-4}$\hspace*{2mm},
    $h_{\text{ext}}=5\cdot10^{-4}$\hspace*{2mm},
    $h_{\text{ext}}=7\cdot10^{-4}$,},
    legend columns = -1,]
    \addlegendimage{red, thick, solid}
    \addlegendimage{blue, thick, dashed}
    \addlegendimage{green, thick, dotted}
    \addlegendimage{orange, thick, dashdotted}
    \addlegendimage{purple, thick, densely dotted}
    \end{customlegend}
\end{tikzpicture}
\caption{Experiment of Section~\ref{sec:applied_magnetic_field}:
Time evolution of the average magnetisation
and displacement components for varied applied magnetic fields.}
\label{fig:appliedfieldfigure}
\end{figure}

We observe the magnetisation aligning with the applied external
field as expected through a precession, yielding an effect on the displacement.
The coupling is clearly visible in Figure~\ref{fig:appliedfieldfigure},
where we plot the time evolution of the average magnetisation and displacement components,
e.g.\ $\langle u_{x} \rangle = (1/|\Omega|)\int_{\Omega}u_{x}$.
The applied field is pointing in the $y$-direction, so the $y$ and $z$ components begin to increase in magnitude as seen in Figures \ref{zeeman_avg_y_mag} and \ref{zeeman_avg_z_mag}, taking from the $x$ component. The displacement on the other hand mirrors the magnetisation in the $y$ and $z$ components, with the $x$ component increasing due to magnetostriction, and then changing slowly as the magnetisation changes.
Moreover, we see that,
with stronger applied magnetic fields,
the average magnetisation in the
$y$ direction increases, displacing
the body in the same direction.

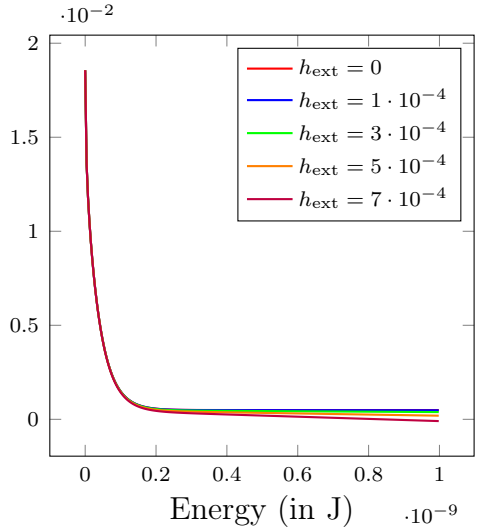
\begin{figure}[ht]
\centering
\begin{tikzpicture}
\pgfplotstableread{data/Applied_Field/0em4_zeeman.dat}{\zerozeeman}
\pgfplotstableread{data/Applied_Field/1em4_zeeman.dat}{\onezeeman}
\pgfplotstableread{data/Applied_Field/3em4_zeeman.dat}{\threezeeman}
\pgfplotstableread{data/Applied_Field/5em4_zeeman.dat}{\fivezeeman}
\pgfplotstableread{data/Applied_Field/7em4_zeeman.dat}{\sevenzeeman}
\begin{axis}
[
width = 0.45\textwidth,
height=0.45\textwidth,
xlabel = {$t$ (in \si{s})},
ylabel = {Energy},
legend pos= north east,
legend cell align= left,
legend style = {fill=none},
]
\addplot[red, thick, solid]	table[x=t, y=totalenergy]{\zerozeeman};
\addplot[blue, thick, dashed]	table[x=t, y=totalenergy]{\onezeeman};
\addplot[green, thick, dotted]	table[x=t, y=totalenergy]{\threezeeman};
\addplot[orange, thick, dashdotted]	table[x=t, y=totalenergy]{\fivezeeman};
\addplot[purple, thick, densely dotted]	table[x=t, y=totalenergy]{\sevenzeeman};
\legend{
$h_{\text{ext}}=0$,
$h_{\text{ext}}=1\cdot10^{-4}$,
$h_{\text{ext}}=3\cdot10^{-4}$,
$h_{\text{ext}}=5\cdot10^{-4}$,
$h_{\text{ext}}=7\cdot10^{-4}$,
}
\end{axis}
\end{tikzpicture}
\caption{Experiment of Section~\ref{sec:applied_magnetic_field}:
Total energy over time for varied applied magnetic fields.}
\label{fig:totalenergyzeeman}
\end{figure}

In Figure~\ref{fig:totalenergyzeeman}, we plot the time evolution of the energy
for all considered applied external fields.
For greater strength applied fields, the energy reaches a lower value at later times.
Importantly, we always see the energy decreasing.

\subsubsection{Inverse magnetostrictive effect}\label{sec:applied_traction}

In this experiment,
we show that changes in the mechanical state of the body
yield changes in the magnetisation.
To this end,
we disable the Zeeman field and apply a traction
on the $x=20\ellex$ plane in the $+y$ direction.
Specifically, we consider a surface force of the form $\boldsymbol{g} = (0,b,0)$ 
for $b \in \{0,1.28\cdot10^{-9},3.19\cdot10^{-9},6.38\cdot10^{-9},1.28\cdot10^{-8}\}$,
which corresponds to forces of strength $0$, $10$, $25$, $50$, $100$ \si{N.m^{-2}}.

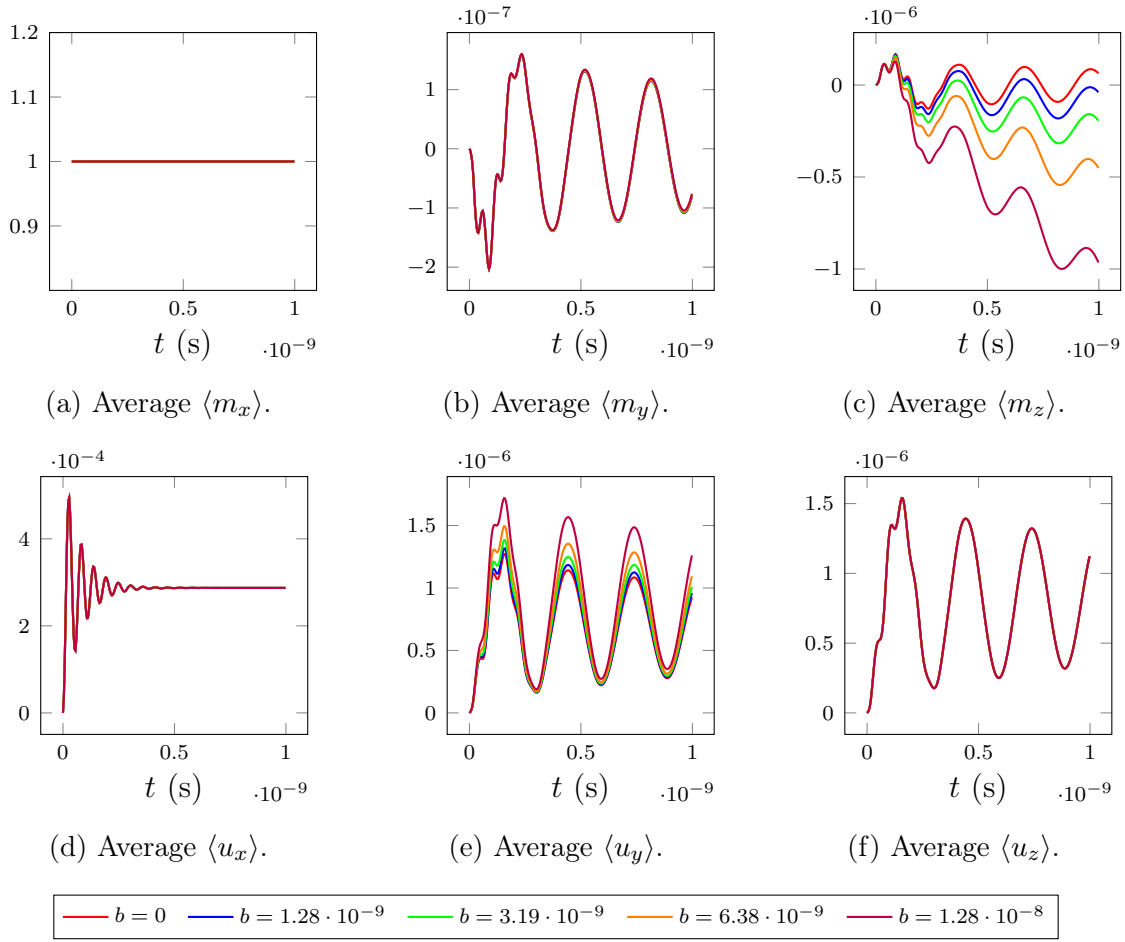
\begin{figure}[ht]
\begin{subfigure}{0.32\textwidth}
\centering
\begin{tikzpicture}
\pgfplotstableread{data/Traction_Test/zero_newton_per_m2.dat}{\zerotraction}
\pgfplotstableread{data/Traction_Test/ten_newton_per_m2.dat}{\tentraction}
\pgfplotstableread{data/Traction_Test/twentyfive_newton_per_m2.dat}{\twentyfivetraction}
\pgfplotstableread{data/Traction_Test/fifty_newton_per_m2.dat}{\fiftytraction}
\pgfplotstableread{data/Traction_Test/hundred_newton_per_m2.dat}{\hundredtraction}
\begin{axis}
[
width = \textwidth, height=5cm,
xlabel = {$t$ (\si{s})},
]
\addplot[red, thick, solid]	table[x=t, y=x_mag_avg]{\zerotraction};
\addplot[blue, thick, dashed]	table[x=t, y=x_mag_avg]{\tentraction};
\addplot[green, thick, dotted]	table[x=t, y=x_mag_avg]{\twentyfivetraction};
\addplot[orange, thick, dashdotted]	table[x=t, y=x_mag_avg]{\fiftytraction};
\addplot[purple, thick, densely dotted]	table[x=t, y=x_mag_avg]{\hundredtraction};
\end{axis}
\end{tikzpicture}
\caption{Average $\langle m_x \rangle$.\label{traction_avg_x_mag}}
\end{subfigure}
\begin{subfigure}{0.32\textwidth}
\centering
\begin{tikzpicture}
\pgfplotstableread{data/Traction_Test/zero_newton_per_m2.dat}{\zerotraction}
\pgfplotstableread{data/Traction_Test/ten_newton_per_m2.dat}{\tentraction}
\pgfplotstableread{data/Traction_Test/twentyfive_newton_per_m2.dat}{\twentyfivetraction}
\pgfplotstableread{data/Traction_Test/fifty_newton_per_m2.dat}{\fiftytraction}
\pgfplotstableread{data/Traction_Test/hundred_newton_per_m2.dat}{\hundredtraction}
\begin{axis}
[
width = \textwidth, height=5cm,
xlabel = {$t$ (\si{s})},
legend pos= north west,
legend cell align= left,
legend style = {fill=none},
]
\addplot[red, thick, solid]	table[x=t, y=y_mag_avg]{\zerotraction};
\addplot[blue, thick, dashed]	table[x=t, y=y_mag_avg]{\tentraction};
\addplot[green, thick, dotted]	table[x=t, y=y_mag_avg]{\twentyfivetraction};
\addplot[orange, thick, dashdotted]	table[x=t, y=y_mag_avg]{\fiftytraction};
\addplot[purple, thick, densely dotted]	table[x=t, y=y_mag_avg]{\hundredtraction};
\end{axis}
\end{tikzpicture}
\caption{Average $\langle m_y \rangle$.\label{traction_avg_y_mag}}
\end{subfigure}
\begin{subfigure}{0.32\textwidth}
\centering
\begin{tikzpicture}
\pgfplotstableread{data/Traction_Test/zero_newton_per_m2.dat}{\zerotraction}
\pgfplotstableread{data/Traction_Test/ten_newton_per_m2.dat}{\tentraction}
\pgfplotstableread{data/Traction_Test/twentyfive_newton_per_m2.dat}{\twentyfivetraction}
\pgfplotstableread{data/Traction_Test/fifty_newton_per_m2.dat}{\fiftytraction}
\pgfplotstableread{data/Traction_Test/hundred_newton_per_m2.dat}{\hundredtraction}
\begin{axis}
[
width = \textwidth, height=5cm,
xlabel = {$t$ (\si{s})},
legend pos= south west,
legend cell align= left,
legend style = {fill=none},
]
\addplot[red, thick, solid]	table[x=t, y=z_mag_avg]{\zerotraction};
\addplot[blue, thick, dashed]	table[x=t, y=z_mag_avg]{\tentraction};
\addplot[green, thick, dotted]	table[x=t, y=z_mag_avg]{\twentyfivetraction};
\addplot[orange, thick, dashdotted]	table[x=t, y=z_mag_avg]{\fiftytraction};
\addplot[purple, thick, densely dotted]	table[x=t, y=z_mag_avg]{\hundredtraction};
\end{axis}
\end{tikzpicture}
\caption{Average $\langle m_z \rangle$.\label{traction_avg_z_mag}}
\end{subfigure}\\
\bigskip
\begin{subfigure}{0.32\textwidth}
\centering
\begin{tikzpicture}
\pgfplotstableread{data/Traction_Test/zero_newton_per_m2.dat}{\zerotraction}
\pgfplotstableread{data/Traction_Test/ten_newton_per_m2.dat}{\tentraction}
\pgfplotstableread{data/Traction_Test/twentyfive_newton_per_m2.dat}{\twentyfivetraction}
\pgfplotstableread{data/Traction_Test/fifty_newton_per_m2.dat}{\fiftytraction}
\pgfplotstableread{data/Traction_Test/hundred_newton_per_m2.dat}{\hundredtraction}
\begin{axis}
[
width = \textwidth, height=5cm,
xlabel = {$t$ (\si{s})},
legend pos= north east,
legend cell align= left,
legend style = {fill=none},
]
\addplot[red, thick, solid]	table[x=t, y=x_disp_avg]{\zerotraction};
\addplot[blue, thick, dashed]	table[x=t, y=x_disp_avg]{\tentraction};
\addplot[green, thick, dotted]	table[x=t, y=x_disp_avg]{\twentyfivetraction};
\addplot[orange, thick, dashdotted]	table[x=t, y=x_disp_avg]{\fiftytraction};
\addplot[purple, thick, densely dotted]	table[x=t, y=x_disp_avg]{\hundredtraction};
\end{axis}
\end{tikzpicture}
\caption{Average $\langle u_x \rangle$.\label{traction_avg_x_disp}}
\end{subfigure}
\begin{subfigure}{0.32\textwidth}
\centering
\begin{tikzpicture}
\pgfplotstableread{data/Traction_Test/zero_newton_per_m2.dat}{\zerotraction}
\pgfplotstableread{data/Traction_Test/ten_newton_per_m2.dat}{\tentraction}
\pgfplotstableread{data/Traction_Test/twentyfive_newton_per_m2.dat}{\twentyfivetraction}
\pgfplotstableread{data/Traction_Test/fifty_newton_per_m2.dat}{\fiftytraction}
\pgfplotstableread{data/Traction_Test/hundred_newton_per_m2.dat}{\hundredtraction}
\begin{axis}
[
width = \textwidth, height=5cm,
xlabel = {$t$ (\si{s})},
legend pos= north west,
legend cell align= left,
legend style = {fill=none},
]
\addplot[red, thick, solid]	table[x=t, y=y_disp_avg]{\zerotraction};
\addplot[blue, thick, dashed]	table[x=t, y=y_disp_avg]{\tentraction};
\addplot[green, thick, dotted]	table[x=t, y=y_disp_avg]{\twentyfivetraction};
\addplot[orange, thick, dashdotted]	table[x=t, y=y_disp_avg]{\fiftytraction};
\addplot[purple, thick, densely dotted]	table[x=t, y=y_disp_avg]{\hundredtraction};
\end{axis}
\end{tikzpicture}
\caption{Average $\langle u_y \rangle$.\label{traction_avg_y_disp}}
\end{subfigure}
\begin{subfigure}{0.32\textwidth}
\centering
\begin{tikzpicture}
\pgfplotstableread{data/Traction_Test/zero_newton_per_m2.dat}{\zerotraction}
\pgfplotstableread{data/Traction_Test/ten_newton_per_m2.dat}{\tentraction}
\pgfplotstableread{data/Traction_Test/twentyfive_newton_per_m2.dat}{\twentyfivetraction}
\pgfplotstableread{data/Traction_Test/fifty_newton_per_m2.dat}{\fiftytraction}
\pgfplotstableread{data/Traction_Test/hundred_newton_per_m2.dat}{\hundredtraction}
\begin{axis}
[
width = \textwidth, height=5cm,
xlabel = {$t$ (\si{s})},
legend pos= south west,
legend cell align= left,
legend style = {fill=none},
]
\addplot[red, thick, solid]	table[x=t, y=z_disp_avg]{\zerotraction};
\addplot[blue, thick, dashed]	table[x=t, y=z_disp_avg]{\tentraction};
\addplot[green, thick, dotted]	table[x=t, y=z_disp_avg]{\twentyfivetraction};
\addplot[orange, thick, dashdotted]	table[x=t, y=z_disp_avg]{\fiftytraction};
\addplot[purple, thick, densely dotted]	table[x=t, y=z_disp_avg]{\hundredtraction};
\end{axis}
\end{tikzpicture}
\caption{Average $\langle u_z \rangle$.\label{traction_avg_z_disp}}
\end{subfigure}\\
\bigskip
\begin{tikzpicture}
    \begin{customlegend}[legend entries={
    $b=0$\hspace*{2mm},
    $b=1.28\cdot10^{-9}$\hspace*{2mm},
    $b=3.19\cdot10^{-9}$\hspace*{2mm},
    $b=6.38\cdot10^{-9}$\hspace*{2mm},
    $b=1.28\cdot10^{-8}$,},
    legend columns = -1,]
    \addlegendimage{red, thick, solid}
    \addlegendimage{blue, thick, dashed}
    \addlegendimage{green, thick, dotted}
    \addlegendimage{orange, thick, dashdotted}
    \addlegendimage{purple, thick, densely dotted}
    \end{customlegend}
\end{tikzpicture}
\caption{
Experiment of Section~\ref{sec:applied_traction}:
Time evolution of the average magnetisation
and displacement components for varied traction strengths.}
\label{fig:tractionfigure}
\end{figure}

The time evolution of the average displacement and magnetisation components is shown in Figure \ref{fig:tractionfigure}.
When more traction is applied, the
average displacement in the $y$ direction
increases. The $z$ component of the magnetisation
in Figure~\ref{traction_avg_z_mag} is the most
interesting, as it decreases more strongly due to
stronger tractions.

\subsubsection{Nutation dynamics}\label{sec:Nutation}
\comment{It has been shown that at extremely short timescales, the
LLG equation is inadequate for the ultrafast dynamics that occur \cite{ciornei2011magnetization,wegrowe2012magnetization}.
At these timescales, the LLG equation~\eqref{eq:llg} should then be replaced by the inertial LLG (iLLG) equation,
given by
\begin{equation}
	\mmt = -\mm\times\heff[\mm] + \alpha \, \mm\times\mmt + \tau \, \mm\times\partial_{tt}\mm,
\end{equation}
where the additional parameter $\tau>0$ is a relaxation time.

The main difference between iLLG and LLG dynamics is the inclusion of
nutation, where the LLG path is not instantly followed
due to the inertia of the magnetisation~\cite{inertia2020}.
A preliminary form of the iLLG equation was initially
derived using an expansion coming from the magnetoelastic coupling \cite{suhl1998theory,suhl2007relaxation}.
Here the momentum is stored by the displacement instead of the magnetisation.
With the experiment in this section, we aim to demonstrate this effect.

We consider the same material parameters as in the previous experiments,
except for the Gilbert damping parameter (for which we choose $\alpha=0.1$).
Moreover, as the nutation effects are small,
we increase the magnetostriction constant $\lambda_\text{100}$ in Table~\ref{tab:valuetable}.
For the ferromagnetic body, we consider a hemisphere of radius $\ellex$
with a clamped planar face with outer normal $(-1,0,0)$.
We use the values $\{20\lambda_{100},50\lambda_{100},100\lambda_{100}\}$,
and compare this to a reference LLG simulation with no
magnetoelastic coupling (computed by the same algorithm with $0\lambda_{100}$).
The initial magnetisation is slightly perturbed from the $x$-direction, specifically
$\mm^0 = (0.9,0.2,0)$ (normalised), and subject to a strong Zeeman field $\hh_{\mathrm{ext}}=(1,0,0) = (3\pi/5,0,0)\si{\tesla}$.
The initial displacement and velocity are zero.
We set $\theta = 0.50000005$, and consider a non-dimensional time-step $k=0.001\approx 3\cdot10^{-15}\si{s}$, for time $T=1\cdot10^{-10}\si{s}$.
The mesh is made of 210 nodes, 672 elements, and
satisfies $h_{\max} \approx 0.6 \ellex$.

\begin{figure}[ht]
	\centering
	\begin{tikzpicture}
		\pgfplotstableread[col sep=comma]{data/Nutation/zerolambda.dat}{\zerolambda}
		\pgfplotstableread[col sep=comma]{data/Nutation/twentylambda.dat}{\twentylambda}
		\pgfplotstableread[col sep=comma]{data/Nutation/fiftylambda.dat}{\fiftylambda}
		\pgfplotstableread[col sep=comma]{data/Nutation/hundredlambda.dat}{\hundredlambda}
		\begin{axis}
			[
			width = 0.45\textwidth,
			height=0.45\textwidth,
			xlabel = {$t$ (in \si{s})},
			legend pos= north east,
			legend cell align= left,
			legend style = {fill=none},
			]
			\addplot[red, thick, dashed]	table[x=t, y=y_mag_avg]{\zerolambda};
			\addplot[blue, thick, dashdotted]	table[x=t, y=y_mag_avg]{\twentylambda};
			\addplot[purple, thick, densely dotted]	table[x=t, y=y_mag_avg]{\fiftylambda};
			\addplot[orange, thick, solid]	table[x=t, y=y_mag_avg]{\hundredlambda};
			\legend{
				LLG,
				$20\lambda_{100}$,
				$50\lambda_{100}$,
				$100\lambda_{100}$
			}
		\end{axis}
	\end{tikzpicture}
	\caption{\comment{Experiment of Section~\ref{sec:Nutation}:
		Time evolution of the average $\langle \mm_y \rangle$ for varying magnetostrain values $\lambda_{100}$, compared with purely LLG without magnetoelastic effects.}}
	\label{fig:nutation_y}
\end{figure}

The results can be seen in Figure~\ref{fig:nutation_y}.
It is easily seen that as the magnetostrain parameter is increased, the dissipative effects decrease, seen as a stretching effect to the right.
As expected, we see nutation effects perturbing the natural LLG
precession behaviour, especially for $50\lambda_{100}$ and $100\lambda_{100}$.

\begin{figure}[ht]
	\centering
	\begin{subfigure}{0.45\textwidth}
		\centering
		\begin{tikzpicture}
			\pgfplotstableread[col sep=comma]{data/Nutation/zerolambda.dat}{\zerolambda}
			\pgfplotstableread[col sep=comma]{data/Nutation/iLLG.dat}{\iLLG}
			\begin{axis}
				[
				width = 1\textwidth,
				height=1\textwidth,
				xlabel = {$t$ (in \si{s})},
				legend pos= south east,
				legend cell align= left,
				legend style = {fill=none},
				]
				\addplot[red, thick, dashed]	table[x=t, y=x_mag_avg]{\zerolambda};
				\addplot[black, thick,densely dotted]	table[x=t, y=x_mag_avg]{\iLLG};
				\legend{
					LLG,
					iLLG
				}
			\end{axis}
		\end{tikzpicture}
		\caption{Average $\langle \mm_x\rangle$ for reference iLLG and LLG.\label{fig:iLLGvsMag1}}
	\end{subfigure}%
	\hfill
	\begin{subfigure}{0.45\textwidth}
		\centering
		\begin{tikzpicture}
			\pgfplotstableread[col sep=comma]{data/Nutation/zerolambda.dat}{\zerolambda}
			\pgfplotstableread[col sep=comma]{data/Nutation/twentylambda.dat}{\twentylambda}
			\pgfplotstableread[col sep=comma]{data/Nutation/fiftylambda.dat}{\fiftylambda}
			\pgfplotstableread[col sep=comma]{data/Nutation/hundredlambda.dat}{\hundredlambda}
			\begin{axis}
				[
				width = 1\textwidth,
				height=1\textwidth,
				xlabel = {$t$ (in \si{s})},
				legend pos= south east,
				legend cell align= left,
				legend style = {fill=none},
				]
				\addplot[red, thick, dashed]	table[x=t, y=x_mag_avg]{\zerolambda};
				\addplot[blue, thick, dashdotted]	table[x=t, y=x_mag_avg]{\twentylambda};
				\addplot[purple, thick, densely dotted]	table[x=t, y=x_mag_avg]{\fiftylambda};
				\addplot[orange, thick, solid]	table[x=t, y=x_mag_avg]{\hundredlambda};
				\legend{
					LLG,
					$20\lambda_{100}$,
					$50\lambda_{100}$,
					$100\lambda_{100}$
				}
			\end{axis}
		\end{tikzpicture}
		\caption{Average $\langle \mm_x\rangle$ for magnetoelastic LLG.\label{fig:iLLGvsMag2}}
	\end{subfigure}
	\caption{Experiment of Section~\ref{sec:Nutation}. Average magnetisation $\langle \mm_x\rangle$ of (a) reference LLG and iLLG, (b) magnetoelastic LLG with varying magnetostrain.\label{fig:iLLGvsMag}}
\end{figure}

For reference, we also consider this same system
without magnetoelastic coupling for the iLLG equation.
The modifications required to extend the tangent plane scheme
presented in \cite{alouges2008a} to the iLLG equation
are described and analysed in \cite{ruggeri2022numerical}.
For consistency, we apply the tangent plane scheme here with the nodal projection step removed.
We choose the relaxation time arbitrarily to be $\tau=0.4 \approx 1.21 \si{ps}$.
The average $x$ component of the magnetisation for iLLG is shown in Figure~\ref{fig:iLLGvsMag1},
and the average $x$ components of the magnetisation for the
magnetoelastic LLG simulations is shown
in Figure~\ref{fig:iLLGvsMag2}.
The qualitative similarities are obvious,
with additional oscillations not seen in the pure LLG case,
along with lessened damping.}

\subsection{Properties of Algorithm~\ref{algorithm}}

In this section, we present three experiments to numerically investigate the properties of Algorithm~\ref{algorithm}.
For all of them, the computational domain will be a cube with edge length equal to $6\ellex$.

\subsubsection{$\theta$-dependence} \label{sec:theta_dependence}

In this experiment,
we investigate the effect on numerical simulations of the parameter $\theta \in (1/2,1]$,
which controls the `degree of implicitness' in the treatment of the exchange contribution in~\eqref{alg1:magnetisation_update}.
We use material parameters for FeCoSiB (cf.\ Table~\ref{tab:valuetable}) except for the Gilbert damping parameter,
for which we use the smaller value $\alpha=0.001$.
The initial condition for the magnetisation is a `hot' magnetic state, i.e.\
the values at the vertices of the mesh (which in this experiment has mesh size $h_{\max}\approx 3 \ellex$) are assigned randomly
to the magnetisation before being normalised.
The displacement and its time derivative are initialised by zero.

We run the simulation for $1\cdot 10^{-11}$ \si{s} using a time step size of $1\cdot 10^{-15}$ \si{s}
and different values of $\theta \in \{ 0.50000005 , 0.505 , 0.6 , 0.7 , 0.8 , 0.9 , 1 \}$.

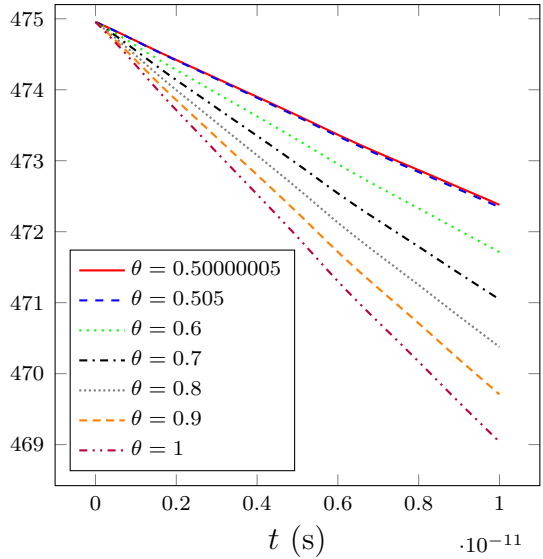
\begin{figure}[ht]
\centering
\begin{tikzpicture}
\pgfplotstableread{data/theta_variance.dat}{\thetaalmosthalf}
\pgfplotstableread{data/theta_variance.dat}{\thetanearlyhalf}
\pgfplotstableread{data/theta_variance.dat}{\thetasix}
\pgfplotstableread{data/theta_variance.dat}{\thetaseven}
\pgfplotstableread{data/theta_variance.dat}{\thetaeight}
\pgfplotstableread{data/theta_variance.dat}{\thetanine}
\pgfplotstableread{data/theta_variance.dat}{\thetaten}
\begin{axis}
[
width = 0.5\textwidth, height=0.5\textwidth,
xlabel = {$t$ (in \si{s})},
legend pos= south west,
legend cell align= left,
legend style = {fill=none},
ymax = 475.2
]
\addplot[red, thick, solid]	table[x=t, y=theta50000005]{\thetaalmosthalf};
\addplot[blue, thick, dashed]	table[x=t, y=theta505]{\thetanearlyhalf};
\addplot[green, thick, dotted]	table[x=t, y=theta6]{\thetasix};
\addplot[black, thick, dashdotted]	table[x=t, y=theta7]{\thetaseven};
\addplot[gray, thick, densely dotted]	table[x=t, y=theta8]{\thetaeight};
\addplot[orange, thick, densely dashed]	table[x=t, y=theta9]{\thetanine};
\addplot[purple, thick, dashdotdotted]	table[x=t, y=theta10]{\thetaten};
\legend{$\theta=0.50000005$,
$\theta=0.505$,
$\theta=0.6$,
$\theta=0.7$,
$\theta=0.8$,
$\theta=0.9$,
$\theta=1$}
\end{axis}
\end{tikzpicture}
\caption{
Experiment of Section~\ref{sec:theta_dependence}:
Time evolution of the total energy for different values of $\theta$.
}
\label{fig:thetatesting}
\end{figure}

The energy-decreasing behaviour can be
seen in Figure~\ref{fig:thetatesting},
with considerably more energy loss associated
with greater $\theta$ values.
So changing the $\theta$-implicitness
parameter away from $1/2$ can yield
considerable amounts of artificial numerical damping,
which can be particularly bad in certain situations
(e.g.\ in the case of long-time simulations).

\subsubsection{Unit length constraint violation} \label{sec:numerics_constraint}

An essential property of the LLG equation at constant temperature 
is the unit length constraint on the magnetisation.
Hence, an essential feature of any approximation algorithm
must be the capability to achieve the unit length constraint. 
For Algorithm~\ref{algorithm}, this property is the subject of Proposition~\ref{prop:stability},
particularly~\eqref{eq:general_constraint}, i.e.
\begin{equation*}
\big\lVert \II_h\big[\abs{\mm_h^j}^2\big]-1 \big\rVert_{L^1(\Omega)}
\le C k,
\end{equation*}
which shows that the unit length constraint is violated at most linearly in time (if measured in the $L^1$-norm).

To see this numerically, we again consider a hot magnet as in Section~\ref{sec:theta_dependence},
a particularly bad case with plenty of rotation by the magnetisation
(note that the constant $C>0$ in~\eqref{eq:general_constraint} depends, among other things, upon the energy of the initial magnetisation
and is large for a random configuration), and use various time-steps with $\theta = 0.50000005$.

\begin{figure}[ht]
\begin{subfigure}{0.45\textwidth}
\centering
\begin{tikzpicture}
\pgfplotstableread{data/constraint_violation_k_integral_error_nodal_max.dat}{\integralconstraint}
\begin{loglogaxis}
[
ymin=0.1,
ymax=250,
width = \textwidth, height=\textwidth,
xlabel = {$k$},
legend pos= north west,
legend cell align= left,
legend style = {fill=none},
]
\addplot[red, thick, mark=*]	table[x=k, y=integral_error]{\integralconstraint};
\legend{$\lVert \II_h[\abs{\mm_h^N}^2]-1 \rVert_{L^1(\Omega)}$}
\end{loglogaxis}
\end{tikzpicture}
\caption{Constraint violation.\label{fig:ConstraintViolation}}
\end{subfigure}
\begin{subfigure}{0.45\textwidth}
\centering
\begin{tikzpicture}
\pgfplotstableread{data/constraint_violation_k_integral_error_nodal_max.dat}{\nodalmax}
\begin{loglogaxis}
[
ymin=0.95,
width = \textwidth, height=\textwidth,
xlabel = {$k$},
legend pos= north west,
legend cell align= left,
legend style = {fill=none},
]
\addplot[red, thick, mark=*]	table[x=k, y=nodalmax]{\nodalmax};
\legend{$\norm[\LL^{\infty}(\Omega)]{\mm_{h}^{N}}$}
\end{loglogaxis}
\end{tikzpicture}
\caption{Nodal maximum.}
\label{fig:Linfty}
\end{subfigure}
\caption{
Experiment of Section~\ref{sec:numerics_constraint}:
(a) Constraint violation at the final iterate against the time-step size.
(b) $L^\infty$-norm of the magnetisation at the final iterate against the time-step size.}
\label{fig:constraintviolationtest}
\end{figure}
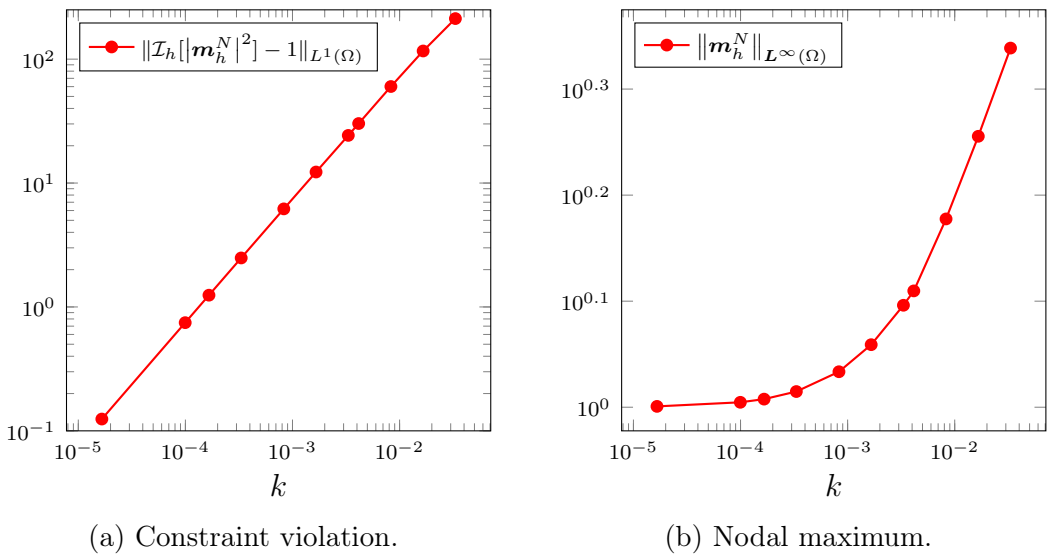

In Figure~\ref{fig:ConstraintViolation}, we plot the constraint violation (measured as the left-hand side of~\eqref{eq:general_constraint})
at the finale iterate $\mm_{h}^{N}$ of Algorithm~\ref{algorithm} against the time-step size $k$.
We observe that the error decays linearly in $k$ as predicted by~\eqref{eq:general_constraint}.
The constraint violation is of the order $10^2$ for $k$ on the order of $10^{-3/2}$ due to the hot initial state,
as the magnetisation at a node may need to rotate several times.
Note that these simulations are run for only $0.01$ \si{ns} as we are only interested in verifying the constraint violation inequalities.

In Figure~\ref{fig:Linfty},
we plot the $L^\infty$-norm of the magnetisation at the final iterate against the time-step size $k$.
We note that while we can control the integral violation with~\eqref{eq:general_constraint}, in our projection-free algorithm
we cannot directly control the maximum norm $\norm[\LL^{\infty}(\Omega)]{\mm_{h}^{j}}$,
which with a projection would be $1$ for each $j$.
We see that the nodal maximum numerically tends to $1$ as desired,
but the decay is not linear.
Using similar methods to those to prove~\eqref{eq:general_constraint}
(see Lemma~\ref{lem:Magnetisation_Estimate} below) and classical inverse estimates \cite[Lemma~3.5]{bartels2015},
one can show that
\begin{equation*}
\begin{split}
    \norm[\LL^{\infty}]{\mm_{h}^{j}}^2 - 1
    &= \max_{z\in\mathcal{N}_{h}}|\mm_{h}^{j}(z)|^2 - 1
    \leq k^2 \sum_{i=0}^{j-1}\max_{z\in\mathcal{N}_{h}}|\vv_{h}^{i}(z)|^2\\
    &= k^2 \sum_{i=0}^{j-1}\norm[\LL^{\infty}(\Omega)]{\vv_{h}^{i}}^2
    \lesssim h_{\min}^{-3}k^2 \sum_{i=0}^{j-1}\norm[\LL^{2}(\Omega)]{\vv_{h}^{i}}^2
    \lesssim h_{\min}^{-3}k,
\end{split}
\end{equation*}
thus
the desired convergence $\norm[\LL^{\infty}(\Omega)]{\mm_{h}^{j}} \to 1$ as $h,k\to 0$
can be obtained assuming the CFL condition $k = o(h^3)$.


\subsubsection{Energy law robustness} \label{sec:numerics_cfl}

In this experiment,
we investigate the robustness of the evolution of the energy
of the approximations generated by Algorithm~\ref{algorithm}
with respect to the discretisation parameters.

We consider a similar setup to the one used in Section~\ref{sec:theta_dependence}.
Specifically,
we keep $\theta = 0.50000005$ and $\alpha=0.001$,
but we add a Zeeman field $\hh_{\mathrm{ext}} = (0.001, 0, 0)\approx(1.9, 0, 0)\si{\milli\tesla}$
to encourage the system to approach the same final state.
To give consistency between mesh refinements we change from a purely random initial state
to the following initial condition for the magnetisation,
\begin{equation*}
    \mm^{0}(x,y,z) = \frac{1}{\sqrt{5}}(2,\sin(x+y+z),\cos(x+y+z))
    \quad \text{for all } (x,y,z)\in\Omega.
\end{equation*}
It is easily shown that the initial condition satisfies 
$\norm{\Grad \mm^{0}}^2/2 = 64.8$ and $\abs{\mm^0}=1$ in $\Omega$.
NGSolve interpolates the initial condition onto the mesh via
an Oswald-type interpolation~\cite{oswald1993}, applying an $L^2$-projection and then averaging for
conformity, thus to enforce the condition
$\mm_{h}^{0}\in\MM_{h,0}$ we apply the nodal projection
to the result of this interpolation.
We then ran the simulation for $T \approx 3.32$
with combinations of
$k=0.01$, $0.005$, $0.0025$, $0.00125$, $0.000625$ as time-step size and
$h = 1.59$, $1.09$, $0.84$, $0.45$ as mesh size.

\begin{figure}[h]
	\centering
	\begin{tikzpicture}
		\begin{axis}
			[
			width = 0.6\textwidth, height=0.5\textwidth,
			xlabel = {$t$ (in \si{s})},
			ymax = 73,
			ymin = 63
			]
			\addplot[red, thick]	table[x=t, y=totalenergy, col sep=comma]{data/CFL_Test/k001h045.dat};
			\label{plot:line1}
			\addplot[red, thick, loosely dotted]	table[x=t, y=totalenergy, col sep=comma]{data/CFL_Test/k001h084.dat};
			\label{plot:line2}
			\addplot[red, thick, dotted]	table[x=t, y=totalenergy, col sep=comma]{data/CFL_Test/k001h109.dat};
			\label{plot:line3}
			\addplot[red, thick, dashed]	table[x=t, y=totalenergy, col sep=comma]{data/CFL_Test/k001h159.dat};
			\label{plot:line4}
			\addplot[blue, thick, loosely dotted]	table[x=t, y=totalenergy, col sep=comma]{data/CFL_Test/k0005h084.dat};
			\label{plot:line5}
			\addplot[blue, thick, dotted]	table[x=t, y=totalenergy, col sep=comma]{data/CFL_Test/k0005h109.dat};
			\label{plot:line6}
			\addplot[blue, thick, dashed]	table[x=t, y=totalenergy, col sep=comma]{data/CFL_Test/k0005h159.dat};
			\label{plot:line7}
			\addplot[brown, thick, loosely dotted]	table[x=t, y=totalenergy, col sep=comma]{data/CFL_Test/k00025h084.dat};
			\label{plot:line8}
			\addplot[brown, thick, dotted]	table[x=t, y=totalenergy, col sep=comma]{data/CFL_Test/k00025h109.dat};
			\label{plot:line9}
			\addplot[brown, thick, dashed]	table[x=t, y=totalenergy, col sep=comma]{data/CFL_Test/k00025h159.dat};
			\label{plot:line10}
			\addplot[cyan, thick, loosely dotted]	table[x=t, y=totalenergy, col sep=comma]{data/CFL_Test/k000125h084.dat};
			\label{plot:line11}
			\addplot[cyan, thick, dotted]	table[x=t, y=totalenergy, col sep=comma]{data/CFL_Test/k000125h109.dat};
			\label{plot:line12}
			\addplot[cyan, thick, dashed]	table[x=t, y=totalenergy, col sep=comma]{data/CFL_Test/k000125h159.dat};
			\label{plot:line13}
			\addplot[violet, thick, loosely dotted]	table[x=t, y=totalenergy, col sep=comma]{data/CFL_Test/k0000625h084.dat};
			\label{plot:line14}
			\addplot[violet, thick, dotted]	table[x=t, y=totalenergy, col sep=comma]{data/CFL_Test/k0000625h109.dat};
			\label{plot:line15}
			\addplot[violet, thick, dashed]	table[x=t, y=totalenergy, col sep=comma]{data/CFL_Test/k0000625h159.dat};
			\label{plot:line16}
			\coordinate (legend) at (axis description cs:0.5,-0.5);
		\end{axis}
		\matrix [
		draw,
		matrix of nodes,
		anchor=center,
		] at (legend) { \scriptsize
			     & \scriptsize $h=1.59$ & \scriptsize $h=1.09$ & \scriptsize $h= 0.84$ & \scriptsize $h= 0.45$ \\
			\scriptsize $k=0.01$ & \ref*{plot:line4} & \ref*{plot:line3} & \ref*{plot:line2} &  \ref*{plot:line1} \\
			\scriptsize $k=0.005$  & \ref*{plot:line7} & \ref*{plot:line6} & \ref*{plot:line5} & \\
			\scriptsize $k=0.0025$ & \ref*{plot:line10} & \ref*{plot:line9} & \ref*{plot:line8} & \\
			\scriptsize $k=0.00125$ & \ref*{plot:line13} & \ref*{plot:line12} & \ref*{plot:line11} & \\
			\scriptsize $k=0.000625$ & \ref*{plot:line16} & \ref*{plot:line15} & \ref*{plot:line14} & \\
		};
	\end{tikzpicture}
	\caption{
		Time evolution of the total energy of different values of $h$ and $k$.}
	\label{fig:CFLGraph}
\end{figure}
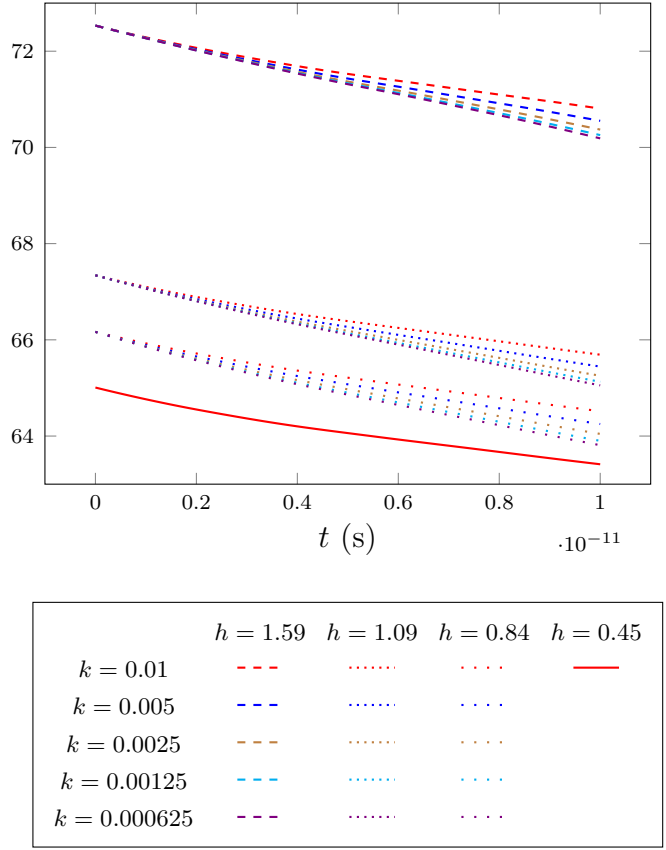

As can be observed in Figure~\ref{fig:CFLGraph},
the energy decay (and thus stability of the algorithm)
occurs for all mesh sizes and time-steps.
The initial energy is
different for each due to the differing underlying mesh,
and the interpolation process mentioned
above (which is also different for each mesh),
however the initial energies approach
the actual energy.
The different energy progressions are clustered
into the four groups with similar
energy decay when the time-step is the same.
When the time-step size is smaller,
the energy decay is slower, likely due to the error term in the discrete energy law
(cf.\ the term $E_{h,k}^i$ in~\eqref{eq:discrete_energy_law}).
With no error term present, the dissipation would always reduce with lower time-step sizes.

These results show that the algorithm behaves energetically well
for all combinations of mesh and time-step size considered,
including the worst case scenario for a CFL condition
(when the finest mesh with $h=0.45$ and the largest time-step size $k=0.01$ are used).
Clearly, this is not a mathematical proof
that the restrictive CFL condition we need to show Theorem~\ref{thm:convergence}(ii)
is not needed, however our numerical experiments seem to corroborate this claim.

\section{Proofs}
\label{sec:proofs}

In this section,
we collect the proofs of the results presented in Section~\ref{sec:main}.
For the convenience of the reader, we start with recalling some well-known results
that will be used multiple times throughout the upcoming analysis.

The norm $\norm[h]{\cdot}$ induced by the mass-lumped $L^2$-product~\eqref{eq:mass-lumping}
satisfies the norm equivalence
\begin{equation} \label{eq:h-scalar-product_equivalence}
\norm{\pphi_h}
\leq \norm[h]{\pphi_h}
\leq \sqrt{5} \, \norm{\pphi_h}
\quad \text{for all } \pphi_h \in \S^1(\T_h)^3,
\end{equation}
and we have the error estimate
\begin{equation} \label{eq:h-scalar-product}
\lvert \inner{\pphi_h}{\ppsi_h} -  \inner[h]{\pphi_h}{\ppsi_h} \rvert
\le C h^2 \norm{\Grad\pphi_h} \norm{\Grad\ppsi_h}
\quad \text{for all } \pphi_h, \ppsi_h \in \S^1(\T_h)^3,
\end{equation}
(cf.\ \cite[Lemma~3.9]{bartels2015}).
For all $K \in \T_h$ and $1 \le r,p \le \infty$,
we have the local inverse estimate
\begin{equation} \label{eq:inverse_estimate}
\norm[\LL^p(K)]{\pphi_h}
\le
C h_K^{3(r-p)/(pr)}
\norm[\LL^r(K)]{\pphi_h}
\quad
\text{for all }
\pphi_h \in \S^1(\T_h)^3
\end{equation}
(see, e.g.\ \cite[Lemma~3.5]{bartels2015}).
For all $1 \le p < \infty$,
the $L^p$-norm of functions in $\S^1(\T_h)^3$
is equivalent with the $\ell^p$-norm of the vector
collecting their nodal values, weighted by the local mesh size,
i.e.
\begin{equation} \label{eq:Lp_equivalence}
C^{-1} \norm[\LL^p(\Omega)]{\pphi_h}
\leq \left( \sum_{z \in \NN_h} h_z^3 \abs{\pphi_h(z)}^p\right)^{1/p}
\leq C \norm[\LL^p(\Omega)]{\pphi_h}
\quad \text{for all } \pphi_h \in \S^1(\T_h)^3,
\end{equation}
where $h_z>0$ denotes the diameter of the node patch of $z \in \NN_h$
(cf.\ \cite[Lemma 3.4]{bartels2015}).
If $p = \infty$, we have that
\begin{equation*}
\norm[\LL^\infty(\Omega)]{\pphi_h}
=
\max_{z \in \NN_h} \abs{\pphi_h(z)}
\quad \text{for all } \pphi_h \in \S^1(\T_h)^3.
\end{equation*}
Finally, the nodal projection is $H^1$-stable,
i.e,
it holds that
\begin{equation} \label{eq:projection_stability}
\norm{\Grad\Pi_h\phi_h} \le C \norm{\Grad\phi_h}
\quad
\text{for all } \pphi \in \S^1(\T_h)^3
\text{ satisfying } 
\abs{\pphi(z)} \ge 1
\text{ for all } z \in \NN_h;
\end{equation}
see~\cite[Lemma~2.2]{bartels2016}.
We recall that~\eqref{eq:projection_stability} holds with $C=1$ if all non-diagonal entries
of the stiffness matrix are non-positive (cf.\ \cite[Proposition~3.2]{bartels2015}).
This assumption, which is satisfied under very restrictive geometric conditions on the mesh in three dimensions,
is \emph{not} required by the upcoming analysis.
In all these inequalities,
the constant $C>0$ (not the same at each occurrence)
depends only on the shape-regularity of $\T_h$.

\subsection{Well-posedness} \label{sec:proofs_wellposedness}

We start by showing an estimate of the
$L^2$-norm of
the discrete elastic field.

\begin{lemma} \label{lem:boundedness_elastic_field}
For all $\uu_h,\mm_h \in \S^1(\T_h)^3$
with $\abs{\mm_h(z)}\ge 1$ for all $z \in \NN_h$,
it holds that
\begin{equation} \label{eq:boundedness_elastic_field}
\norm{\hhm[\uu_{h},\Pi_h\mm_{h}]}^2
\le
8 \norm[\LL^{\infty}(\Omega)]{\Z}^2
\norm[\LL^{\infty}(\Omega)]{\C}^2
\left(\norm{\boldvar(\uu_{h})}^2 + \norm[\LL^{\infty}(\Omega)]{\Z}^2 |\Omega| \right).
\end{equation}
\end{lemma}

\begin{proof}
Using the expression
of the discrete elastic field, we have
\begin{equation*}
\begin{split}
& \norm{\hhm[\uu_{h},\Pi_h\mm_{h}]}^2 \\
& \quad \stackrel{\eqref{eqn:MagnetostrictiveEffectiveField}}{=}
\norm{2 (\Z^{\top}: \{ \C :[\boldvar(\uu_{h}) - \boldvarm(\Pi_h\mm_{h})]\} ) \Pi_{h} \mm_{h}}^2 \\
& \quad \stackrel{\phantom{\eqref{eqn:MagnetostrictiveEffectiveField}}}{\le}
4 \norm[\LL^{\infty}(\Omega)]{\Z}^2
\norm[\LL^{\infty}(\Omega)]{\C}^2
\norm{\boldvar(\uu_{h}) - \Z:(\Pi_h\mm_{h} \otimes \Pi_h\mm_{h})}^2
\norm[\LL^{\infty}(\Omega)]{\Pi_{h} \mm_{h}}^2 \\
& \quad \stackrel{\phantom{\eqref{eqn:MagnetostrictiveEffectiveField}}}{\le}
8 \norm[\LL^{\infty}(\Omega)]{\Z}^2
\norm[\LL^{\infty}(\Omega)]{\C}^2
\left(\norm{\boldvar(\uu_{h})}^2 + \norm[\LL^{\infty}(\Omega)]{\Z}^2 |\Omega| \right),
\end{split}
\end{equation*}
where we have used the boundedness of the fourth-order tensors
and $\Pi_{h} \mm_{h}$.
\end{proof}

We can now show the well-posedness of Algorithm~\ref{algorithm}.

\begin{proof}[Proof of Proposition~\ref{prop:wellposedness}]
The proof is basically identical to the one given in~\cite{bppr2013}
for the algorithm proposed therein.
We restate it here including other terms.

For the magnetisation term, define the family of bilinear form
$a^{i}_{1}(\cdot,\cdot):\KK_{h}[\mm_{h}^{i}]\times \KK_{h}[\mm_{h}^{i}]\to\R$ for $i=0,\dots,N-1$, by
\begin{equation*}
a^{i}_{1}(\pphi_h,\ppsi_h)
:=
\alpha\inner[h]{\pphi_h}{\ppsi_h}
+ \theta  k \inner{\Grad \pphi_h}{\Grad\ppsi_h}
+ \inner{\mm_{h}^{i}\times \pphi_h}{\ppsi_h}
\end{equation*}
and the family of linear (and bounded by Lemma~\ref{lem:boundedness_elastic_field}) functionals $L_{1}^{i}$ for $i=0,\ldots,N-1$ by
\begin{equation*}
L_{1}^{i}(\pphi_h)
:=
-\inner{\Grad\mm_{h}^{i}}{\Grad \pphi_h}
+ \inner{\hhm[\uu_{h}^{i},\Pi_h\mm_{h}^{i}]}{\pphi_h}.
\end{equation*}
Then \eqref{alg1:magnetisation_update} can be rewritten as
$a^{i}_{1}(\vv_h^i,\ppsi_h) = L_{1}^{i}(\ppsi_h)$
for all $\ppsi_h\in\KK_{h}[\mm_{h}^{i}]$.
We can see that $a_{1}^{i}(\cdot,\cdot)$ is positive definite
(in $\LL^2(\Omega)$ and $\HH^1(\Omega)$),
as letting $\pphi_h=\ppsi_h$ eliminates the final term,
leaving a combination of the $L^2$-norm and $H^1$-seminorm.
It follows by the finite dimensionality
that \eqref{alg1:magnetisation_update} has a unique solution
$\vv_{h}^{i}\in \KK_{h}[\mm_{h}^{i}]$.

For the displacement term, define the bilinear form
$a_{2}:\S^{1}_{D}(\T_{h})^3\times \S^{1}_{D}(\T_{h})^3\to\R$ by
\begin{equation*}
a_{2}(\pphi_h,\ppsi_h)
:=
\inner{\pphi_h}{\ppsi_h}
+ k^2\inner{\C:\boldvar(\pphi_h)}{\boldvar(\ppsi_h)}.
\end{equation*}
As $\mathbb{C}$ is positive definite by assumption,
applying Korn's inequality
(see, e.g.\ \cite[~Theorem 11.2.6]{brenner2008mathematical})
yields positive definiteness of $a_{2}(\cdot,\cdot)$ in $\HH^1(\Omega)$.
Furthermore, defining the family of linear functionals
\begin{equation*}
L^{i}_{2}(\ppsi_h)
:=
k^2 \inner{\C:\boldvarm(\Pi_h\mm_{h}^{i+1})}{\boldvar(\ppsi_h)}
+ k\inner{\dt\uu_{h}^{i}}{\ppsi_h}
+ \inner{\uu_{h}^{i}}{\ppsi_h}
+ k^2 \inner{\ff}{\ppsi_h}
+ k^2 \inner[\Gamma_{N}]{\gg}{\ppsi_h},
\end{equation*}
we have that \eqref{alg1:displacement_update}
is equivalent to $a_{2}(\uu_h^{i+1},\ppsi_h) = L^{i}_{2}(\ppsi_h)$ for all $\ppsi_h\in\S^{1}_{D}(\T_{h})^3$, for each $i=0,\ldots,N-1$.
Again exploiting the finite dimension,
we have existence and uniqueness of a solution
$\uu_{h}^{i+1} \in \S^{1}_{D}(\T_{h})^3$ to~\eqref{alg1:displacement_update}.
\end{proof}

\subsection{Discrete energy law} \label{sec:proofs_energy_law}

We now prove the discrete energy law satisfied by the iterates
of Algorithm~\ref{algorithm}.

\begin{proof}[Proof of Proposition~\ref{prop:energy}]
Let $0 \le i \le N-1$ be an arbitrary integer.
Choosing the test function
$\pphi_h = \vv_{h}^{i} \in \KK_h[\mm_{h}^{i}]$ in~\eqref{alg1:magnetisation_update},
we obtain
\begin{equation*}
    \alpha \norm[h]{\vv_{h}^{i}}^2
    + \theta k\norm{\Grad \vv_{h}^{i}}^2
    = - \inner{\Grad\mm_{h}^{i}}{\Grad \vv_{h}^{i}}
    + \inner{\hhm[\uu_{h}^{i},\Pi_h\mm_{h}^{i}]}{\vv_{h}^{i}}.
\end{equation*}
Moreover, we have
\begin{equation*}
\frac{1}{2}\norm{\Grad\mm_{h}^{i+1}}^2
= \frac{1}{2}\norm{\Grad\mm_{h}^{i}}^2
+ k\inner{\Grad \mm_{h}^{i}}{\Grad \vv_{h}^{i}}
+\frac{k^2}{2}\norm{\Grad\vv_{h}^{i}}^2.
\end{equation*}
Combining the two above equations, we obtain
\begin{equation} \label{eqn:MagneticNextStepEnergyRelation}
    \E_{\mathrm{m}}[\mm_{h}^{i+1}]
    - \E_{\mathrm{m}}[\mm_{h}^i]
    = - \alpha k \norm[h]{\vv_{h}^{i}}^2
    - k^2(\theta - 1/2) \norm{\Grad\vv_{h}^{i}}^2
    + k\inner{\hhm[\uu_{h}^{i},\Pi_h\mm_{h}^{i}]}{\vv_{h}^{i}}.
\end{equation}
Choosing the test function
$\ppsi_h = \uu_h^{i+1} - \uu_h^{i} = k \, \dt \uu_h^{i+1}$ in~\eqref{alg1:displacement_update}
yields
\begin{multline*}
    \inner{\dt \uu_{h}^{i+1} - \dt\uu_{h}^{i}}{\dt\uu_h^{i+1}} +  \inner{\mathbb{C}[\boldvar(\uu_{h}^{i+1}) - \boldvarm(\Pi_h\mm_{h}^{i+1})]}{\boldvar(\uu_h^{i+1}) - \boldvar(\uu_h^{i})} \\
    = \inner{\ff}{\uu_h^{i+1} - \uu_h^{i}}
    + \inner[\Gamma_{N}]{\gg}{\uu_h^{i+1} - \uu_h^{i}}.
\end{multline*}
%
%
%
Using Lemma~\ref{lem:abel}, the first term on the left-hand side can be reformulated as
\begin{equation*}
\inner{\dt \uu_{h}^{i+1} - \dt\uu_{h}^{i}}{\dt\uu_h^{i+1}}
=
\frac{1}{2}
\norm{\dt \uu_{h}^{i+1}}^2
-
\frac{1}{2}
\norm{\dt \uu_{h}^i}^2
+
\frac{1}{2}
\norm{\dt \uu_{h}^{i+1} - \dt\uu_{h}^i}^2
\end{equation*}
which yields
\begin{multline}\label{eqn:AfterAbelFirstSummation}
    \frac{1}{2}\norm{\dt \uu_{h}^{i+1}}^2-\frac{1}{2}\norm{\dt \uu_{h}^i}^2+\frac{1}{2}\norm{\dt \uu_{h}^{i+1} - \dt\uu_{h}^i}^2 \\
    + \inner{\mathbb{C}[\boldvar(\uu_{h}^{i+1}) - \boldvarm(\Pi_h\mm_{h}^{i+1})]}{\boldvar(\uu_h^{i+1}) - \boldvar(\uu_h^{i})} \\
    = \inner{\ff}{\uu_h^{i+1} - \uu_h^{i}}
    + \inner[\Gamma_{N}]{\gg}{\uu_h^{i+1} - \uu_h^{i}}.
\end{multline}
Similarly, we have
\begin{equation*}
\begin{split}
&\inner{\C:[\boldvar(\uu_{h}^{i+1}) - \boldvarm(\Pi_h\mm_{h}^{i+1})]}{\boldvar(\uu_h^{i+1}) - \boldvar(\uu_h^{i})} \\
& \quad
\stackrel{\phantom{\eqref{eq:abel}}}{=}
\inner{\C:[\boldvar(\uu_{h}^{i+1}) - \boldvarm(\mm_{h}^{i+1})]}{\boldvar(\uu_h^{i+1}) - \boldvar(\uu_h^{i})}\\
& \qquad\quad
+
\inner{\C:[\boldvarm(\mm_{h}^{i+1}) - \boldvarm(\Pi_h\mm_{h}^{i+1})]}{\boldvar(\uu_h^{i+1}) - \boldvar(\uu_h^{i})} \\
& \quad
\stackrel{\phantom{\eqref{eq:abel}}}{=}
\inner{\C:[\boldvar(\uu_{h}^{i+1}) - \boldvarm(\mm_{h}^{i+1})]}{[\boldvar(\uu_h^{i+1}) -
\boldvarm(\mm_{h}^{i+1})] - [\boldvar(\uu_h^{i})
- \boldvarm(\mm_{h}^i)]}\\
& \qquad\quad
+
\inner{\C:[\boldvar(\uu_{h}^{i+1}) - \boldvarm(\mm_{h}^{i+1})]}{\boldvarm(\mm_h^{i+1}) - \boldvarm(\mm_h^{i})} \\
& \qquad\quad
+
\inner{\C:[\boldvarm(\mm_{h}^{i+1}) - \boldvarm(\Pi_h\mm_{h}^{i+1})]}{\boldvar(\uu_h^{i+1}) - \boldvar(\uu_h^{i})} \\
& \quad
\stackrel{\eqref{eq:abel}}{=}
\frac{1}{2} \norm[\C]{\boldvar(\uu_h^{i+1}) -
\boldvarm(\mm_{h}^{i+1})}^2
-
\frac{1}{2} \norm[\C]{\boldvar(\uu_h^{i})
- \boldvarm(\mm_{h}^i)}^2 \\
& \qquad\quad
+
\frac{1}{2} \norm[\C]{[\boldvar(\uu_h^{i+1}) -
\boldvarm(\mm_{h}^{i+1})] - [\boldvar(\uu_h^{i})
- \boldvarm(\mm_{h}^i)]}^2\\
& \qquad\quad
+
\inner{\C:[\boldvar(\uu_{h}^{i+1}) - \boldvarm(\mm_{h}^{i+1})]}{\boldvarm(\mm_h^{i+1}) - \boldvarm(\mm_h^{i})}\\
& \qquad\quad
+
\inner{\C:[\boldvarm(\mm_{h}^{i+1}) - \boldvarm(\Pi_h\mm_{h}^{i+1})]}{\boldvar(\uu_h^{i+1}) - \boldvar(\uu_h^{i})},
\end{split}
\end{equation*}
respectively.
Altogether, we thus obtain
\begin{equation} \label{eqn:ElasticNextStepEnergyRelation}
\begin{split}
& \E_{\mathrm{el}}[\uu_{h}^{i+1},\mm_{h}^{i+1}]
+ \frac{1}{2} \norm{\dt \uu_{h}^{i+1}}^2
- \E_{\mathrm{el}}[\uu_{h}^i,\mm_{h}^i]
- \frac{1}{2} \norm{\dt \uu_{h}^i}^2 \\
& \quad =
- \frac{1}{2}
\norm{\dt \uu_{h}^{i+1} - \dt\uu_{h}^i}^2    
- \frac{1}{2} \norm[\C]{[\boldvar(\uu_h^{i+1}) -
\boldvarm(\mm_{h}^{i+1})] - [\boldvar(\uu_h^{i})
- \boldvarm(\mm_{h}^i)]}^2 \\
& \qquad
- \inner{\C:[\boldvar(\uu_{h}^{i+1}) - \boldvarm(\mm_{h}^{i+1})]}{\boldvarm(\mm_h^{i+1}) - \boldvarm(\mm_h^{i})}\\
& \qquad
-
\inner{\C:[\boldvarm(\mm_{h}^{i+1}) - \boldvarm(\Pi_h\mm_{h}^{i+1})]}{\boldvar(\uu_h^{i+1}) - \boldvar(\uu_h^{i})}.
\end{split}
\end{equation}
Combining~\eqref{eqn:MagneticNextStepEnergyRelation}
and~\eqref{eqn:ElasticNextStepEnergyRelation} yields
\begin{equation*}
\begin{split}
& \E[\uu_{h}^{i+1},\mm_{h}^{i+1}]
+ \frac{1}{2} \norm{\dt \uu_{h}^{i+1}}^2
- \E[\uu_{h}^i,\mm_{h}^i]
- \frac{1}{2} \norm{\dt \uu_{h}^i}^2 \\
& \quad =
- \alpha k \norm[h]{\vv_{h}^{i}}^2
    - k^2(\theta - 1/2) \norm{\Grad\vv_{h}^{i}}^2
    + k\inner{\hhm[\uu_{h}^{i},\Pi_h\mm_{h}^{i}]}{\vv_{h}^{i}} \\
& \qquad - \frac{1}{2}
\norm{\dt \uu_{h}^{i+1} - \dt\uu_{h}^i}^2    
- \frac{1}{2} \norm[\C]{[\boldvar(\uu_h^{i+1}) -
\boldvarm(\mm_{h}^{i+1})] - [\boldvar(\uu_h^{i})
- \boldvarm(\mm_{h}^i)]}^2 \\
& \qquad
- \inner{\C:[\boldvar(\uu_{h}^{i+1}) - \boldvarm(\mm_{h}^{i+1})]}{\boldvarm(\mm_h^{i+1}) - \boldvarm(\mm_h^{i})} \\
& \qquad
-
\inner{\C:[\boldvar(\mm_{h}^{i+1}) - \boldvarm(\Pi_h\mm_{h}^{i+1})]}{\boldvarm(\uu_h^{i+1}) - \boldvar(\uu_h^{i})}\\
& \quad =
- \alpha k \norm[h]{\vv_{h}^{i}}^2
- D_{h,k}^i
- E_{h,k}^i,
\end{split}
\end{equation*}
where, in the last identity,
we have used the expression of $D_{h,k}^i$ in~\eqref{eq:dissipation}
and
we have defined
\begin{align*}
E_{h,k}^i
&:= \inner{\C:[\boldvar(\uu_{h}^{i+1}) - \boldvarm(\mm_{h}^{i+1})]}{\boldvarm(\mm_h^{i+1}) - \boldvarm(\mm_h^{i})}\\
& \quad - k\inner{\hhm[\uu_{h}^{i},\Pi_h\mm_{h}^{i}]}{\vv_{h}^{i}}\\
& \quad + \inner{\C:[\boldvarm(\mm_{h}^{i+1}) - \boldvarm(\Pi_h\mm_{h}^{i+1})]}{\boldvar(\uu_h^{i+1}) - \boldvar(\uu_h^{i})}
\end{align*}
To conclude the proof of~\eqref{eq:discrete_energy_law},
it remains to show that the latter coincides with~\eqref{eq:error_energy}.
To this end,
using the expression of the elastic field and Lemma~\ref{lem:tensor},
we obtain
\begin{equation*}
\begin{split}
k\inner{\hhm[\uu_{h}^{i},\Pi_h\mm_{h}^{i}]}{\vv_{h}^{i}}
& \,\stackrel{\eqref{eqn:MagnetostrictiveEffectiveField}}{=}
2k\inner{\Z^{\top}\C:[\boldvar(\uu_{h}^{i}) - \boldvarm(\Pi_h\mm_{h}^{i})]\Pi_h \mm_h^i}{\vv_{h}^{i}} \\
& \stackrel{\eqref{eq:tensor_identity}}{=}
2k\inner{\C:[\boldvar(\uu_{h}^{i}) - \boldvarm(\Pi_h\mm_{h}^{i})]}{\Z(\Pi_h \mm_h^i \otimes \vv_{h}^{i})}.
\end{split}
\end{equation*}
Moreover,
from~\eqref{eq:alg_update_m} and the minor symmetry of $\Z$, we get the expansion
\begin{equation} \label{eq:expansion_magnetostrain}
\boldvarm(\mm_{h}^{i+1})
= \boldvarm(\mm_{h}^{i})
+ 2k \, \Z(\mm_{h}^{i} \otimes \vv_{h}^{i})
+ k^2 \boldvarm(\vv_{h}^{i}).
\end{equation}
Altogether, it follows that
\begin{equation*}
\begin{split}
E_{h,k}^i
& =
k^2 \inner{\C:[\boldvar(\uu_{h}^{i+1}) - \boldvarm(\mm_{h}^{i+1})]}{\boldvarm(\vv_{h}^{i})} \\
& \quad
+ 2k \inner{\C:[\boldvar(\uu_{h}^{i+1}) - \boldvarm(\mm_{h}^{i+1})]}{\Z(\mm_{h}^{i} \otimes \vv_{h}^{i})} \\
& \quad
- 2k\inner{\C:[\boldvar(\uu_{h}^{i}) - \boldvarm(\mm_{h}^{i})]}{\Z(\Pi_h \mm_h^i \otimes \vv_{h}^{i})} \\
& \quad
+ \inner{\C:[\boldvarm(\mm_{h}^{i+1}) - \boldvarm(\Pi_h\mm_{h}^{i+1})]}{\boldvar(\uu_h^{i+1}) - \boldvar(\uu_h^{i})}\\
& =
k^2 \inner{\C:[\boldvar(\uu_{h}^{i+1}) - \boldvarm(\mm_{h}^{i+1})]}{\boldvarm(\vv_{h}^{i})} \\
& \quad
+ 2k \inner{\C:\{[\boldvar(\uu_{h}^{i+1}) - \boldvarm(\mm_{h}^{i+1})] - [\boldvar(\uu_{h}^{i}) - \boldvarm(\mm_{h}^{i})] \}}{\Z(\mm_{h}^{i} \otimes \vv_{h}^{i})} \\
& \quad
+ 2k\inner{\C:[\boldvar(\uu_{h}^{i}) - \boldvarm(\mm_{h}^{i})]}{\Z[(\mm_h^i - \Pi_h \mm_h^i) \otimes \vv_{h}^{i}]}\\
& \quad
+ \inner{\C:[\boldvarm(\mm_{h}^{i+1}) - \boldvarm(\Pi_h\mm_{h}^{i+1})]}{\boldvar(\uu_h^{i+1}) - \boldvar(\uu_h^{i})}.
\end{split}
\end{equation*}
This shows~\eqref{eq:error_energy} and concludes the proof.
\end{proof}

\subsection{Stability} \label{sec:proofs_stability}

We now prove Proposition~\ref{prop:stability} showing
unconditional stability of Algorithm~\ref{algorithm}
and an estimate of the violation
of the unit length constraint.
For the sake of clarity, we split the proof into several lemmas.

An immediate consequence of the projection-free update~\eqref{alg1:magnetisation_update}
is the following $L^2$-bound for the approximate magnetisations.

\begin{lemma}\label{lem:Magnetisation_Estimate}
For every integer $1 \le j \le N$,
it holds that
\begin{align}\label{alg2:m_bound_estimate}
    \norm{\mm_{h}^j}^2
    &\le
    C_1 \left(
    1 + k^2 \sum_{i=0}^{j-1}\norm{\vv_{h}^{i}}^2
    \right),\\
\label{eq:auxiliary_constraint}
\big\lVert \II_h\big[\abs{\mm_h^j}^2\big]-1 \big\rVert_{L^1(\Omega)}
&\le C_2 k^2 \sum_{i=0}^{j-1}\norm{\vv_{h}^{i}}^2.
\end{align}
where $C_1,C_2>0$ are constants depending 
the shape-regularity parameter of $\T_h$
($C_1$ depends also on $\abs{\Omega}$).
\end{lemma}

\begin{proof}
We follow~\cite{bartels2016}.
Starting from~\eqref{eq:alg_update_m}
and noting that $\vv_{h}^{i}\in\KK_{h}[\mm_{h}^{i}]$,
we have for each $z\in\NN_{h}$ that for every $0 \le i \le j-1$
\begin{equation*}
\abs{\mm_{h}^{i+1}(z)}^2
= \abs{\mm_{h}^{i}(z)}^2 + k^2\abs{\vv_{h}^{i}(z)}^2.
\end{equation*}
Inductively, starting with $|\mm_{h}^{0}(z)| = 1$, we deduce that
\begin{equation*}
    |\mm_{h}^{j}(z)|^2= 1+ k^2\sum_{i=0}^{j-1}|\vv_{h}^{i}(z)|^2.
\end{equation*}
Then, noting that $\norm{1}=\abs{\Omega}^{1/2}$ and using~\eqref{eq:Lp_equivalence} yields~\eqref{alg2:m_bound_estimate}
(for a suitable constant $C_1>0$ we do not explicitly compute).
The same argument shows~\eqref{eq:auxiliary_constraint}.
%
%
%
%
%
%
\end{proof}

We also have the following estimate of all quantities involving the magnetisation.

\begin{lemma}\label{lem:StabilityOfMagnetisation}
For every integer $1 \le j \le N$, it holds that
\begin{multline} \label{eqn:MagneticStability}
\norm{\Grad\mm_{h}^{j}}^2
+ k \sum_{i=0}^{j-1} \norm{\vv_{h}^{i}}^2
+ \left(\theta - \frac{1}{2}\right) k^2\sum_{i=0}^{j-1} \norm{\Grad\vv_{h}^{i}}^2 \\
\leq C_3 \left[ \norm{\Grad\mm_{h}^0}^2
+ k \sum_{i=0}^{j-1}\left(
1 + \norm{\boldvar(\uu_{h}^{i})}^2 \right)
\right],
\end{multline}
where $C_3>0$ depends only on $\alpha$, $\abs{\Omega}$, $\norm[\LL^{\infty}(\Omega)]{\Z}$, and $\norm[\LL^{\infty}(\Omega)]{\C}$.
\end{lemma}

\begin{proof}
Let $1 \le j \le N$ be an integer.
Starting from~\eqref{eqn:MagneticNextStepEnergyRelation}
(cf.\ the proof of Proposition~\ref{prop:energy}),
we sum up from $0$ to $j-1$ to obtain
\begin{multline*}
\frac{1}{2}\norm{\Grad\mm_{h}^{j}}^2
+ \alpha k\sum_{i=0}^{j-1} \norm[h]{\vv_{h}^{i}}^2
+ \left(\theta - \frac{1}{2}\right) k^2\sum_{i=0}^{j-1} \norm{\Grad\vv_{h}^{i}}^2 \\
= \frac{1}{2}\norm{\Grad\mm_{h}^0}^2
+ k\sum_{i=0}^{j-1}
\inner{\hhm[\uu_{h}^{i},\Pi_h\mm_{h}^{i}]}{\vv_{h}^{i}}.
\end{multline*}
Using Lemma~\ref{lem:boundedness_elastic_field},
we can estimate the term involving the elastic field
for some $\nu >0$ by
\begin{equation*}
\begin{split}
\abs{\inner{\hhm[\uu_{h}^{i},\Pi_h\mm_{h}^{i}]}{\vv_{h}^{i}}}
& \stackrel{\phantom{\eqref{eq:boundedness_elastic_field}}}{\leq}
\frac{1}{4\nu}\norm{\hhm[\uu_{h}^{i},\Pi_h\mm_{h}^{i}]}^2
+ \nu\norm{\vv_{h}^{i}}^2\\
& \stackrel{\eqref{eq:boundedness_elastic_field}}{\leq}
\frac{2}{\nu}
\norm[\LL^{\infty}(\Omega)]{\Z}^2
\norm[\LL^{\infty}(\Omega)]{\C}^2
\left(\norm{\boldvar(\uu_{h}^{i})}^2 + \norm[\LL^{\infty}(\Omega)]{\Z}^2 |\Omega| \right)
+ \nu\norm{\vv_{h}^{i}}^2.
\end{split}
\end{equation*}
Then we get
\begin{multline*}
\frac{1}{2}\norm{\Grad\mm_{h}^{j}}^2
+ \alpha k\sum_{i=0}^{j-1} \norm[h]{\vv_{h}^{i}}^2
+ \left(\theta - \frac{1}{2}\right) k^2\sum_{i=0}^{j-1} \norm{\Grad\vv_{h}^{i}}^2
\le \frac{1}{2}\norm{\Grad\mm_{h}^0}^2 \\
+ \frac{2k}{\nu}\sum_{i=0}^{j-1}
\norm[\LL^{\infty}(\Omega)]{\Z}^2
\norm[\LL^{\infty}(\Omega)]{\C}^2
\left(\norm{\boldvar(\uu_{h}^{i})}^2
+ \norm[\LL^{\infty}(\Omega)]{\Z}^2 |\Omega| \right)
+ \nu k\sum_{i=0}^{j-1} \norm{\vv_{h}^{i}}^2.
\end{multline*}
Using~\eqref{eq:h-scalar-product_equivalence}
and choosing $\nu = \alpha/2$ yields~\eqref{eqn:MagneticStability}
(for a suitable constant $C_3>0$ which we do not compute explicitly).
\end{proof}

In the following lemma, we show that the magnetostrain is Lipschitz continuous
with respect to the magnetisation
(the use of the nodal projection is exploited here).

\begin{lemma} \label{lem:lipschitz_magnetostrain}
For all $\mm_{h,1} , \mm_{h,2} \in \S^1(\T_h)^3$
satisfying $\abs{\mm_{h,\ell}(z)} \ge 1$ for all $\ell=1,2$ and $z \in \NN_h$,
it holds that
\begin{equation} \label{eq:lipschitz_magnetostrain}
\norm{\boldvarm(\Pi_h\mm_{h,1}) - \boldvarm(\Pi_h\mm_{h,2})}
\le C_{\mathrm{m}}
\norm{\mm_{h,1} - \mm_{h,2}},
\end{equation}
where $C_{\mathrm{m}}>0$ depends only on $\norm[\LL^{\infty}(\Omega)]{\Z}$
and the shape-regularity parameter of $\T_h$.
\end{lemma}

\begin{proof}
Straightforward calculations exploiting the boundedness guaranteed by the nodal projection,
i.e.\ $\norm[\LL^\infty(\Omega)]{\Pi_h\mm_{h,1}} = \norm[\LL^\infty(\Omega)]{\Pi_h\mm_{h,2}} = 1$,
show that
\begin{equation*}
\norm{\boldvarm(\Pi_h\mm_{h,1}) - \boldvarm(\Pi_h\mm_{h,2})}
\lesssim
\norm{\Pi_h\mm_{h,1} - \Pi_h\mm_{h,2}},
\end{equation*}
where the hidden constant depends on $\norm[\LL^{\infty}(\Omega)]{\Z}$.
From the norm equivalence in~\cite[Lemma 3.4]{bartels2015})
and the fact that the projection onto the sphere is non-expanding
(i.e.\ Lipschitz continuous with constant 1),
it follows that
\begin{equation*}
\norm{\Pi_h\mm_{h,1} - \Pi_h\mm_{h,2}}
\lesssim
\norm{\mm_{h,1} - \mm_{h,2}},
\end{equation*}
where the hidden constant depends on the shape-regularity of the mesh.
Combining the above two estimates yields the desired result,
where $C_{\mathrm{m}}>0$ is the product of the two constants hidden above.
\end{proof}

\begin{lemma} \label{lem:StrainBoundedness}
For every integer $1 \le j \le N$, the following estimate holds
\begin{multline} \label{eq:StrainBoundedness}
\norm{\dt \uu_{h}^j}^2
+ \norm{\boldvar(\uu_{h}^{j})}^2
+ \sum_{i=0}^{j-1}\norm{\dt \uu_{h}^{i+1} - \dt\uu_{h}^i}^2
+ \sum_{i=0}^{j-1}\norm{\boldvar(\uu_{h}^{i+1})-\boldvar(\uu_{h}^{i})}^2 \\
\le
C_4 \left[ 1
+ \norm{\dot\uu_{h}^0}^2
+ \norm{\boldvar(\uu_{h}^{0})}^2
+ \norm{\Grad\mm_{h}^0}^2
+ k \sum_{i=0}^{j-1}\left(
1 + \norm{\boldvar(\uu_{h}^{i})}^2 \right)
\right],
\end{multline}
where $C_4>0$ depends only on the shape-regularity parameter of $\T_h$
and the problem data $\alpha$, $\Omega$, $\C$, $\Z$, $\ff$ and $\gg$.
\end{lemma}

\begin{proof}
Let $1 \le j \le N$ be an integer.
Starting from~\eqref{eqn:AfterAbelFirstSummation}
(cf.\ the proof of Proposition~\ref{prop:energy}),
summing up from $0$ to $j-1$, we have
\begin{multline*}
\frac{1}{2}\norm{\dt \uu_{h}^{j}}^2
-\frac{1}{2}\norm{\dt \uu_{h}^0}^2
+\frac{1}{2}\sum_{i=0}^{j-1}\norm{\dt \uu_{h}^{i+1}
- \dt\uu_{h}^i}^2
+ \sum_{i=0}^{j-1} \inner{\C:\boldvar(\uu_{h}^{i+1})}{\boldvar(\uu_h^{i+1}) - \boldvar(\uu_h^{i})} \\
= \sum_{i=0}^{j-1} \inner{\C:\boldvarm(\Pi_h\mm_{h}^{i+1})}{\boldvar(\uu_h^{i+1}) - \boldvar(\uu_h^{i})} \\
+\inner{\ff}{\uu_h^{j}}
-\inner{\ff}{\uu_h^{0}}
+ \inner[\Gamma_{N}]{\gg}{\uu_h^{j}}
-\inner[\Gamma_{N}]{\gg}{\uu_h^{0}}.
\end{multline*}
%
%
%
Applying Lemma~\ref{lem:abel} to the last term on the left-hand side and rearranging we have
\begin{multline*}
\frac{1}{2}\norm{\dt \uu_{h}^j}^2
+ \frac{1}{2}\norm[\C]{\boldvar(\uu_{h}^{j})}^2
+ \frac{1}{2}\sum_{i=0}^{j-1}\norm{\dt \uu_{h}^{i+1} - \dt\uu_{h}^i}^2
+ \frac{1}{2} \sum_{i=0}^{j-1}\norm[\C]{\boldvar(\uu_{h}^{i+1})-\boldvar(\uu_{h}^{i})}^2 \\
= \frac{1}{2}\norm{\dt \uu_{h}^0}^2
+ \frac{1}{2}\norm[\C]{\boldvar(\uu_{h}^{0})}^2
+ \sum_{i=0}^{j-1}\inner{\C:\boldvarm(\Pi_h\mm_{h}^{i+1})}{\boldvar(\uu_h^{i+1}) - \boldvar(\uu_h^{i})}\\
+ \inner{\ff}{\uu_h^{j}}
- \inner{\ff}{\uu_h^{0}}
+ \inner[\Gamma_{N}]{\gg}{\uu_h^{j}}
- \inner[\Gamma_{N}]{\gg}{\uu_h^{0}}.
\end{multline*}
The term involving the magnetostrain can be estimated as
\begin{equation*}
\begin{split}
&
\sum_{i=0}^{j-1}\inner{\C:\boldvarm(\Pi_h\mm_{h}^{i+1})}{\boldvar(\uu_h^{i+1}) - \boldvar(\uu_h^{i})} \\
& \ =
\inner{ \boldvarm(\Pi_{h}\mm_{h}^{j})}{\C:\boldvar(\uu_{h}^{j})}
- \inner{\boldvarm(\Pi_{h}\mm_{h}^{1})}{\C:\boldvar(\uu_{h}^{0})}
- k \sum_{i=1}^{j-1} \inner{\dt\boldvarm(\Pi_{h}\mm_{h}^{i+1})}{\C:\boldvar(\uu_{h}^{i})} \\
& \
\leq
\norm{ \boldvarm(\Pi_{h}\mm_{h}^{j})}^2
+ \frac{1}{4}\norm[\C]{\boldvar(\uu_{h}^{j})}^2
+ \frac{1}{2}\norm[\C]{\boldvarm(\Pi_{h}\mm_{h}^{1})}^2
+ \frac{1}{2}\norm[\C]{\boldvar(\uu_{h}^{0})}^2 \\
& \quad
+ \frac{k}{2}\sum_{i=1}^{j-1}\norm{\dt\boldvarm(\Pi_{h}\mm_{h}^{i+1})}^2
+ \frac{k}{2}\sum_{i=1}^{j-1}\norm[\C]{\boldvar(\uu_{h}^{i})}^2 \\
& \ \leq
\frac{3}{2} \norm[\LL^\infty(\Omega)]{\Z}^2 \abs{\Omega}
+ \frac{1}{4}\norm[\C]{\boldvar(\uu_{h}^{j})}^2
+ \frac{1}{2}\norm[\C]{\boldvar(\uu_{h}^{0})}^2 \\
& \quad
+ \frac{k}{2}\sum_{i=1}^{j-1}\norm{\dt\boldvarm(\Pi_{h}\mm_{h}^{i+1})}^2
+ \frac{k}{2}\sum_{i=1}^{j-1}\norm[\C]{\boldvar(\uu_{h}^{i})}^2.
\end{split}
\end{equation*}
%
%
%
Using Lemma~\ref{lem:lipschitz_magnetostrain}, we get
\begin{equation*}
\begin{split}
  \norm{\dt\boldvarm(\Pi_{h}\mm_{h}^{i+1})}
    &\stackrel{\phantom{\eqref{eq:lipschitz_magnetostrain}}}{=} \frac{1}{k}\norm{\boldvarm(\Pi_{h}\mm_{h}^{i+1})-\boldvarm(\Pi_{h}\mm_{h}^{i})} \\
    &\stackrel{\eqref{eq:lipschitz_magnetostrain}}{\le} \frac{1}{k}C_{\mathrm{m}} \norm{\mm_{h}^{i+1}-\mm_{h}^{i}}
    = C_{\mathrm{m}} \norm{\vv_{h}^{i}}.
\end{split}
\end{equation*}
It follows that
\begin{equation*}
\frac{k}{2}\sum_{i=1}^{j-1}\norm{\dt\boldvarm(\Pi_{h}\mm_{h}^{i+1})}^2
\leq \frac{C_{\mathrm{m}}^2 \, k}{2}\sum_{i=1}^{j-1}\norm{\vv_{h}^{i}}^2.
\end{equation*}
Moreover, for every $\delta >0$ we have that
\begin{align*}
\begin{split}
\abs{\inner{\ff}{\uu_h^j}
+ \inner[\Gamma_{N}]{\gg}{\uu_h^j}}
& \leq C_{\mathrm{KPC}} (\norm{\ff} + \norm[\Gamma_{N}]{\gg}) \norm[\C]{\boldvar(\uu_h^j)} \\
&\leq \frac{C_{\mathrm{KPC}}^2}{4\delta}(\norm{\ff} + \norm[\Gamma_{N}]{\gg})^2+ \delta\norm[\C]{\boldvar(\uu_h^j)}^2,
\end{split}
\end{align*}
where $C_{\mathrm{KPC}}>0$ is a constant depending only on $\abs{\Omega}$ and $\C$
(a combination of the continuity constant of the trace operator $\HH^1(\Omega) \to \LL^2(\Gamma_N)$,
the constants appearing in Poincar\'e's and Korn's inequalities, and the equivalence constant in the
norm equivalence $\norm{\cdot} \simeq \norm[\C]{\cdot}$).
Overall, choosing $\delta = 1/8$ and recalling that $\dt \uu_{h}^0 = \dot\uu_{h}^0$, we obtain
\begin{equation*}
\begin{split}
& \frac{1}{2}\norm{\dt \uu_{h}^j}^2
+ \frac{1}{8}\norm[\C]{\boldvar(\uu_{h}^{j})}^2
+ \frac{1}{2}\sum_{i=0}^{j-1}\norm{\dt \uu_{h}^{i+1} - \dt\uu_{h}^i}^2
+ \frac{1}{2} \sum_{i=0}^{j-1}\norm[\C]{\boldvar(\uu_{h}^{i+1})-\boldvar(\uu_{h}^{i})}^2 \\
& \quad\le
\frac{3}{2} \norm[\LL^\infty(\Omega)]{\Z}^2 \abs{\Omega}
+ 2 C_{\mathrm{KPC}}^2(\norm{\ff} + \norm[\Gamma_{N}]{\gg})^2
+ \frac{1}{2}\norm{\dot\uu_{h}^0}^2 \\
& \qquad + \left[1 + C_{\mathrm{KPC}} \left(\norm{\ff} + \norm[\Gamma_{N}]{\gg}\right)\right]
\norm[\C]{\boldvar(\uu_{h}^{0})}^2
+ \frac{k}{2}\sum_{i=1}^{j-1}\norm[\C]{\boldvar(\uu_{h}^{i})}^2
+ \frac{C_{\mathrm{m}}^2 \, k}{2}\sum_{i=1}^{j-1}\norm{\vv_{h}^{i}}^2.
\end{split}
\end{equation*}
Applying Lemma~\ref{lem:StabilityOfMagnetisation} to estimate the last term on the right-hand side,
we obtain~\eqref{eq:StrainBoundedness}
(for a suitable constant $C_4>0$ which we do not compute explicitly).
\end{proof}

We are now in a position to prove Proposition~\ref{prop:stability}.

\begin{proof}[Proof of Proposition~\ref{prop:stability}]
We apply Lemmas~\ref{lem:Magnetisation_Estimate}, \ref{lem:StabilityOfMagnetisation} and~\ref{lem:StrainBoundedness}.
Combining~\eqref{alg2:m_bound_estimate}, \eqref{eqn:MagneticStability} and~\eqref{eq:StrainBoundedness}, we obtain
\begin{equation*} 
\begin{split}
& \norm{\dt \uu_{h}^j}^2
+ \norm{\boldvar(\uu_{h}^{j})}^2
+ \sum_{i=0}^{j-1}\norm{\dt \uu_{h}^{i+1} - \dt\uu_{h}^i}^2
+ \sum_{i=0}^{j-1}\norm{\boldvar(\uu_{h}^{i+1})-\boldvar(\uu_{h}^{i})}^2 \\
& \qquad
+ \norm[\HH^1(\Omega)]{\mm_{h}^{j}}^2
+ (1-C_1 k) k \sum_{i=0}^{j-1} \norm{\vv_{h}^{i}}^2
+ \left(\theta - \frac{1}{2}\right) k^2\sum_{i=0}^{j-1} \norm{\Grad\vv_{h}^{i}}^2 \\
& \qquad\qquad
\leq
C_1 + C_4
+ (C_3 + C_4) \norm{\Grad\mm_{h}^0}^2
+ C_4 \norm{\dot\uu_{h}^0}^2
+ C_4 \norm{\boldvar(\uu_{h}^{0})}^2 \\
& \qquad\qquad\qquad
+ (C_3 + C_4)k \sum_{i=0}^{j-1}\left(
1 + \norm{\boldvar(\uu_{h}^{i})}^2 \right) \\
& \qquad\qquad
\leq
C_1 + C_4
+ (C_3 + C_4) \norm{\Grad\mm_{h}^0}^2
+ C_4 \norm{\dot\uu_{h}^0}^2
+ C_4 \norm{\boldvar(\uu_{h}^{0})}^2 \\
& \qquad\qquad\qquad
+ (C_3 + C_4)T
+ (C_3 + C_4)k \sum_{i=0}^{j-1}\norm{\boldvar(\uu_{h}^{i})}^2,
\end{split}
\end{equation*}
where in the last estimate we have used that $kj \le T$.
If the time-step size $k$ is sufficiently small, the coefficients in front of all terms on the left-hand side are strictly positive.
Given the boundedness of the approximate initial data guaranteed by assumption~\eqref{eq:convergence_initial_data}, 
the desired stability estimate \eqref{eq:stability} then follows from the
discrete Gr\"onwall lemma; see e.g.\ \cite[Lemma~10.5]{thomee2006}.
Finally,
\eqref{eq:general_constraint}
follows from~\eqref{eq:stability}
and~\eqref{eq:auxiliary_constraint}.
This concludes the proof.
\end{proof}

\subsection{Convergence} \label{sec:proofs_convergence}

The proof of convergence of Algorithm~\ref{algorithm}
(Theorem~\ref{thm:convergence}(i)) follows the standard argument
to prove existence of weak solutions for parabolic equations
(uniform boundedness of Galerkin approximations,
extraction of subsequences with suitable convergence properties,
identification of the limit with a weak solution of the problem;
see, e.g., \cite[Section~7.1]{evans2010})
and thus has the same structure as the one which proves
the convergence of~\cite[Algorithm~4.1]{bppr2013}.
Therefore, in the upcoming analysis, we will provide only a sketch
of the steps of the proof that can be found in~\cite{bppr2013}.
However, we will present in detail the (non-obvious) steps
that we have to perform to cope with the partial omission of the nodal projection
(for which we borrow ideas from~\cite{abert_spin-polarized_2014,bartels2016})
and to prove our novel energy estimate.

We start the proof with showing the following lemma,
which provides an estimate of the $L^p$-norm ($p \ge 1$)
of the difference between the approximate magnetisations generated by Algorithm~\ref{algorithm}
and their nodal projections.

\begin{lemma} \label{lem:projection_error_Lp}
Let $p \in [1,\infty)$.
For all integers $1 \le j \le N$, it holds that
\begin{equation} \label{eq:projection_error_Lp}
\norm[\LL^p(\Omega)]{\mm_h^j - \Pi_h \mm_h^j}
\le
C \frac{T^{1-1/p}}{2} k^{1+1/p} \sum_{i=0}^{j-1} \norm[\LL^{2p}(\Omega)]{\vv_h^i}^{2},
\end{equation}
where $C>0$ depends only on the shape-regularity of $\T_h$.
\end{lemma}

\begin{proof}
Let $1 \le j \le N$ be an integer.
For all $z \in \NN_h$, we have that
\begin{equation*}
\begin{split}
\abs{\mm_h^j(z) - \Pi_h \mm_h^j(z)}
& = \abs{\mm_h^j(z) - \frac{\mm_h^j(z)}{\abs{\mm_h^j(z)}}}
= \abs{\mm_h^j(z)} - 1 \\
& = \frac{\abs{\mm_h^j(z)}^2 - 1}{\abs{\mm_h^j(z)} + 1}
\le \frac{1}{2} \left( \abs{\mm_h^j(z)}^2 - 1 \right)
= \frac{k^2}{2} \sum_{i=0}^{j-1} \abs{\vv_h^i(z)}^2.
\end{split}
\end{equation*}
If $p=1$, the norm equivalence~\eqref{eq:Lp_equivalence} immediately yields~\eqref{eq:projection_error_Lp}.
If $p>1$, applying~\eqref{eq:Lp_equivalence} twice and using the convexity of $x^p$ for $x>0$
as well as $jk \le T$, we obtain
\begin{equation*}
\begin{split}
\norm[\LL^p(\Omega)]{\mm_h^j - \Pi_h \mm_h^j}^p
&\lesssim
\sum_{z \in \NN_h} h_z^3 \abs{\mm_h^j(z) - \Pi_h \mm_h^j(z)}^p
\le
\sum_{z \in \NN_h} h_z^3 \left(\frac{k^2}{2} \sum_{i=0}^{j-1} \abs{\vv_h^i(z)}^2 \right)^p \\
& \le
\sum_{z \in \NN_h} h_z^3 \frac{k^{2p}}{2^p} j^{p-1} \sum_{i=0}^{j-1} \abs{\vv_h^i(z)}^{2p}
\le
\sum_{z \in \NN_h} h_z^3 \frac{k^{p+1}}{2^p} T^{p-1} \sum_{i=0}^{j-1} \abs{\vv_h^i(z)}^{2p} \\
& \lesssim
\frac{k^{p+1}}{2^p} T^{p-1} \sum_{i=0}^{j-1} \norm[\LL^{2p}(\Omega)]{\vv_h^i}^{2p},
\end{split}
\end{equation*}
where the hidden constants depend only on the shape-regularity of $\T_h$.
Then, \eqref{eq:projection_error_Lp} for $p>1$ follows from the inequality
$\norm[\ell^p]{\cdot} \le \norm[\ell^1]{\cdot}$ satisfied by the $p$-norms in finite dimensions.
This concludes the proof.
\end{proof}

Now,
let $\{ \mm_{hk}\}$, $\{ \mm_{hk}^{\pm} \}$, $\{ \vv_{hk}^- \}$,
$\{ \uu_{hk}\}$, $\{ \uu_{hk}^\pm \}$,
$\{ \dot\uu_{hk}\}$, $\{ \dot\uu_{hk}^\pm \}$
be the time reconstructions
defined according to~\eqref{eq:reconstructions}
using the approximations $\{ (\uu_h^i,\mm_h^i) \}_{0 \le i \le N}$
generated by Algorithm~\ref{algorithm}.
In the following lemma,
we show that the uniform stability established in Proposition~\ref{prop:stability}
allows us to extract convergent subsequences from the sequences of time reconstructions.

\begin{lemma} \label{lem:convergent_subsequences}
Under the assumptions of Theorem~\ref{thm:convergence}(i),
there exist
$\uu \in L^\infty(0,T;\HH^1_D(\Omega))$ with $\uut \in L^\infty(0,T;\LL^2(\Omega))$
and
$\mm \in L^\infty(0,T;\HH^1(\Omega;\sphere))$ with $\mmt \in L^2(0,T;\LL^2(\Omega))$
such that, upon extraction of (non-relabeled) subsequences,
we have the following convergence results:
\begin{subequations} \label{eq:convergences}
\begin{align}
\label{eq:conv_weakH1_u}
\uu_{hk} \weakto \uu \quad &\text{in } \HH^1(\Omega_T), \\
\label{eq:conv_weakstarLinftyH1_u}
\uu_{hk}, \uu_{hk}^{\pm} \weakstarto \uu \quad &\text{in } L^{\infty}(0,T;\HH^1(\Omega)), \\
\label{eq:conv_weakL2H1_u}
\uu_{hk}, \uu_{hk}^{\pm} \weakto \uu \quad &\text{in } L^2(0,T;\HH^1(\Omega)), \\
\label{eq:conv_strongL2_u}
\uu_{hk}, \uu_{hk}^{\pm} \to \uu \quad &\text{in } \LL^2(\Omega_T), \\
\label{eq:conv_weakLinftyL2_ut}
\dot\uu_{hk}, \dot\uu_{hk}^{\pm} \weakstarto \uut \quad &\text{in } L^\infty(0,T;\LL^2(\Omega)), \\
\label{eq:conv_weakL2_ut}
\dot\uu_{hk}, \dot\uu_{hk}^{\pm} \weakto \uut \quad &\text{in } \LL^2(\Omega_T), \\
\label{eq:conv_weakH1}
\mm_{hk} \weakto \mm \quad &\text{in } \HH^1(\Omega_T), \\
\label{eq:conv_strongHs}
\mm_{hk} \to \mm \quad &\text{in } \HH^s(\Omega_T) \text{ for all } s \in (0,1), \\
\label{eq:conv_weakstarLinftyH1}
\mm_{hk}, \mm_{hk}^{\pm} \weakstarto \mm \quad &\text{in } L^{\infty}(0,T;\HH^1(\Omega)), \\
\label{eq:conv_weakL2H1}
\mm_{hk}, \mm_{hk}^{\pm} \weakto \mm \quad &\text{in } L^2(0,T;\HH^1(\Omega)), \\
\label{eq:conv_strongL2Hs}
\mm_{hk}, \mm_{hk}^{\pm} \to \mm \quad &\text{in } L^2(0,T;\HH^s(\Omega)) \text{ for all } s \in (0,1),\\
\label{eq:conv_strongL2}
\mm_{hk}, \mm_{hk}^{\pm} \to \mm \quad &\text{in } \LL^2(\Omega_T), \\
\label{eq:conv_pointwise}
\mm_{hk}, \mm_{hk}^{\pm} \to \mm \quad &\text{pointwise a.e.\ in } \Omega_T, \\
\label{eq:conv_weakL2_v}
\vv_{hk}^- \weakto \mmt \quad &\text{in } \LL^2(\Omega_T),
\end{align}
\end{subequations}
as $h,k \to 0$.
\end{lemma}

\begin{proof}
Using the boundedness expressed in Proposition~\ref{prop:stability},
we can successively extract weakly(-star) convergent subsequences
(non-relabeled, with possibly different limits)
from $\{\uu_{hk}\}$ and $\{\uu_{hk}^{\pm}\}$,
from $\{\dot{\uu}_{hk}\}$ and $\{\dot{\uu}_{hk}^{\pm}\}$,
from $\{\mm_{hk}\}$ and $\{\mm_{hk}^{\pm}\}$,
and
from $\{\vv_{hk}^{-}\}$.

Let $\uu\in\HH^1(\Omega_T)$ satisfy the weak convergence~\eqref{eq:conv_weakH1_u}.
Owing to the continuous inclusions
$\HH^1(\Omega_T) \subset L^2(0,T;\HH^1(\Omega)) \subset \LL^2(\Omega_T)$
and
the compact inclusion $\HH^1(\Omega_T) \Subset \LL^2(\Omega_T)$,
we obtain convergences~\eqref{eq:conv_weakL2H1_u} and~\eqref{eq:conv_strongL2_u}.
Moreover, from the continuous inclusion
$L^{\infty}(0,T;\HH^1(\Omega)) \subset L^2(0,T;\HH^1(\Omega))$,
we can identify the weak-star limit of $\{\uu_{hk}\}$ in $L^\infty(0,T;\HH^1(\Omega))$
with the weak limit in $L^2(0,T;\HH^1(\Omega))$,
which shows~\eqref{eq:conv_weakstarLinftyH1_u} for $\{\uu_{hk}\}$.

Let $\mm\in\HH^1(\Omega_T)$ satisfy the weak convergence~\eqref{eq:conv_weakH1}.
Arguing as above and using a well-known result for convergence in $L^p$-spaces,
we obtain convergences~\eqref{eq:conv_weakL2H1},
\eqref{eq:conv_strongL2} and (upon extraction of a further subsequence)
\eqref{eq:conv_pointwise} for $\{\mm_{hk}\}$.
The continuous inclusion
$L^{\infty}(0,T;\HH^1(\Omega)) \subset L^2(0,T;\HH^1(\Omega))$,
shows~\eqref{eq:conv_weakstarLinftyH1} for $\{\mm_{hk}\}$.

Let $0 < s < 1$ be arbitrary.
Since
$\HH^s(\Omega_T) = [\LL^2(\Omega_T),\HH^1(\Omega_T)]_s$
and $L^2(0,T;\HH^s(\Omega)) = [\LL^2(\Omega_T),L^2(0,T;\HH^1(\Omega))]_s$,
well-known results from interpolation theory
(see, e.g., \cite[Theorem~6.4.5 and Theorem~3.8.1]{bl1976} and \cite[Theorem~5.1.2]{bl1976})
yield the compact embedding $\HH^1(\Omega_T) \Subset \HH^s(\Omega_T)$
and the continuous inclusion $\HH^s(\Omega_T) \subset L^2(0,T;\HH^s(\Omega))$.
These in turn show
convergences~\eqref{eq:conv_strongHs} and~\eqref{eq:conv_strongL2Hs}
for $\{\mm_{hk}\}$.
Furthermore,
\eqref{eq:conv_weakL2_v} follows directly from $\mmt_{hk}=\vv_{hk}^{-}$.

Overall, this shows the convergence results~\eqref{eq:conv_weakH1_u}--\eqref{eq:conv_strongL2_u} and \eqref{eq:conv_weakH1}--\eqref{eq:conv_weakL2_v}
for the sequences $\{\uu_{hk}\}$, $\{\mm_{hk}\}$ and $\{\vv_{hk}^-\}$.
Using the same argument, one can obtain the same results
for $\{ \uu_{hk}^{\pm} \}$
and $\{ \mm_{hk}^{\pm} \}$.
Since the quantity
\begin{equation} \label{eq:quantity_to_identify_limits}
    \sum_{i=0}^{j-1} \norm{\mm_{h}^{i+1}- \mm_{h}^{i}}^2 + \sum_{i=0}^{j-1}\norm{\uu_{h}^{i+1} - \uu_{h}^{i}}^2 + \sum_{i=0}^{j-1} \norm{\dt \uu_{h}^{i+1} - \dt \uu_{h}^{i}}
\end{equation}
is uniformly bounded,
arguing as in \cite[Lemma~5.7]{bppr2013}
we can show that the limits of $\{ \uu_{hk} \}$ and $\{ \uu_{hk}^{\pm} \}$
(resp.\ $\{ \mm_{hk} \}$ and $\{ \mm_{hk}^{\pm} \}$)
coincide.
The continuous inclusion 
$L^{\infty}(0,T;\LL^2(\Omega)) \subset L^2(0,T;\LL^2(\Omega)) = \LL^2(\Omega_T)$,
the boundedness of the third term in~\eqref{eq:quantity_to_identify_limits}
and the identity $\uut_{hk} = \dot\uu_{hk}^+$
imply \eqref{eq:conv_weakLinftyL2_ut}--\eqref{eq:conv_weakL2_ut}.
Finally, the fact that $\mm$ satisfies $\abs{\mm}=1$ a.e.\ in $\Omega$
follows from the available convergence results and~\eqref{eq:general_constraint}.
For the details, we refer to Step~3 of the proof of~\cite[Proposition~6]{hpprss2019}.
This concludes the proof.
\end{proof}

Let $\{ \hat\mm_{hk}^\pm\}$ be the piecewise constant time reconstructions
defined using the projection of the approximate magnetisations, i.e,
$\hat\mm_{hk}^-(t) := \Pi_h \mm_h^i$ and $\hat\mm_{hk}^+(t) := \Pi_h \mm_h^{i+1}$
for all $i=0,\dots,N_1$ and $t \in [t_i,t_{i+1})$
(cf.\ \eqref{eq:reconstructions}).

In the following lemma, we establish further convergence results that
will be needed to identify the limit $(\uu,\mm)$ constructed in Lemma~\ref{lem:convergent_subsequences} 
with a weak solution of~\eqref{eq:newton}--\eqref{eq:ibc}.

\begin{lemma}[auxiliary convergences] \label{lem:convergent_subsequences_aux}
Under the assumptions of Theorem~\ref{thm:convergence}(i),
upon extraction of a further (non-relabeled) subsequence,
we have the following convergence results:
\begin{subequations}
\begin{align}
\label{eq:conv_mhat_LinftyH1}
\hat\mm_{hk}^\pm \weakstarto \mm \quad &\text{in } L^\infty(0,T;\HH^1(\Omega)), \\
\hat\mm_{hk}^\pm \weakto \mm \quad &\text{in } L^2(0,T;\HH^1(\Omega)), \\
\label{eq:conv_mhat_L2}
\hat\mm_{hk}^\pm \to \mm \quad &\text{in } \LL^2(\Omega_T), \\
\label{eq:conv_m2_L2}
\mm_{hk}^\pm\otimes\mm_{hk}^\pm \to \mm\otimes\mm \quad &\text{in } \LL^2(\Omega_T), \\
\label{eq:conv_mhat2_L2}
\hat\mm_{hk}^\pm\otimes\hat\mm_{hk}^\pm \to \mm\otimes\mm \quad &\text{in } \LL^2(\Omega_T),
\end{align}
\end{subequations}
as $h,k \to 0$.
\end{lemma}

\begin{proof}
Firstly, we note that $\norm[\LL^\infty(\Omega)]{\Pi_h \mm_h^i} = 1$
and $\norm{\Grad\Pi_h \mm_h^i} \lesssim \norm{\Grad \mm_h^i} \lesssim 1$
for all $i=0,\dots,N$ (the estimate of the gradient follows from~\eqref{eq:stability} and~\eqref{eq:projection_stability}).
We infer that the sequences $\{ \hat\mm_{hk}^\pm \}$ are uniformly bounded in $L^\infty(0,T;\HH^1(\Omega))$
and arguing as in the proof of Lemma~\ref{lem:convergent_subsequences},
we can extract subsequences satisfying the convergence properties
in~\eqref{eq:conv_mhat_LinftyH1}--\eqref{eq:conv_mhat_L2}.
The fact that the limit is exactly the function $\mm \in L^\infty(0,T;\HH^1(\Omega;\sphere))$
constructed in Lemma~\ref{lem:convergent_subsequences} follows from Lemma~\ref{lem:projection_error_Lp}
(applied $p=1$), which guarantee that $\hat\mm_{hk}^{\pm} \to \mm$ in $\LL^1(\Omega_T)$,
which in turn implies that the limit functions in $\LL^2(\Omega_T)$, $L^2(0,T;\HH^1(\Omega))$ and $L^\infty(0,T;\HH^1(\Omega))$ must necessarily be the same.

To show~\eqref{eq:conv_m2_L2}--\eqref{eq:conv_mhat2_L2},
we note that for $\boldsymbol{x}, \boldsymbol{y} \in\R^3$ we have
\begin{equation*}
\boldsymbol{x}\otimes \boldsymbol{x} - \boldsymbol{y}\otimes \boldsymbol{y}
= \frac{1}{2}
[(\boldsymbol{x}+\boldsymbol{y})\otimes (\boldsymbol{x}-\boldsymbol{y})
+ (\boldsymbol{x}-\boldsymbol{y})\otimes(\boldsymbol{x}+\boldsymbol{y})].
\end{equation*}
Let $3/4 \le s < 1$.
Using the above identity and the continuous inclusion $\HH^s(\Omega) \subset \LL^4(\Omega)$ for all $s \ge 3/4$,
for arbitrary $t \in (0,T)$, we have
\begin{equation*}
\begin{split}
\norm{\mm_{hk}^\pm(t)\otimes\mm_{hk}^\pm(t) - \mm(t)\otimes\mm(t)}
& \leq
\norm[\LL^4(\Omega)]{\mm_{hk}^\pm(t) + \mm(t)} \norm[\LL^4(\Omega)]{\mm(t) - \mm_{hk}^\pm(t)} \\
& \lesssim \norm[\HH^1(\Omega)]{\mm_{hk}^\pm(t) + \mm(t)} \norm[\HH^{s}(\Omega)]{\mm(t) - \mm_{hk}^\pm(t)}.
\end{split}
\end{equation*}
It follows that
\begin{equation*}
\begin{split}
\norm[\LL^2(\Omega_T)]{\mm_{hk}^\pm\otimes\mm_{hk}^\pm - \mm\otimes\mm}
& \lesssim \norm[L^\infty(0,T;\HH^1(\Omega))]{\mm_{hk}^\pm + \mm} \norm[L^2(0,T;\HH^{s}(\Omega))]{\mm - \mm_{hk}^\pm}.
\end{split}
\end{equation*}
Convergence~\eqref{eq:conv_m2_L2} then follows
from the uniform boundedness of both $\mm_{hk}^\pm$ and $\mm$ in $L^\infty(0,T;\HH^1(\Omega))$
and the strong convergence~\eqref{eq:conv_strongL2Hs} from Lemma~\ref{lem:convergent_subsequences}.
The proof of~\eqref{eq:conv_mhat2_L2} is identical
(due to the use of the nodal projection, one can use
the H\"older inequality $\norm[L^2]{\cdot} \le \norm[L^\infty]{\cdot}\norm[L^2]{\cdot}$).
This concludes the proof.
\end{proof}

Now, we are in a position to prove Theorem~\ref{thm:convergence}(i).

\begin{proof}[Proof of Theorem~\ref{thm:convergence}(i)]
We apply Lemma~\ref{lem:convergent_subsequences}, which yields
$\uu \in L^\infty(0,T;\HH^1_D(\Omega))$ with $\uut \in L^\infty(0,T;\LL^2(\Omega))$
and
$\mm \in L^\infty(0,T;\HH^1(\Omega;\sphere))$ with $\mmt \in L^2(0,T;\LL^2(\Omega))$
as well as subsequences of $\{\uu_{hk}\}$ and $\{\mm_{hk}\}$
satisfying the desired convergence properties.
This already shows that $\uu$ and $\mm$ satisfy property~(i) of Definition~\ref{def:weak}.
Property~(iii) follows from
the available convergence results,
the continuity of the trace operator $\HH^1(\Omega_T) \to \HH^{1/2}(\Omega)$,
and assumption~\eqref{eq:convergence_initial_data} on the discrete initial data.
To conclude the proof, it remains to show that property~(ii) holds,
i.e., that $\uu$ and $\mm$ satisfy the variational formulations~\eqref{eq:weak_u} and~\eqref{eq:weak_m}, respectively. 
The result follows from the convergence properties established in Lemmas~\ref{lem:convergent_subsequences}--\ref{lem:convergent_subsequences_aux}.
We omit the details, because
\begin{itemize}
\item the proof that $\uu$ satisfies~\eqref{eq:weak_u}
is identical to the one presented in~\cite[page~1378]{bppr2013},
which is a consequence of the fact that in the displacement update~\eqref{alg1:displacement_update}
we employ the nodal projection for the magnetisation appearing on the right-hand side
(our generalised setting involving a more general expression for the magnetostrain, body forces and traction
does not pose further mathematical challenges here).
\item the proof that $\mm$ satisfies~\eqref{eq:weak_m} can be obtained combining the argument of~\cite[pages~1376--1378]{bppr2013}
(which show convergence of the method with nodal projection towards a variational formulation of the LLG equation with magnetoelastic term)
with the one of~\cite[page~1363]{hpprss2019}
(where the modifications due to the omission of the nodal projection are presented).
\end{itemize}
This concludes the proof.
\end{proof}

\subsection{Energy inequality} \label{sec:proofs_energy_inequality}

In this section, we use the compact notation $\hat\mm_h = \Pi_h\mm_h$ to denote
the nodal projection of a general magnetisation approximation $\mm_h$.

To start with,
in the following proposition,
we state a variant of Proposition~\ref{prop:energy}
for the discrete energy
\begin{equation*}
\hat\E_h[\uu_h,\mm_h]
=
\frac{1}{2}\norm{\Grad \mm_h}^2
+ \frac{1}{2} \norm[\C]{\boldvar(\uu_h) - \boldvarm(\hat\mm_h)}^2
- \inner{\ff}{\uu_h}
- \inner[\Gamma_{N}]{\gg}{\uu_h},
\end{equation*}
which is obtained from~\eqref{eqn:definition_of_energy} 
by applying the nodal projection to the discrete magnetisation
appearing in the elastic energy.
We omit the proof since it is very similar to the one of Proposition~\ref{prop:energy}.

\begin{proposition} \label{prop:energy_modified}
For every integer $0 \le i \le N-1$,
the iterates of Algorithm~\ref{algorithm} satisfy the discrete energy law
\begin{equation} \label{eq:modified_discrete_energy_law}
\hat\E_h[\uu_h^{i+1},\mm_h^{i+1}] + \frac{1}{2} \norm{\dt\uu_h^{i+1}}^2
-
\hat\E_h[\uu_h^i,\mm_h^i] - \frac{1}{2} \norm{\dt\uu_h^i}^2
=
- \alpha k \norm[h]{\vv_h^i}^2
- \hat D_{h,k}^i
- \hat E_{h,k}^i,
\end{equation}
where $\hat D_{h,k}^i$
and $\hat E_{h,k}^i$ are given by
\begin{multline*}
\hat D_{h,k}^i
:= k^2 (\theta-1/2) \norm{\Grad\vv_h^i}^2
+ \frac{1}{2} \norm{\dt\uu_h^{i+1} - \dt\uu_h^i}^2 \\
+ \frac{1}{2} \norm[\C]{[\boldvar(\uu_{h}^{i+1})-\boldvarm(\hat\mm_{h}^{i+1})]
- [\boldvar(\uu_{h}^i)-\boldvarm(\hat\mm_{h}^i)]}^2 \ge 0
\end{multline*}
and
\begin{equation*}
\begin{split}
\hat E_{h,k}^i & := \sum_{\ell=1}^5 \hat E_{h,k,\ell}^i \\
& := \inner{\boldsig(\uu_{h}^{i+1},\hat\mm_{h}^{i+1}) - \boldsig(\uu_{h}^i,\hat\mm_{h}^i)}{\boldvarm(\hat\mm_{h}^{i+1}) - \boldvarm(\hat\mm_{h}^i)} \\
& \quad 
+ \inner{\boldsig(\uu_{h}^i,\hat\mm_{h}^i)}{\boldvarm(\hat\mm_{h}^{i+1}) - \boldvarm(\mm_{h}^{i+1})}
+ \inner{\boldsig(\uu_{h}^i,\hat\mm_{h}^i)}{\boldvarm(\mm_{h}^i) - \boldvarm(\hat\mm_{h}^i)} \\
& \quad
+ 2k \inner{\boldsig(\uu_{h}^i,\hat\mm_{h}^i)}{\Z:[(\hat\mm_{h}^{i} - \mm_{h}^{i}) \otimes \vv_{h}^{i}]}
+ k^2 \inner{\boldsig(\uu_{h}^i,\hat\mm_{h}^i)}{\boldvarm(\vv_{h}^i)}.
\end{split}
\end{equation*}
respectively.
\end{proposition}

Now, we are in a position to prove Theorem~\ref{thm:convergence}(ii).

\begin{proof}[Proof of Theorem~\ref{thm:convergence}(ii)]
Let $t' \in (0,T)$.
Let $1 \le j \le N$ such that $t' \in (t_{j-1},t_j)$.
Summing~\eqref{eq:modified_discrete_energy_law} for $i=0,\dots,j-1$ yields
\begin{equation*}
\hat\E_h[\uu_h^j,\mm_h^j] + \frac{1}{2} \norm{\dt\uu_h^j}^2
-
\hat\E_h[\uu_h^0,\mm_h^0] - \frac{1}{2} \norm{\dt\uu_h^0}^2
+ \alpha k \sum_{i=0}^{j-1}\norm[h]{\vv_h^i}^2
+ \sum_{i=0}^{j-1} \hat D_{h,k}^i
=
- \sum_{i=0}^{j-1} \hat E_{h,k}^i.
\end{equation*}
Using the Cauchy--Schwarz inequality, the weighted Young inequality,
and Lemma~\ref{lem:lipschitz_magnetostrain},
we obtain the estimate
\begin{equation*}
\begin{split}
\lvert E_{h,k,1}^i \rvert
& =
\lvert
\inner{\boldsig(\uu_{h}^{i+1},\hat\mm_{h}^{i+1}) - \boldsig(\uu_{h}^i,\hat\mm_{h}^i)}{\boldvarm(\hat\mm_{h}^{i+1}) - \boldvarm(\hat\mm_{h}^i)}
\rvert \\
& \le
\norm[\C]{[\boldvar(\uu_{h}^{i+1}) - \boldvarm(\hat\mm_{h}^{i+1})] - [\boldvar(\uu_{h}^i) - \boldvarm(\hat\mm_{h}^i)]}
\norm[\C]{\boldvarm(\hat\mm_{h}^{i+1}) - \boldvarm(\hat\mm_{h}^i)} \\
& \le
\frac{1}{4} \norm[\C]{[\boldvar(\uu_{h}^{i+1}) - \boldvarm(\hat\mm_{h}^{i+1})] - [\boldvar(\uu_{h}^i) - \boldvarm(\hat\mm_{h}^i)]}^2 \\
& \quad + C_{\mathrm{m}}^2 \norm[\LL^{\infty}(\Omega)]{\C} k^2
\norm{\vv_{h}^i}^2.
\end{split}
\end{equation*}
We now estimate $E_{h,k,4}^i$ (assuming $i \ge 1$, because $E_{h,k,4}^0=0$ as $\hat\mm_{h}^0 = \mm_{h}^0$ by assumption).
Using the Cauchy--Schwarz inequality,
the H\"older inequality
(for $p=2/(1-2\eps)$ and $p' = 2p/(p-2)$ with $0 < \eps \ll 1/2$ arbitrary),
Lemma~\ref{lem:projection_error_Lp},
and classical inverse estimates
(see, e.g.\ \cite[Lemma~3.5]{bartels2015}),
we obtain
\begin{equation*}
\begin{split}
\lvert E_{h,k,4}^i \rvert
& =
2k \lvert \inner{\boldsig(\uu_{h}^i,\hat\mm_{h}^i)}{\Z:[(\hat\mm_{h}^{i} - \mm_{h}^{i}) \otimes \vv_{h}^{i}]} \rvert \\
& \le
2 \norm[\LL^\infty(\Omega)]{\Z} k \norm{\boldsig(\uu_{h}^i,\hat\mm_{h}^i)}
\norm[\LL^p(\Omega)]{\hat\mm_{h}^{i} - \mm_{h}^{i}}\norm[\LL^{p'}(\Omega)]{\vv_{h}^{i}} \\
& \lesssim
k \norm{\boldsig(\uu_{h}^i,\hat\mm_{h}^i)} k^{(p+1)/p}
\left( \sum_{\ell=0}^{i-1} \norm[\LL^{2p}(\Omega)]{\vv_{h}^{\ell}}^2  \right)
\norm[\LL^{p'}(\Omega)]{\vv_{h}^{i}} \\
& \lesssim
k^{2+1/p} \norm{\boldsig(\uu_{h}^i,\hat\mm_{h}^i)}
h_{\min}^{3(1-p)/p} \left( \sum_{\ell=0}^{i-1} \norm{\vv_{h}^{\ell}}^2  \right)
h_{\min}^{3(2-p')/(2p')} \norm{\vv_{h}^{i}} \\
& \lesssim
h_{\min}^{-3 } k^{5/2-\eps} \norm{\boldsig(\uu_{h}^i,\hat\mm_{h}^i)}
\left( \sum_{\ell=0}^{i-1} \norm{\vv_{h}^{\ell}}^2  \right)
\norm{\vv_{h}^{i}}.
\end{split}
\end{equation*}
Similarly, we obtain
\begin{equation*}
\lvert E_{h,k,5}^i \rvert
=
\lvert k^2 \inner{\boldsig(\uu_{h}^i,\hat\mm_{h}^i)}{\boldvarm(\vv_{h}^i)} \rvert
\lesssim
h_{\min}^{-3/2} k^2 \norm{\boldsig(\uu_{h}^i,\hat\mm_{h}^i)} \norm{\vv_{h}^{i}}^2.
\end{equation*}
Moreover, noting $\mm_{h}^{0} = \hat{\mm}_{h}^{0}$ we have that
\begin{equation*}
\begin{split}
& \sum_{i=0}^{j-1} \left(\hat E_{h,k,2}^i + \hat E_{h,k,3}^i \right) \\
& \quad =
\sum_{i=0}^{j-1} \inner{\boldsig(\uu_{h}^i,\hat\mm_{h}^i)}{\boldvarm(\hat\mm_{h}^{i+1}) - \boldvarm(\mm_{h}^{i+1})}
+ \sum_{i=0}^{j-1} \inner{\boldsig(\uu_{h}^i,\hat\mm_{h}^i)}{\boldvarm(\mm_{h}^i) - \boldvarm(\hat\mm_{h}^i)} \\
& \quad =
\inner{\boldsig(\uu_{h}^0,\hat\mm_{h}^0)}{\boldvarm(\mm_{h}^0) - \boldvarm(\hat\mm_{h}^0)}
- \inner{\boldsig(\uu_{h}^{j-1},\hat\mm_{h}^{j-1})}{\boldvarm(\mm_{h}^j) - \boldvarm(\hat\mm_{h}^j)} \\
& \qquad
+ \sum_{i=0}^{j-2} \inner{\boldsig(\uu_{h}^{i+1},\hat\mm_{h}^{i+1}) - \boldsig(\uu_{h}^i,\hat\mm_{h}^i)}{\boldvarm(\mm_{h}^i) - \boldvarm(\hat\mm_{h}^i)}.\\
& \quad =
- \inner{\boldsig(\uu_{h}^{j-1},\hat\mm_{h}^{j-1})}{\boldvarm(\mm_{h}^j) - \boldvarm(\hat\mm_{h}^j)} \\
& \qquad
+ \sum_{i=0}^{j-2} \inner{\boldsig(\uu_{h}^{i+1},\hat\mm_{h}^{i+1}) - \boldsig(\uu_{h}^i,\hat\mm_{h}^i)}{\boldvarm(\mm_{h}^i) - \boldvarm(\hat\mm_{h}^i)}
\end{split}
\end{equation*}
Using inverse estimates, Lemma~\ref{lem:lipschitz_magnetostrain} and Lemma~\ref{lem:projection_error_Lp},
we obtain the estimate
\begin{equation*}
\begin{split}
& \lvert \inner{\boldsig(\uu_{h}^{i+1},\hat\mm_{h}^{i+1}) - \boldsig(\uu_{h}^i,\hat\mm_{h}^i)}{\boldvarm(\mm_{h}^i) - \boldvarm(\hat\mm_{h}^i)} \rvert \\
& \quad \lesssim
\left(\norm{\boldvar(\uu_h^{i+1}) - \boldvar(\uu_h^i)} + \norm{\boldvarm(\hat\mm_h^{i+1}) - \boldvarm(\hat\mm_h^i)} \right)
\norm{\boldvarm(\mm_{h}^i) - \boldvarm(\hat\mm_{h}^i)} \\
& \quad \lesssim
\left(h_{\min}^{-1} k \norm{\dt\uu_h^{i+1}} + k^2 \norm{\vv_h^i}^2 \right)
\norm[\HH^1(\Omega)]{\mm_h^i + \hat\mm_h^i} k^{4/3} h_{\min}^{-2} \sum_{\ell=0}^{i-1}  \norm{\vv_h^\ell}^2.
\end{split}
\end{equation*}
Altogether,
omitting all non-negative dissipative terms
and using the stability from Proposition~\ref{prop:stability}, we thus obtain
\begin{equation*}
\begin{split}
& \hat\E_h[\uu_h^j,\mm_h^j] + \frac{1}{2} \norm{\dt\uu_h^j}^2
-
\hat\E_h[\uu_h^0,\mm_h^0] - \frac{1}{2} \norm{\dt\uu_h^0}^2
+ \alpha k \sum_{i=0}^{j-1}\norm[h]{\vv_h^i}^2 \\
& \qquad\quad - \inner{\boldsig(\uu_{h}^{j-1},\hat\mm_{h}^{j-1})}{\boldvarm(\mm_{h}^j) - \boldvarm(\hat\mm_{h}^j)} \\
& \quad \lesssim
\sum_{i=0}^{j-1}
\Big(
k^2 \norm{\vv_h^i}^2
+ h_{\min}^{-3} k^{4/3} \norm{\dt\uu_h^{i+1}}\norm[\HH^1(\Omega)]{\mm_h^i + \hat\mm_h^i}
+ h_{\min}^{-2} k^{7/3} \norm{\vv_h^i}^2 \\
& \qquad\quad
+ h_{\min}^{-3} k^{3/2-\eps} \norm{\boldsig(\uu_{h}^i,\hat\mm_{h}^i)}\norm{\vv_{h}^{i}}
+ h_{\min}^{-3/2} k^2 \norm{\vv_{h}^{i}}^2
\Big) \\
& \quad \lesssim
k
+ h_{\min}^{-3} k^{1/3}
+ h_{\min}^{-2} k^{4/3}
+ h_{\min}^{-3} k^{1/2-\eps}
+ h_{\min}^{-3/2} k.
\end{split}
\end{equation*}
Using~\eqref{eq:h-scalar-product_equivalence},
rewriting the above using the time reconstructions~\eqref{eq:reconstructions}
and integrating in time over an arbitrary measurable set $\mathfrak{T} \subset [0,T]$,
we obtain
\begin{equation*}
\begin{split}
& \int_{\mathfrak{T}} \left(
\E[\uu_{hk}^+(t'),\mm_{hk}^+(t')]
+ \frac{1}{2} \norm{\dot\uu_{hk}^+(t')}^2
-
\hat\E_h[\uu_{hk}^-(0),\mm_{hk}^-(0)] - \frac{1}{2} \norm{\dot\uu_{hk}^-(0)}^2
\right) \mathrm{d}t' \\
& \qquad
+ \int_{\mathfrak{T}} \left( \alpha \int_0^{t'} \norm{\vv_{hk}^-(t)}^2 \mathrm{d}t \right) \mathrm{d}t'
+ \int_{\mathfrak{T}} \left( \alpha \int_{t'}^{t_j} \norm{\vv_{hk}^-(t)}^2 \mathrm{d}t \right) \mathrm{d}t' \\
& \qquad - \int_{\mathfrak{T}} \inner{\boldsig(\uu_{hk}^-(t'),\hat\mm_{hk}^-(t'))}{\boldvarm(\mm_{hk}^+(t')) - \boldvarm(\hat\mm_{hk}^+(t'))} \mathrm{d}t' \\
& \quad \lesssim
k
+ h_{\min}^{-3} k^{1/3}
+ h_{\min}^{-2} k^{4/3}
+ h_{\min}^{-3} k^{1/2-\eps}
+ h_{\min}^{-3/2} k.
\end{split}
\end{equation*}
We now consider the limit of this inequality as $h,k \to 0$.
The assumed CFL condition $k=o(h^9)$ implies that the right-hand
side converges to $0$ in the limit as $h,k \to 0$.
The last two terms on the left-hand side converge to $0$:
the first one by no concentration of Lebesgue functions,
the other thanks to the available convergence results
(cf. the convergences guaranteed by Lemmas~\ref{lem:convergent_subsequences}--\ref{lem:convergent_subsequences_aux}).
Weak lower semicontinuity guarantees
\begin{multline*}
\int_{\mathfrak{T}} \left(
\E[\uu(t'),\mm(t')]
+ \frac{1}{2} \norm{\partial_t\uu(t')}^2
+ \alpha \int_0^{t'} \norm{\mmt(t)}^2
\mathrm{d}t\right) \mathrm{d}t' \\
\le \liminf_{h,k \to 0}
\int_{\mathfrak{T}} \left(
\hat\E_h[\uu_{hk}^+(t'),\mm_{hk}^+(t')]
+ \frac{1}{2} \norm{\dot\uu_{hk}^+(t')}^2
+ \alpha \int_0^{t'} \norm{\vv_{hk}^-(t)}^2
\mathrm{d}t\right) \mathrm{d}t'.
\end{multline*}
Assumption~\eqref{eq:convergence_initial_data} yields
\begin{equation*}
\lim_{h,k \to 0} \left( \hat\E_h[\uu_{hk}^-(0),\mm_{hk}^-(0)] + \frac{1}{2} \norm{\dot\uu_{hk}^-(0)}^2 \right)
=
\E[\uu^0,\mm^0] + \frac{1}{2} \norm{\dot\uu^0}^2.
\end{equation*}
Since $\mathfrak{T} \subset [0,T]$ was arbitrary, this shows that the energy inequality~\eqref{eqn:EnergyLawDefinition}
holds a.e.\ in $(0,T)$ and concludes the proof.
\end{proof}

\section{Code Availability Statement}
\comment{The code used to generate the findings of this study is openly available in Zenodo at \url{https://doi.org/10.5281/zenodo.14641594}.}

\section*{Acknowledgements}

MR is a member of the `Gruppo Nazionale per il Calcolo Scientifico (GNCS)'
of the Italian `Istituto Nazionale di Alta Matematica (INdAM)'
\comment{and
was partially supported by GNCS (research project GNCS 2024 on \emph{Advanced numerical methods 
for nonlinear problems in materials science} -- CUP E53C23001670001)}.
The authors thank Martin Kru\v{z}\'ik
(Institute of Information Theory and Automation, Czech Academy of Sciences)
for several interesting discussions on magnetoelasticity.
The support of the Royal Society (grant IES{\textbackslash}R2{\textbackslash}222118)
and the Czech Ministry of Education, Youth and Sports
(M\v{S}MT \v{C}R project 8J22AT017)
is thankfully acknowledged.

\bibliographystyle{habbrv}
\bibliography{ref}

\appendix


\section{Linear algebra definitions and identities} 
\label{sec:linear_algebra}

In this section, for the convenience of the reader,
we collect some definitions and vector/matrix/tensor identities from linear algebra that are used throughout the work.

\begin{definition}
Let $\A \in \R^{3^4}$ be a fourth-order tensor (4-tensor) with components $\A_{ij\ell m}$, where $i,j,\ell,m = 1,2,3$.
We say that
\begin{enumerate}
\item $\A$ is minorly symmetric if $\A_{ij\ell m} = \A_{ji\ell m}=\A_{ij m \ell}$,
\item $\A$ is majorly symmetric if $\A_{ij\ell m} = \A_{\ell mji}$,
\item and $\A$ is (fully) symmetric if the aforementioned two conditions hold together.
\end{enumerate}
The transpose of $\A$ is the 4-tensor $\A^\top \in \R^{3^4}$ given by $(\A^\top)_{ij\ell m} = \A_{\ell mji}$.
In particular, $\A$ is majorly symmetric if $\A^\top = \A$.
\end{definition}

\begin{remark}
In three dimensions, 4-tensors have $3^4=81$ components.
Minorly symmetric 4-tensors have 36 independent components,
majorly symmetric 4-tensors have 45 independent components,
and fully symmetric 4-tensors have 21 independent components.
Throughout this work the stiffness tensor $\C$ is assumed to be fully symmetric,
whereas the magnetostriction tensor $\Z$ is assumed to be only minorly symmetric.
In the numerical experiments of Section~\ref{sec:numerics},
we consider the isotropic case,
in which $\C$ and $\Z$ have only two (the so-called Lam\'e constants)
and one (the so-called saturation magnetostriction) independent components, respectively.
\end{remark}

In the following definition, we recall some operations between tensors.

\begin{definition}
Let $\A,\mathbb{B} \in \R^{3^4}$ be 4-tensors,
let $\boldsymbol\nu,\boldsymbol\mu \in \R^{3 \times 3}$ be 2-tensors (matrices),
and let $\mm,\ww \in \R^3$ be 1-tensors (vectors).
    \begin{itemize}
        \item We denote the double contraction between $\A$ and $\mathbb{B}$
        as the 4-tensor $\A:\mathbb{B} \in \R^{3^4}$ given by
        \begin{equation*}
            (\A:\mathbb{B})_{ij\ell m} = \sum_{p,q}\A_{ijpq}\mathbb{B}_{pq\ell m}.
        \end{equation*}
        \item We denote the double contraction between $\A$ and $\boldsymbol\nu$
        as the 2-tensor $\A : \boldsymbol\nu \in \R^{3 \times 3}$ given by
        \begin{equation*}
            (\A : \boldsymbol\nu)_{ij} = \sum_{\ell,m}\A_{ij\ell m}\nu_{\ell m}.
        \end{equation*}
        \item We denote the Frobenius product of $\boldsymbol\mu$ and $\boldsymbol\nu$
        as the scalar $\boldsymbol\mu:\boldsymbol\nu \in \R$ given by
        \begin{equation*}
             \boldsymbol\mu : \boldsymbol\nu = \sum_{i,j}\mu_{ij} \nu_{ij}.
        \end{equation*}
        \item We denote the tensor product of $\mm$ and $\ww$
        as the 2-tensor $\mm\otimes\ww \in \R^{3\times 3}$ given by
        \begin{equation*}
            (\mm\otimes\ww)_{ij} = m_{i}w_{j}.
        \end{equation*}
    \end{itemize}
\end{definition}

The following result is useful for manipulation of the magnetostrain terms.

\begin{lemma}\label{lem:tensor}
Let $\Z \in \R^{3^4}$ be a minorly symmetric 4-tensor,
let $\boldsig \in \R^{3 \times 3}$ be a symmetric 2-tensor,
and let $\mm,\ww\in\R^3$ be two 1-tensors.
We have the identity
\begin{equation} \label{eq:tensor_identity}
[(\Z^\top:\boldsig)\ww]\cdot\mm=[\Z^\top:\boldsig)\mm]\cdot\ww = \boldsig:[\Z:(\mm\otimes\ww)].
\end{equation}
\end{lemma}

\begin{proof}
We have by the minor symmetry of $\Z$ that
\begin{align*}
    (\mathbb{Z}^\top:\boldsig) \mm \cdot \ww
    &= \sum_{i,j,\ell,m}\mathbb{Z}_{ijk\ell}^\top\sigma_{\ell m}m_{j}w_{i}\\
    &= \sum_{i,j,\ell,m}\mathbb{Z}_{jik\ell}^\top\sigma_{\ell m}m_{i}w_{j}\quad\textnormal{by relabelling}\\
    &= \sum_{i,j,\ell,m} \mathbb{Z}_{ijk\ell}^\top\sigma_{\ell m}w_{j}m_{i}\quad\textnormal{via minor symmetry}\\
    &= (\mathbb{Z}^\top:\boldsig) \ww \cdot \mm.
\end{align*}
Furthermore, we have that
\begin{align*}
    (\mathbb{Z}^\top:\boldsig) \mm \cdot \ww
    &= \sum_{i,j,\ell,m} \mathbb{Z}_{ij\ell m}^\top\sigma_{\ell m}m_{j}w_{i}
    = \sum_{i,j,\ell,m} \sigma_{\ell m}\Z_{ij\ell m}^\top m_{j}w_{i}\\
    &= \sum_{i,j,\ell,m} \sigma_{\ell m}\Z_{\ell mji}m_{j}w_{i}\quad\textnormal{via minor symmetry}\\
    &= \sum_{\ell,m} \sigma_{\ell m}[\Z:(\mm\otimes\ww)]_{\ell m}
    = \boldsig:[\Z:(\mm\otimes\ww)].
\end{align*}
This show both identities in~\eqref{eq:tensor_identity}.
\end{proof}

The following identity is useful to show the stability of numerical schemes.

\begin{lemma} \label{lem:abel}
Let $\{\nu_i\}$ be a sequence in an inner product space
with inner product $\inner{\cdot}{\cdot}$
and associated norm $\norm{\cdot}$.
We have the identity
\begin{equation} \label{eq:abel}
    \inner{\nu_{i+1} - \nu_{i}}{\nu_{i+1}}
    = \frac{1}{2}\norm{\nu_{i+1}}^2
    - \frac{1}{2}\norm{\nu_{i}}^2
    + \frac{1}{2}\norm{\nu_{i+1} - \nu_{i}}^2.
\end{equation}
\end{lemma}

\section{Nondimensionalisation}
\label{sec:physics}

Let $\Omega \subset \R^3$ denote the volume occupied by a ferromagnetic body
(with the spatial variable $x \in \Omega$ measured in \si{meter}).
Consider the magnetisation $\boldsymbol{M}$ (measured in \si{A/m}),
which satisfies the length constraint $\abs{\boldsymbol{M}} = M_s$,
where the saturation magnetisation $M_s>0$ is also measured in \si{A/m},
and the displacement $\UU$ (measured in \si{m}).
We denote by $\boldvar(\UU)$ the total strain given by
\begin{equation*}
    \boldvar(\UU) = (\Grad \UU + \Grad \UU^\top)/2
\end{equation*}
and by $\boldvar_{\mathrm{m}}(\boldsymbol{M})$
the magnetostrain given by
\begin{equation*}
    \boldvarm(\boldsymbol{M})=\Z:(\boldsymbol{M}\otimes \boldsymbol{M}/M_s^2),
\end{equation*}
where $\Z$ is a dimensionless fourth-order tensor.

The total energy of the system (measured in \si{\joule}) is given by
\begin{multline*}
\mathcal{E}[\UU,\boldsymbol{M}]
= \frac{A}{M_s^2}\int_{\Omega}|\Grad \boldsymbol{M}|^2 
- \mu_0 \int_{\Omega} \HH_{\mathrm{ext}} \cdot \boldsymbol{M} \\
+ \frac{1}{2}\int_{\Omega}[\boldvar(\UU) -\boldvarm(\boldsymbol{M})]:\{\CC:[\boldvar(\UU) -\boldvarm(\boldsymbol{M})]\}
- \int_{\Omega}\FF\cdot\UU
- \int_{\Gamma_N}\boldsymbol{G}\cdot\UU,
\end{multline*}
where
$A$ is the exchange constant (measured in \si{J.m^{-1}}),
$\mu_0$ is the permeability of free space (measured in \si{N.A^{-2}}),
$\HH_{\mathrm{ext}}$ is an applied external field (measured in \si{\ampere\per\meter}),
$\CC$ is the fourth-order stiffness tensor (measured in \si{N.m^{-2}}),
$\FF$ is a body force (measured in \si{N.m^{-3}}),
and
$\boldsymbol{G}$ is a surface force (measured in \si{N.m^{-2}}).

The dynamics of $\boldsymbol{M}$ is described by the LLG equation:
\begin{equation*}
\partial_t \boldsymbol{M}
= -\gamma \mu_{0} \, \boldsymbol{M} \times \HH_{\mathrm{eff}}[\UU,\boldsymbol{M}]
+ \frac{\alpha}{M_s}\boldsymbol{M} \times \partial_t \boldsymbol{M},
\end{equation*}
where $\gamma$ is the gyromagnetic ratio (measured in \si{rad.s^{-1}T^{-1}}),
$\alpha>0$ is the dimensionless Gilbert damping parameter,
and the effective field $\HH_{\mathrm{eff}}$ (measured in \si{A/m}) reads as
\begin{equation*}
\HH_{\mathrm{eff}}[\UU,\boldsymbol{M}]
= - \frac{1}{\mu_0} \frac{\delta\mathcal{E}[\UU,\boldsymbol{M}]}{\delta\boldsymbol{M}}
= \frac{2A}{\mu_0 M_s^2}\boldsymbol\Delta \boldsymbol{M} + \HH_{\mathrm{ext}} + \frac{2}{\mu_0 M_s^2}[\Z^\top :\boldsymbol\Sigma(\UU,\boldsymbol{M})]\boldsymbol{M},
\end{equation*}
where $\boldsymbol{\Sigma}(\UU,\boldsymbol{M}) = \CC:[\boldvar(\UU) -\boldvarm(\boldsymbol{M})]$ is the stress (measured in \si{N.m^{-2}}).
The LLG equation is coupled with the conservation of momentum equation satisfied by the displacement:
\begin{equation*}
    \rho \, \partial_{tt} \UU = \Grad \cdot \boldsymbol{\Sigma}(\UU,\boldsymbol{M}) + \FF,
\end{equation*}
where $\rho$ is the mass density (measured in \si{kg.m^{-3}}).

Let $\mm = \boldsymbol{M} / M_s$ denote the normalised magnetisation.
We define the exchange length $\ellex^2 = 2A / \mu_0 M_s^2$ (measured in \si{m})
and use it to rescale the spatial variable and the displacement according to $x' = x/\ellex$ and $\uu = \UU / \ellex$, respectively.
Additionally we introduce the dimensionless domain $\Omega'=\Omega/\ellex$,
the dimensionless time $t' = \gamma \mu_0 M_s t$,
the dimensionless coupling parameter $\kappa = \rho \ellex^2\gamma^2 \mu_0$,
as well as
the dimensionless differential operators
$\Grad=\Grad'/\ellex$ and $\boldsymbol\Delta=\boldsymbol\Delta'/\ellex^2$.
Further we define the dimensionless energy as
\begin{equation*}
\begin{split}
\E'[\uu,\mm]
& = \frac{\mathcal{E}[\comment{\ellex\uu},M_s \mm]}{\mu_0 M_s^2\ellex^3} \\
& = \frac{1}{2}\int_{\Omega'}|\Grad' \mm|^2
- \int_{\Omega'} \hh_{\mathrm{ext}} \cdot\mm \\
& \qquad
+ \frac{\kappa}{2}\int_{\Omega'}[\boldvar(\uu) - \boldvarm'(\mm)]: \{\C:[\boldvar(\uu) -\boldvarm'(\mm)]
- \kappa\int_{\Omega'}\ff\cdot\uu 
- \kappa\int_{\Gamma_N'}\boldsymbol{g}\cdot\uu,
\end{split}
\end{equation*}
where
$\hh_{\mathrm{ext}} = \HH_{\mathrm{ext}} / M_s$,
$\boldvarm(\mm) = \Z:(\mm\otimes\mm)$,
$\boldsymbol\Sigma = \kappa \mu_0 M_s^2\boldsig$,
$\boldsymbol{G} = \kappa \mu_0 M_s^2\boldsymbol{g}$,
$\FF = \kappa(\mu_0 M_s^2/\ellex)\ff$,
and $\CC = \kappa \mu_0 M_s^2\C$
(all dimensionless).
The dimensionless effective field, defined by
$\hh_{\mathrm{eff}}[\uu,\mm] = \HH_{\mathrm{eff}}[\ellex\uu,M_s\mm]/M_s$,
satisfies
the relation
\begin{equation*}
    \hh_{\mathrm{eff}}[\uu,\mm]
    = -\frac{\delta \mathcal{E}'[\uu,\mm]}{\delta \mm}
    = \boldsymbol\Delta'\mm + 2\kappa[\Z^\top :\boldsig(\uu,\mm)]\mm + \hh_{\mathrm{ext}}.
\end{equation*}
With all these definitions,
we retrieve the coupled system
\begin{align*}
    \partial_{t't'} \uu &= \comment{\Grad'} \cdot \boldsig(\uu,\mm) + \ff,\\
    \partial_{t'} \mm &= -\mm \times \hh_{\textnormal{eff}}[\uu,\mm] + \alpha \, \mm \times \partial_{t'} \mm.
\end{align*}
Altogether, we thus obtain the dimensionless model problem discussed throughout this work.
Note that, to simplify the notation, in Sections~\ref{sec:model}--\ref{sec:proofs}
we omit all `primes' from the dimensionless quantities,
we assume $\kappa=1$,
and we neglect the applied external field (unless otherwise mentioned).

\end{document}